\documentclass[a4paper,10pt,oneside]{article}
\usepackage{amsfonts,amsmath,amssymb,amsthm,verbatim,pstricks-add,mathtools}
\usepackage{fullpage}
\usepackage{tikz}
\usepackage{authblk}
\usepackage{url}

\usepackage[normalem]{ulem}

\newcommand{\Z}{\mathbb{Z}}
\newcommand{\D}{\Delta}
\newcommand{\hD}{\hat{\Delta}}
\newcommand{\hc}{\hat{c}}
\newcommand{\pres}[2]{\langle {#1}\ |\ {#2} \rangle}
\newcommand{\gpres}[1]{\langle {#1} \rangle}
\newtheorem{theorem}{Theorem}
\newtheorem{lemma}[theorem]{Lemma}

\newtheorem{maintheorem}{Theorem}

\newtheorem{maincorollary}[maintheorem]{Corollary}
\newtheorem{claim}{Claim}
\newtheorem*{lemma*}{Main Lemma}

\theoremstyle{definition}
\newtheorem{remark}[theorem]{Remark}
\newtheorem*{remark*}{Remark}
\newtheorem*{notation*}{Notation}

\newcommand{\textoverline}[1]{\ensuremath{\overline{\text{#1}}}}

\usepackage{textcomp}

\newcommand{\encircle}[1]{%
  \tikz[baseline=(X.base)] 
    \node (X) [draw, shape=circle, inner sep=0] {\strut $#1$};}

\newcommand{\LL}{\mathbf{L}}
\newcommand{\DD}{\mathbf{D}}
\renewcommand{\SS}{\mathbf{S}}

\title{All hyperbolic cyclically presented groups with positive length three relators}
\author[1]{Ihechukwu Chinyere}
\author[2]{Martin Edjvet}
\author[3]{Gerald Williams}
\affil[1]{Department of Mathematics and Applied Mathematics, University of Pretoria,  Hatfield 0028, South Africa, \texttt{i.chinyere@up.ac.za}}
\affil[2]{School of Mathematical Sciences, Mathematical Sciences Building, University Park, Nottingham, NG7 2RD, United Kingdom, \texttt{martin.edjvet@nottingham.ac.uk}}
\affil[3]{School of Mathematics, Statistics and Actuarial Science, University of Essex, Colchester, Essex CO4 3SQ, United Kingdom (corresponding author), \texttt{gerald.williams@essex.ac.uk}}
\begin{document}
\maketitle

\begin{abstract}
We consider the cyclically presented groups defined by cyclic presentations with $2m$ generators $x_i$ whose relators are the $2m$ positive length three relators $x_ix_{i+1}x_{i+m-1}$. We show that they are hyperbolic if and only if $m\in \{1,2,3,6,9\}$. This completes the classification of the hyperbolic cyclically presented groups with positive length three relators.
\end{abstract}

\bigskip

\noindent {\bf  MSC:} 20F05; 20F06; 20F67.

\noindent {\bf Keywords:} hyperbolic group; cyclically presented group; curvature.

\section{Introduction}

We prove that the groups
\[\Gamma_{2m}=\pres{x_0,\ldots, x_{2m-1}}{x_ix_{i+1}x_{i+m-1}\ (0\leq i<2m)}\]
(subscripts mod $2m$) are hyperbolic if and only if $m\in \{1,2,3,6,9\}$. These groups form a family of \em cyclically presented groups, \em that is, of groups 
\[G_n(w)=\pres{x_0,\ldots, x_{n-1}}{w(x_i,x_{i+1},\ldots , x_{i+n-1})\ (0\leq i<n)}\]
(subscripts mod $n$), where $w(x_0,\ldots ,x_{n-1})$ is some word in generators $x_0,\ldots ,x_{n-1}$. In particular, the groups $\Gamma_{2m}=G_{2m}(x_0x_1x_{m-1})$ form a subfamily of the class of cyclically presented groups with positive length three relators, considered in \cite{BogleyShift,CRS05,ChinyereWilliamsT6,EdjvetWilliams,MohamedWilliams}. The hyperbolicity status of the remaining groups in this class is already known by the results of \cite{ChinyereWilliamsT6,EdjvetWilliams} and so the present article completes the classification of the hyperbolic cyclically presented groups with positive length three relators. This complements the main result of \cite{ChinyereWilliamsT5} which (except for two groups) classifies the hyperbolic cyclically presented groups with non-positive, non-negative length three relators (the so-called \em groups of Fibonacci type\em) and continues the programme of research of classifying the hyperbolic groups within certain classes of cyclically presented groups \cite{Chalk, ChalkEdjvet, ChalkEdjvetJuhasz, ChinyereWilliamsT5,ChinyereWilliamsT6,  HowieWilliamsTadpole}. It is readily verified that cyclically presented groups with length two relators are free products of either $\Z_2$ or $\Z$, and hence are hyperbolic. The abelianization $\Gamma_{2m}^{\mathrm{ab}}$, when finite, is calculated in \cite[Theorem 4.10]{NoferiniWilliams3}.

We prove the following result:
\begin{maintheorem}\label{thm:summarytheorem}
For $m\geq 1$ the group $\Gamma_{2m}$ is hyperbolic if and only if $m=1,2,3,6$ or $9$, in which case $\Gamma_{2m}$ is isomorphic to $\Z_3,\Z_{15},\Z*\Z$, $\Z_5*\Z*\Z$, $\Z_{19}*\Z*\Z$, respectively.
\end{maintheorem}

The essence of the proof is to show that, for $m\neq 1,2,3,6,9$, the group $\Gamma_{2m}$ is not hyperbolic, which we do by showing that a finite extension of $\Gamma_{2m}$ contains a free abelian subgroup of rank 2. In corollaries we consider the cyclically presented groups with positive length three  relators; that is, the groups $G_n(x_0x_kx_l)$, and for these we recall the following system of congruences introduced in \cite{EdjvetWilliams}:
%
%
\begin{itemize}
    \item[(A)] $n\equiv 0 \bmod 3$ and $k+l\equiv 0 \bmod 3$;
    \item[(B)] $k+l\equiv 0 \bmod n$ or $2l-k\equiv 0 \bmod n$ or $2k-l\equiv 0 \bmod n$;
    \item[(C)] $3l\equiv 0 \bmod n$ or $3k\equiv 0 \bmod n$ or $3(l-k)\equiv 0 \bmod n$;
    \item[(D)] $2(k+l)\equiv 0 \bmod n$ or $2(2l-k)\equiv 0 \bmod n$ or $2(2k-l)\equiv 0 \bmod n$.
\end{itemize}
We write $A=T$ (True) if condition (A) holds and $A=F$ (False) otherwise, and similarly for conditions (B),(C),(D). Note that if $A=F$ then $(B,C)\neq (T,T)$.

In Corollary \ref{cor:hyperbolicityoflength3positive} we classify the elementary and non-elementary hyperbolic groups $G_n(x_0x_kx_l)$. A free product $H*K$ is hyperbolic if and only if $H$ and $K$ are hyperbolic (see \cite[Theorem H]{BaumslagGerstenShapiroShort}); and a hyperbolic group is elementary hyperbolic if and only if it is virtually cyclic. As explained in \cite{EdjvetWilliams} the group $G_n(x_0x_kx_l)$ is isomorphic to the free product of $d=\gcd(n,k,l)$ copies of $G_{n/d}(x_0x_{k/d}x_{l/d})$, and contains a non-abelian free subgroup if $d>1$, in which case it is therefore non-elementary hyperbolic when hyperbolic, so we may assume $d=1$. If, in addition, $k=l$ then $G_n(x_0x_kx_l)\cong \Z_{2^n-(-1)^n}$ so we may also assume $k\neq l$.

\begin{maincorollary}\label{cor:hyperbolicityoflength3positive}
Let $n>0, 0\leq k,l <n$, where $k\neq l$ and $\gcd (n,k,l)=1$, and let $G=G_n(x_0x_kx_l)$. Then exactly one of the following holds:
\begin{itemize}
    \item[(a)] $(B,C,D)=(F,F,T)$, in which case $G\cong \Gamma_{2m}$, where $m=n/2$, and is finite, and therefore elementary hyperbolic, if $n \in \{2,4\}$, non-elementary hyperbolic if $n\in \{6,12,18\}$, and is non-hyperbolic otherwise;

    \item[(b)] $(B,C,D)=(F,F,F)$, in which case the defining presentation of $G$ satisfies the small cancellation condition T(6), and $G$ is not hyperbolic if $n=7$ or $n=8$, or 
    \begin{itemize}
        \item[(i)] $n=21$ and ($l\equiv 5k$ or $k\equiv 5l\bmod n$), or
        \item[(ii)] $n=24$ and ($l\equiv 5k \bmod n$ or $k\equiv -4l\bmod n$ or $l\equiv -4k \bmod n$),
    \end{itemize}
    and is non-elementary hyperbolic otherwise;

    \item[(c)] either
    \begin{itemize}
        \item[(i)] $(A,B)=(T,T)$, in which case $G\cong \Z*\Z$; or
        \item[(ii)] $(A,B,C)=(T,F,T)$, in which case $G\cong \Z*\Z*\Z_{(2^{n/3}-(-1)^{n/3})/3}$;
    \end{itemize}
    and so $G$ is non-elementary hyperbolic;
    
    \item[(d)] either
    \begin{itemize}
        \item[(i)] $(A,B,C)=(F,F,T)$, in which case $G$ is a finite metacyclic group of order $2^{n}-(-1)^n$; or
        \item[(ii)] $(A,B,C)=(F,T,F)$, in which case $G\cong \Z_3$;
    \end{itemize}
    and so $G$ is elementary hyperbolic.
\end{itemize}
\end{maincorollary}

As noted in Corollary \ref{cor:hyperbolicityoflength3positive}(a), if $(B,C,D)=(F,F,T)$ then the corresponding group $G_n(x_0x_kx_l)$ is isomorphic to $G_n(x_0x_1x_{n/2-1})$. In \cite[Conjecture 7.12]{MohamedWilliams} it was conjectured that for even $n\geq 20$, if  $(B,C,D)=(F,F,F)$ then $G_n(x_0x_kx_l)$ is not isomorphic to $G_n(x_0x_1x_{n/2-1})$. In Corollary \ref{cor:SFFFF(n)bothparts} we prove this conjecture.

\begin{maincorollary}\label{cor:SFFFF(n)bothparts}
Suppose $n\geq 20$ is even, $0\leq k,l <n$, where $\gcd(n,k,l)=1$ and let $G=G_n(x_0x_kx_l)$. If $(B,C,D)=(F,F,F)$ then $G\not \cong G_n(x_0x_1x_{n/2-1})$.
\end{maincorollary}

\section{The proofs of the main results}

The proof of Theorem \ref{thm:summarytheorem} proceeds as follows. Writing $x_i=t^iyt^{-i}$, the split extension $E_{2m}=\Gamma_{2m}\rtimes \pres{t}{t^{2m}}$ is given by the presentation $\pres{y,t}{t^{2m}, y tyt^{m-2}yt^{-(m-1)}}$. Introducing the generator $x=yt$ and eliminating $y$ gives
\begin{equation}
    E_{2m}= \pres{x,t}{t^{2m}, x^2t^{m-3}xt^m}.\label{eq:E2m}
\end{equation}
We introduce the following elements:
\begin{equation}
 A=xt^{-3}, B=  t^m xt^{m+3}xt^{m-3}.\label{eq:AandB}
\end{equation}
First observe the following:

\begin{lemma}\label{lem:ABcommute}
The elements ${A},{B}$ commute in $E_{2m}$.
\end{lemma}

\begin{proof}
We first note some consequences of the relators of $E_{2m}$:
\begin{alignat}{1}
t^mx^{-2}t^m&=t^{-3}x,\label{eq:tmX2tm=T3x}\\
t^3x^{-2}&=t^mxt^m,\label{eq:t3X2=tmxtm}\\
t^mx^2t^m&=x^{-1}t^3,\label{eq:tmx2tm=Xt3}\\
xt^{-3}&=(x^{-1}t^m)^2.\label{eq:xT3=XtmXtm}
\end{alignat}
Then
\begin{alignat*}{1}
    {B}{A}
    &= t^m xt^{m+3}x t^m (t^{-3} x)t^{-3}\\
    &= t^m xt^{m+3}x t^m (t^mx^{-2}t^m )t^{-3}\quad \mathrm{by}~\eqref{eq:tmX2tm=T3x}\\
    &= t^m xt^m (t^{3}x^{-2}) x t^{m-3}\\
    &= t^m xt^m (t^mxt^m) x t^{m-3}\quad \mathrm{by}~\eqref{eq:t3X2=tmxtm}\\
    &= (t^m x^2t^m) x t^{m-3}\\
    &= (x^{-1}t^3) x t^{m-3}\quad \mathrm{by}~\eqref{eq:tmx2tm=Xt3}\\
    &=(x^{-1}t^m)^2 \cdot t^m(xt^{m+3}xt^{-3})t^m\\
    &= (xt^{-3}) \cdot t^m(xt^{m+3}xt^{-3})t^m\quad \mathrm{by}~\eqref{eq:xT3=XtmXtm}\\
    &= {A}{B}.
\end{alignat*}
\end{proof}

The main business of the proof of Theorem \ref{thm:summarytheorem} is in showing that, for $m\not \in \{1,2,3,4,6,9,12\}$, $A$ and $B$ generate a non-cyclic, free abelian subgroup of $E_{2m}$, and so $E_{2m}$ and $\Gamma_{2m}$ (being a finite index subgroup of $E_{2m})$ are not hyperbolic. This is done in the following lemma, whose proof we defer to Section \ref{sec:mainlemmaproof}. (Non-hyperbolicity in the cases $m=4,12$ is proved by separate methods, and hyperbolicity in the cases $m=1,2,3,6$ or $9$ is known by prior results.)

\begin{lemma*}
Suppose $m\not \in \{1,2,3,4,6,9,12\}$ then ${A}^\alpha {B}^\beta=1$ in $E_{2m}$ if and only if $\alpha=\beta=0$.
\end{lemma*}

Subject to this, we are now in a position to prove Theorem \ref{thm:summarytheorem} and its corollaries.

\begin{proof}[Proof of Theorem \ref{thm:summarytheorem}]
The cases $m= 1,2,3,6,9$ follow from \cite{EdjvetWilliams} (or see \cite[Table 1]{MohamedWilliams}). If $m=4$ the group $\Gamma_{2m}=\Gamma_8$ contains a free abelian subgroup of rank 2 (see \cite[Example 3(i)]{BW2}) so is not hyperbolic. Now consider the case $m=12$. The group
\begin{alignat*}{1}
    E_{24}&=\pres{x,t}{t^{24},x^2t^9xt^{12}}\\
          &=\pres{x,t,z}{t^{24},x^2t^9xt^{12},z=xt^{-9}}\\
          &=\pres{z,t}{t^{24},zt^9zt^{18}zt^{21}},
\end{alignat*}
which is isomorphic to the split extension $G_{24}(x_0x_9x_3)\rtimes \pres{t}{t^{24}}$ (by writing $x_i=t^izt^{-i})$. The group $G_{24}(x_0x_9x_3)$ is isomorphic to the free product of three copies of $G_8(x_0x_3x_1)\cong \Gamma_8$ (see, for example \cite{EdjvetIrreducible}), which as noted above contains a free abelian subgroup of rank 2. Therefore $E_{24}$ also contains that subgroup so it, and hence $\Gamma_{24}$, is not hyperbolic. Thus we may assume $m\neq 1,2,3,4,6,9,12$.

Since $\Gamma_{2m}$ is a finite index subgroup of $E_{2m}$ the group $\Gamma_{2m}$ is hyperbolic if and only if $E_{2m}$ is hyperbolic. By Lemma \ref{lem:ABcommute} the elements $A,B$ commute, and by Main Lemma each have infinite order and so $A,B$ generate a free abelian subgroup of $E_{2m}$. If $A,B$ generate $\Z$ then there exist non-zero $\alpha,\beta \in \Z$ such that $A^\alpha B^\beta=1$ in $E_{2m}$, a contradiction to Main Lemma. Therefore $A,B$ generate a free abelian subgroup of rank 2 and so $E_{2m}$, and hence $\Gamma_{2m}$, are not hyperbolic.
\end{proof}

\begin{proof}[Proof of Corollary \ref{cor:hyperbolicityoflength3positive}]
If $(B,C,D)=(F,F,T)$ then by \cite[Lemma 2.2]{MohamedWilliams} $G\cong \Gamma_{2m}$, where $m=n/2$, as in (a), so the result follows from Theorem \ref{thm:summarytheorem}. If $(B,C,D)=(F,F,F)$ then the defining cyclic presentation for $G$ satisfies the T(6) small cancellation condition (\cite[Lemma 5.1]{EdjvetWilliams}), as in (b), and the result follows from \cite[Theorem A]{ChinyereWilliamsT6}. Thus we may assume $(B,C)\neq (F,F)$. If $A=T$ then either $(A,B)=(T,T)$ or $(A,B,C)=(T,F,T)$, as in (c). If $(A,B)=(T,T)$ then $G\cong \Z*\Z$ by \cite[Lemma 2.4]{EdjvetWilliams}, and if $(A,B,C)=(T,F,T)$ then $G\cong \Z*\Z*\Z_{(2^{n/3}-(-1)^{n/3})/3}$ by \cite[Lemma 2.5]{EdjvetWilliams} or \cite[Corollary D]{BW2}. Thus we may assume $A=F$, and it follows from the defining congruences that $(B,C)\neq (T,T)$; therefore $(A,B,C)=(F,F,T)$ or $(F,T,F)$, as in (d). If $(A,B,C)=(F,F,T)$ then $G$ is a finite metacyclic group of order $2^n-(-1)^n$ by \cite[Lemma 3.2]{EdjvetWilliams}, \cite[Lemma 5.5]{BogleyShift}, or \cite[Corollary D]{BW2}, and if $(A,B,C)=(F,T,F)$ then $G\cong \Z_3$ by \cite[Lemma 2.4]{EdjvetWilliams}.
\end{proof}

\begin{proof}[Proof of Corollary \ref{cor:SFFFF(n)bothparts}]
Let $m=n/2$ so that $G_n(x_0x_1x_{n/2-1})=\Gamma_{2m}$, and suppose for contradiction that $G\cong \Gamma_{2m}$. The computations carried out in \cite[Section 7]{MohamedWilliams} show, in particular, that if $20\leq n\leq 24$ then $G^\mathrm{ab}\not \cong \Gamma_{2m}^\mathrm{ab}$, a contradiction, so we may assume $n>24$. If $(B,C,D)=(F,F,F)$ then $G$ satisfies the small cancellation condition $T(6)$ so by \cite[Theorem A]{ChinyereWilliamsT6} it is hyperbolic. But by Theorem \ref{thm:summarytheorem} the group $\Gamma_{2m}$ is not hyperbolic, a contradiction.
\end{proof}

\section{The proof of Main Lemma}\label{sec:mainlemmaproof}

\subsection{Relative presentations, diagrams, and curvature}\label{sec:generalsetup}

We must show that, for $m\not \in \{1,2,3,4,6,9,12\}$, if $(\alpha,\beta)\neq (0,0)$ then $A^\alpha B^\beta\neq 1$ in the group $E_{2m}$, as defined at \eqref{eq:E2m}. Without loss of generality we may assume $\beta\geq 0$ and, by way of contradiction, that $A^\alpha B^\beta=1$. The group $E_{2m}$ has a one-relator relative presentation (\cite{BogleyPride}) 
\begin{alignat}{1}
E_{2m}=\pres{H,x}{x^2t^{m-3}xt^m}\label{eq:relpresE2m}
\end{alignat}
where $H=\pres{t}{t^{2m}}$. The condition $A^\alpha B^\beta=1$ implies that there exists a relative diagram $\LL$ (\cite{Howie83}) over the relative presentation for $E_{2m}$ having boundary label $A^\alpha B^\beta$. The proof proceeds by using curvature arguments to show that no such relative diagram can exist. Remark \ref{rem:Huebschmann} illuminates the strategy.

In analysing such relative diagrams it is convenient to introduce the following notation:
\[ \lambda = t^0, a= t^{m-3}, b=t^m, c=t^{-3}, d=t^{m+3}, e=t^{-3}, f=t^{m-3}, g=t^m, h=t^m. \]
Then the relator of the relative presentation \eqref{eq:relpresE2m} can be written
\begin{alignat}{1}
    x^2t^{m-3}xt^m=x\lambda x a x b\label{eq:relator}
\end{alignat}
and, with $A,B$ as defined at \eqref{eq:AandB}, we have 
\begin{alignat}{1}
    {A}^\alpha {B}^\beta \sim \begin{cases}
        (xc)^{\alpha -1} xf (xdxe)^{\beta -1} xdxf & \mathrm{if}~\alpha>0, \beta>0,\\
        x^{-1} (c^{-1}x^{-1})^{|\alpha|-1} g(xdxe)^{\beta -1} xdxh& \mathrm{if}~\alpha<0,\beta>0,\\
        (xc)^{\alpha} & \mathrm{if}~\alpha>0, \beta=0,\\
        (xdxe)^{\beta}& \mathrm{if}~\alpha=0,\beta>0,\label{eq:starandstarstar}
    \end{cases}
\end{alignat}
where $\sim$ denotes cyclic permutation.

Therefore the boundary of $\LL$ is given by the relevant expression in \eqref{eq:starandstarstar}, each interior corner of $\LL$ has label $a,b$ or $\lambda$, and by \eqref{eq:starandstarstar}, each exterior corner has label $c,d,e,f,g$ or $h$. The sum of the powers of $t$ read around any vertex is congruent to $0\bmod 2m$ and the product of the directed corner labels and edge labels of any given region of $\LL$ yields (up to cyclic permutation and inversion) the relator $x^2t^{m-3}xt^m$. Thus the regions of $\LL$ are given by Figure {{A}}(i),(ii), which we draw as shown in Figure {{A}}(iii),(iv), with the understanding that vertex labels $l(v)$ are read anti-clockwise and, as throughout, the labels $\hat{a},\hat{b},\ldots$ denote $a^{-1},b^{-1},\ldots$ respectively, and the label $\mu$ denotes $\lambda^{-1}$. Often, the vertex labels are considered up to cyclic permutation and inversion; when this is clear from the context we will not explicitly state this.

Without loss of generality we make the following assumptions:

\begin{itemize}
    \item[(A1)] $\LL$ is minimal with respect to the number of regions (which implies that $\LL$ is reduced).
    \item[(A2)] Subject to (A1), $\LL$ is maximal with respect to the number of interior vertices of degree $2$.
    \item[(A3)] Subject to (A1) and (A2), $\LL$ is maximal with respect to the number of vertices having label $b\mu b\mu$.
\end{itemize}

We now describe how to construct a spherical diagram $\SS$ from $\LL$. If $\alpha \beta\neq 0$ and $\LL$ consists of a single disc $\DD$ then it contains exactly two boundary vertices $v_1,v_2$ whose exterior corner labels are both $f$ (in the case $\alpha>0,\beta>0$) or $g$ and $h$ (in the case $\alpha<0,\beta >0$), which we call \em exceptional vertices, \em denoted $u_f,u_g,u_h$, respectively. Note that there are no such exceptional vertices in the cases $\alpha>0,\beta=0$ or $\alpha=0$, $\beta>0$. If $\LL$ does not consist of a single disc then it has a vertex $u$ with at least two corners whose labels are boundary labels and which is a vertex of an extremal disc $\DD$ containing at most one exceptional vertex. In this case detach $\DD$ from $\LL$ and in doing so create another boundary vertex $v_3$ (namely the vertex $u$, above). We also call $v_3$ an \em exceptional vertex \em and denote it $u_*$.  The exterior corner of vertex $u_*$ in the detached $\DD$ now has multiple labels, consisting of the labels of the exterior corners of $u_*$ in $\LL$. Note that in this case the boundary of $\DD$ does not spell the relator \eqref{eq:starandstarstar}, and that $\DD$ has at least one exceptional vertex (namely $v_3$) and at most two exceptional vertices. We use $u^*$ to denote an exceptional vertex of any type, i.e. $u^*\in \{u_f,u_g,u_h,u_*\}$.

We form a spherical diagram $\SS$ whose southern hemisphere is $\DD$ (and so is tesselated by triangles) and whose northern hemisphere consists of a single region whose boundary label is a cyclic subword of \eqref{eq:starandstarstar} (which is a proper subword if $\LL$ does not consist of a single disc). This exceptional region is denoted by $\D^*$. A region of $\SS$ sharing an edge with $\D^*$ is called a \em boundary region\em; otherwise it is called \em interior\em. There are three types of vertices: the (at most) two exceptional vertices; vertices whose label involves only $a,b,\lambda$ and these are called \em interior vertices\em; and the remaining vertices are called \em boundary vertices. \em Thus the label of a boundary vertex is of the form $\theta w$ where $\theta^{\pm 1}\in \{c,d,e\}$ and $w$ involves only $a,b$ or $\lambda$.

We now turn to curvature. If $v$ is a vertex of $\SS$ having degree $d(v)=d$, assign $2\pi/d$ to each corner at $v$. This way the curvature of each vertex is $0$. The \em curvature \em of a region $\D$ of $\SS$ of degree $k$ and whose vertices have degree $d_i\geq 2$ ($1\leq i\leq k$) is defined to be
\[ c(\D) = c(d_1,\ldots ,d_k) = (2-k)\pi +2\pi \sum_{i=1}^k \frac{1}{d_i}. \]
Observe that if $v$ is not an exceptional vertex then $d(v) \geq 4$ by Lemma \ref{lem:stargraphlemma}(i),(ii) and, moreover, $l(v)$ corresponds to a reduced closed path of even length in the star graph $\Gamma$. Thus if $d(\D)=3$ and $c(\D)>0$ then 
\[c(\D) \in \{ c(4,4,k)=2\pi/k, c(4,6,6)=\pi/6, c(4,6,8)=\pi/12, c(4,6,10)=\pi/30 \}.\]
It follows from the Gauss-Bonnet theorem that the total curvature of $\SS$, $\sum_{\D \in \LL} c(\D) = 4\pi$ (see \cite[Section 4]{McCammondWise} and the references therein).  Our contradiction will be obtained by showing that this cannot occur. 

In our curvature analysis of the relative diagram $\LL$, in order to locate regions of positive curvature it will be important to understand the labelling of its low degree vertices. In Section \ref{sec:lowdegreelabels} we use bridge moves on relative diagrams and star graphs to obtain restrictions on the possible labellings that can occur. We are then in a position to define a \em curvature distribution scheme \em in Sections \ref{sec:StageI} and \ref{sec:StageII}. That is, we locate each region $\D$ ($\D \neq \D^*$) satisfying $c(\D)>0$ and distribute its positive curvature $c(\D)$ to negatively curved near regions $\hD$ of $\D$ so that the curvature of $\D$ is reduced to zero. For such regions $\hD$, define $c^*(\hD)$ to equal $c(\hD)$ plus all the positive curvature $\hD$ receives minus all the positive curvature distributed from $\hD$ by the application of the curvature distribution scheme. We complete the proof by showing that $\sum c(\D)=\sum c^*(\hD)<4 \pi$, a contradiction.

The curvature distribution scheme is divided into two stages. In Section \ref{sec:StageI} we define the curvature distribution scheme for Stage I. This considers interior regions $\D$ with $c(\D)>0$ that do not contain any exceptional vertices, and distributes their curvature to near regions (which are either boundary or interior regions). In Section \ref{sec:StageIconsequences} we record detailed implications of how curvature has been transferred in Stage I, but we note here that a fundamental consequence of Stage I is that by its conclusion the only regions with positive curvature are boundary regions or interior regions with at least one exceptional vertex.

Next, in Section \ref{sec:StageII}, we define the curvature distribution scheme for Stage II, which considers these regions. Such regions $\hD$ may have received curvature in Stage I, so instead of distributing $c(\hD)$ we distribute their new curvature (which we will denote $c^*(\hD)$), when positive, to the exceptional region $\D^*$. To conclude the proof, in Section \ref{sec:concludingtheproof} we analyse the curvature of the resulting regions to reach our desired contradiction that the total curvature $4\pi$ of $\SS$ cannot be obtained. 

\begin{remark}
Curvature redistribution methods have been used (for example in \cite{ChalkEdjvet,ChalkEdjvetJuhasz,ChinyereWilliamsT5,ChinyereWilliamsT6}) to prove hyperbolicity of certain classes of cyclically presented groups. The curvature redistribution methods used in the current article, to prove non-hyperbolicity, are applied in a different way, as here the global negative curvature obtained is a property of the specific relative diagrams that need to be considered in order to prove Main Lemma. It does not prove any negative curvature property of the groups themselves.
\end{remark}

\begin{remark}[J.Huebschmann \cite{HuebschmannPersonalCommunication}]\label{rem:Huebschmann}
Let $r=x^2t^{m-3}xt^m$ denote the relator of the relative presentation \eqref{eq:relpresE2m} and let $G=\gpres{x}*H$ where, as before, $H=\pres{t}{t^{2m}}$. Main Lemma shows that for $m\not \in \{1,2,3,4,6,9,12\}$, $(\alpha,\beta)\neq (0,0)$ the equation
\begin{alignat}{1}
A^\alpha B^\beta =\prod_{j=1}^u y_jr^{\epsilon_j} y_j^{-1}\label{eq:equation}
\end{alignat}
does not admit a solution with $y_j\in G$, $\epsilon_j=\pm 1$, $u\geq 1$.

For a $G$-group $K$, write the $G$-action on $K$ as
\[ G \times K \longrightarrow K, \ (y,b) \mapsto {}^y b ,\ y \in G,\ b \in K.\]
The member $r$ of $G$ generates the free $G$-crossed module $C\stackrel{\partial} \rightarrow G$ having $\mathrm{coker}(\partial)\cong E_{2m}$: Let $\widehat C$ be the free $G$-group generated by $r$, let $\widehat \partial \colon \widehat C \to G$ be the canonical homomorphism that sends ${}^y r\in \widehat C$ to $yry^{-1}\in G$, for $y \in G$, and let $C$ be the quotient $\widehat C/P$ of $\widehat C$ modulo the normal $G$-subgroup $P$ that the Peiffer elements 
\[aba^{-1}  \left({}^{\widehat \partial a} b\right)^{-1},\ a, b \in \widehat C,\]
generate. A little thought reveals that the subgroup $P$ of $\widehat C$ generated by the Peiffer elements is normal and $G$-invariant. The $G$-homomorphism $\widehat \partial$ passes to a $G$-homomorphism $\partial \colon C \to G$.

The relative diagrams used in the proof of Main Lemma are dual to pictures over the relative presentation, described in \cite{BogleyPride}. The idea of a diagram over a presentation (``Randwegaggregat") goes back at least to \cite{Peiffer49}. The paper \cite{BrownHuebschmann} discusses (ordinary) diagrams and pictures. Section 10 of that paper reveals that classes of diagrams form $C$ in such a way that the assignment to a diagram of its boundary path induces the homomorphism $\partial \colon C\rightarrow G$. This observation explains why a non-trivial solution to \eqref{eq:equation} implies that there exists a relative diagram over the relative presentation for $E_{2m}$ having boundary label $A^\alpha B^\beta$. While \cite{BrownHuebschmann} handles only the absolute case, the reasoning there carries over to the relative case. An account of the history of these ideas is presented in \cite{Huebschmann23}. 
\end{remark}

\subsection{Labelling of low degree vertices}\label{sec:lowdegreelabels}

The following lemma provides some important restrictions on the possible labellings of vertices of $\LL$.

\begin{lemma}
Let $v$ be a vertex of $\LL$ and let $k>0$. Then (up to cyclic permutation and inversion):
\begin{itemize}
    \item[(i)] $l(v)$ does not have a sublabel of the form $\theta w \theta^{-1}$ where $w = 1$;
    \item[(ii)] if $l(v)= (b\mu b\mu)^k$ then $k=1$.
\end{itemize}
\end{lemma}

\begin{proof}
(i) If $l(v)$ does have such a sublabel then a bridge move (see, for example, \cite[Figure 2]{Howie83}) followed by cancellation of inverse regions contradicts (A1). (ii) If $k>1$ then the number of vertices having label $b\mu b\mu$ can be increased by a bridge move at $v$ without affecting the number of interior vertices of degree 2, contradicting (A3).
\end{proof}

The vertex labels of a reduced relative diagram correspond to reduced closed paths in its star graph, and so analysis of such closed paths places restrictions on the possible vertex labels of that diagram (see \cite[Section 2.1]{BogleyPride}, \cite[Section 2]{EdjvetHowie91}). The star graph $\Gamma$ of $\LL$ is the (directed, labelled) graph with vertices $x,x^{-1}$ and with a directed edge from $x^{\epsilon_1}$ to $x^{-\epsilon_2}$ ($\epsilon_1, \epsilon_2\in \{\pm 1\}$) for each cyclic subword $x^{\epsilon_1} u x^{\epsilon_2}$ of either \eqref{eq:relator} or \eqref{eq:starandstarstar}, where $u\in H$; such a directed edge is labelled $u$. Thus, $\Gamma$ is given by Figure {{A}}(v),(vi) according to ($\alpha>0$ and $\beta\geq 0$) and ($\alpha\leq 0,\beta > 0)$ respectively. For ease of presentation we have introduced the inverse edge (from $x^{-1}$ to $x$, labelled $\mu$) to the edge from $x$ to $x^{-1}$, labelled $\lambda$, with the understanding that the edges labelled $\lambda,\mu$ (unlike the others) are only traversed in the direction indicated. For the case $\alpha>0, \beta=0$ the star graph is the subgraph of Figure {{A}}(v) containing only the edges labelled $a,b,c,\lambda, \mu$; and for case $\alpha=0, \beta>0$ the star graph is the subgraph of Figure {{A}}(vi) containing only the edges labelled $a,b,d,e,\lambda,\mu$. Analysis of short reduced closed paths in $\Gamma$ yields the following result that provides further restrictions on the possible labelling of low degree vertices of $\LL$. We will use this result throughout, usually without explicit reference.

\begin{lemma}\label{lem:stargraphlemma}
Let $v$ be a vertex of $\LL$.
\begin{itemize}
    \item[(i)] If $v$ is an interior vertex and $d(v)\leq 8$ then (up to cyclic permutation and inversion) $l(v)\in \{b\mu b\mu,$ $ab^{-1}\lambda a^{-1}b\mu, ab^{-1}\lambda a^{-1}\lambda b^{-1}\}$.

    \item[(ii)] If $l(v) \in \{cw,dw,ew\}$ where $w$ involves only $a,b$ or $\lambda$ or their inverses and $d(v)\leq 6$ then
    \begin{alignat*}{1}
        l(v) &\in \{ ca^{-1}b\mu, ca^{-1}\lambda b^{-1}, cb^{-1}\lambda a^{-1}, c\mu ba^{-1}, ea^{-1}b\mu, ea^{-1}\lambda b^{-1}, eb^{-1}\lambda a^{-1}, e\mu ba^{-1},\\
             &\qquad db^{-1}ab^{-1}, d\mu a\mu, d\mu ab^{-1}\lambda b^{-1}, d\mu b\mu ab^{-1}, db^{-1}a\mu b\mu, db^{-1}\lambda b^{-1}a\mu \}.
    \end{alignat*}

    \item[(iii)] If $l(v)\in \{fw,gw,hw\}$ where $w$ involves only $a,b$ or $\lambda$ or their inverses and $d(v)\leq 4$ then $l(v) \in \{fa^{-1},gb^{-1}\lambda , g\mu b, hb\mu, h\lambda b^{-1}\}$.

    \item[(iv)] If $v$ is not exceptional, $l(v)\in \{a\mu a\mu w_1, a^{-1}ba^{-1}bw_2, ba^{-1}ba^{-1}w_3, \mu a \mu a w_4\}$ and $d(v)=8$ then
        \begin{alignat*}{1}
            w_1&\in \{ba^{-1}\lambda c^{-1}, ba^{-1}\lambda e^{-1}, a\mu ad^{-1}\ (m=5,15), ba^{-1}db^{-1}, ad^{-1}a\mu\ (m=5,15)\},\\
            w_2&\in \{\mu ab^{-1}c, \mu ab^{-1}e, a^{-1}\lambda a^{-1}d \ (m=5,15), \mu ad^{-1}\lambda, a^{-1}da^{-1}\lambda\ (m=5,15)\},\\
            w_3&\in \{cb^{-1}a\mu, eb^{-1}a\mu, da^{-1}\lambda a^{-1}\ (m=5,15), \lambda d^{-1}a\mu, \lambda a^{-1}da^{-1}\ (m=5,15)\},\\
            w_4&\in \{ c^{-1}\lambda a^{-1}b, e^{-1}\lambda a^{-1} b, d^{-1}a\mu a\ (m=5,15), \lambda d^{-1}ab, \mu ad^{-1}a\ (m=5,15)\}.
        \end{alignat*}
\end{itemize}
\end{lemma}

Note that the restrictions $m\not \in \{1,2,3,4,6,9,12\}$ are necessary for Lemma \ref{lem:stargraphlemma}. For example $l(v)\neq ab^{-1}a\mu$ in part (i) since $m-6\equiv 0\bmod 2m$ implies $m=2$ or $6$. Note also that, in particular, interior vertices of degree 4 have label $b\mu b\mu$ (up to cyclic permutation and inversion).

\subsection{Stage I: $\D$ is interior and does not contain an exceptional vertex}\label{sec:StageI}

We set out the curvature distribution scheme for Stage I as follows. In Section \ref{sec:generalrule} we introduce a general rule for distributing curvature from interior regions that share at least one edge with a boundary region. In Section \ref{sec:StageIa} we define the curvature distribution from interior regions that have exactly one boundary vertex of degree 4 (Stage I(a)). We then do the same for interior regions that have no boundary vertices of degree 4 (Stage I(b)). Stage I(b) splits into two largely analogous cases. In Section \ref{sec:StageIb} we give the curvature distribution for one of these cases, and in Section \ref{sec:symmetricases} we describe the differences that are needed for the second case.

\subsubsection{General Rule}\label{sec:generalrule}

\noindent \textbf{General Rule (GR):} If $\D$ shares an edge with a boundary region $\hD$ then distribute $\min \{ c(\D), \pi/12 \}$ from $\D$ to $\hD$ as in Figure {{B}}(i); if with two boundary regions $\hD_1,\hD_2$ then distribute $\min \{ c(\D)/2, \pi/12 \}$ to each of $\hD_1,\hD_2$; or if with three boundary regions $\hD_1,\hD_2,\hD_3$ then distribute $\min \{ c(\D)/3, \pi/12 \}$ to each of $\hD_1,\hD_2,\hD_3$.

\begin{remark*}
\begin{itemize}
    \item[(i)] Having distributed curvature according to GR, it may be the case that $\D$ still retains some positive curvature. In what follows we give details of how this extra curvature is distributed.    
    \item[(ii)] Exceptions to GR will be clearly indicated.
\end{itemize}
\end{remark*}

An exception to GR is when $\D$ shares an edge with a boundary region $\hD$, shown in Figure {{B}}(ii), where $u^*$ is an exceptional vertex of degree 2. In this case distribute $c(\D)$ from $\D$ to $\hD$, as shown. If $\D$ shares an edge with two such $\hD$ then distribute $c(\D)/2$ to each of them. In particular, this occurs when the interior region $\D$ contains two boundary vertices of degree 4. (To see this observe that if $\D$ contains two boundary vertices $u_1,u_2$ of degree 4, then this forces a boundary vertex, $v$, say of degree 2 between them. Either $v$ is an exceptional vertex $u_*$ created in the detaching process or, by Lemma \ref{lem:stargraphlemma}, $l(v)=fa^{-1}$ and so $v$ is an exceptional vertex $u_f$. Note that this argument rules out the possibility of an interior region having three boundary vertices of degree 4, for otherwise it would force three exceptional vertices.) Assume, therefore, that $\D$ does not contain two boundary vertices of degree 4; then either  $\D$ contains exactly one boundary vertex of degree 4 or $\D$ contains no boundary vertices of degree 4. We consider these situations in Stage I(a) and Stage I(b), respectively.
\begin{notation*}
In the figures, a circled value $\encircle{\geq k}$ indicates that the vertex is exceptional (of any type) of some unspecified degree, or not exceptional but of degree at least $k$.
\end{notation*}

\subsubsection{Stage I(a) : $\D$ is interior and contains exactly one boundary vertex of degree 4}\label{sec:StageIa}

Throughout this section, suppose that $\D$ is interior and contains exactly one boundary vertex of degree 4, $v$ say. If the two remaining (interior) vertices of $\D$ each have degree at least 6 then, in accordance with GR, distribute $c(\D)/2\leq \pi/12$ to each of $\hD_1$ and $\hD_2$ as in Figure {{B}}(iii); or if $\hD_3$ in Figure {{B}}(iii) is also a boundary region then distribute $c(\D)/3\leq \pi/18$ to each of $\hD_1,\hD_2,\hD_3$. Otherwise $\D$ must contain at least one interior vertex of degree less than $6$. Such vertices have label $b \mu b\mu$ by Lemma \ref{lem:stargraphlemma}(i), and so only one interior vertex of $\D$ can have degree 4. Given this it follows from Lemma \ref{lem:stargraphlemma} that (up to cyclic permutation and inversion)
\[l(v) \in \{ca^{-1}\lambda b^{-1}, c\mu b a^{-1}, db^{-1}ab^{-1}, d\mu a \mu, ea^{-1}\lambda b^{-1}, e\mu ba^{-1} \},\]
the remaining degree 4 labels resulting in an impossible configuration. We discuss these possible labellings over the following bullet points:
\begin{itemize}
    \item Let $l(v)= ca^{-1}\lambda b^{-1}$. Completing the four regions near $\D$ we obtain Figure {{B}}(iv) without the $\encircle{\geq 6}$. The labelling so far completed implies that the corresponding vertex is not interior of degree 4 (by Lemma \ref{lem:stargraphlemma}(i)), so either it is exceptional or of degree at least 6, as in Figure {{B}}(iv), or it has degree 4, in which case it is boundary with label $ca^{-1}\lambda b^{-1}$ or $ea^{-1}\lambda b^{-1}$ as in Figure {{B}}(v).
    
    In Figure {{B}}(iv) distribute $\min \{c(\D)/2,\pi/12\} \leq \pi/12$ from $\D$ to each of $\hD_1$ and $\hD_2$; and distribute the remaining $\max\{0, c(\D)-\pi/6\} \leq \pi/6$ to $\hD_2$ (giving a total of at most $\pi/4$, as shown). Note that if $\hD_3$ is a boundary region then this represents an exception to GR. In Figure {{B}}(v) (in accordance with GR) distribute $\min\{c(\D)/3,\pi/12\} \leq \pi/12$ from $\D$ to each of $\hD_j$ ($1\leq j \leq 3$) and distribute the remaining $\max\{0, c(\D)-\pi/4\} \leq \pi/12$ to $\hD_2$ (giving a total of at most $\pi/6$, as shown).
  
    \item Let $l(v)= c\mu b a^{-1}$. Similar to the case $l(v)=ca^{-1}\lambda b^{-1}$, the possibilities are as shown in Figures {{B}}(vi),(vii). In Figure {{B}}(vi) distribute $\min\{c(\D)/2,\pi/12\} \leq \pi/12$ from $\D$ to each of $\hD_1$ and $\hD_2$; and distribute the remaining $\max \{0, c(\D)-\pi/6\} \leq \pi/6$ to $\hD_1$ (giving a total of at most $\pi/4$, as shown). Note that if $\hD_3$ is a boundary region then this represents an exception to GR. In Figure {{B}}(vii) distribute $\min\{c(\D)/3,\pi/12\} \leq \pi/12$ from $\D$ to each of $\hD_j$ ($1\leq j \leq 3$) and distribute the remaining $\max\{0, c(\D)-\pi/4\} \leq \pi/12$ to $\hD_1$ (giving a total of at most $\pi/6$, as shown).
    
    \item Let $l(v)= db^{-1}ab^{-1}$. Then the configuration is as in Figure {{B}}(viii). Distribute $c(\D)\leq c(4,4,8)= \pi/4$ to $\hD_2$ as in Figure {{B}}(viii), an exception to GR. Suppose that the region $\D_1$ of Figure {{B}}(viii) is interior of positive curvature and does not contain an exceptional vertex. Then the vertex $w$ has degree 4. If $d(u)\geq 6$ then distribute $c(\D_1)\leq c(4,6,8)=\pi/12$ to the region $\hD$ as in Figure {{B}}(ix), an exception to GR. If $l(u)=d\mu a\mu$ distribute $c(\D_1)\leq c(4,4,8)=\pi/4$ to $\hD$ as shown in Figure {{B}}(x), an exception to GR, noting that $u^*$ is forced to be an exceptional vertex of degree 2 as shown (which is in fact a special case of {{B}}(ii)). 
    
    \item Let $l(v)= d\mu a \mu$. Then the configuration is as in Figure {{B}}(xi). Distribute $c(\D)\leq c(4,4,8)= \pi/4$ to $\hD_1$ as in Figure {{B}}(xi), an exception to GR. Suppose that the region $\D_1$ of Figure {{B}}(xi) is interior of positive curvature and does not contain an exceptional vertex. Then the vertex $w$ has degree 4. If $d(u)\geq 6$ then  distribute $c(\D_1)\leq c(4,6,8)=\pi/12$ to the region $\hD$ as in Figure {{B}}(xii), an exception to GR. If $l(u)=db^{-1} ab^{-1}$ distribute $c(\D_1)\leq c(4,4,8)=\pi/4$ to $\hD$ as shown in Figure {{B}}(xiii), an exception to GR, noting that $u^*$ is forced to be an exceptional vertex of degree 2 as shown (which is in fact a special case of {{B}}(ii)).
        
    \item Let $l(v)= ea^{-1}\lambda b^{-1}$. Then the configuration is as in Figure {{B}}(xiv). Distribute $c(\D)\leq \pi/3$ to $\hD$ as shown, an exception to GR.
    
    \item Let $l(v)= e\mu ba^{-1}$. Then the configuration is as in Figure {{B}}(xv). Distribute $c(\D)\leq \pi/3$ to $\hD$ as shown, an exception to GR.
\end{itemize}

\subsubsection{Stage I(b) : $\D$ is interior and contains no boundary vertices of degree 4}\label{sec:StageIb}

Suppose now that $\D$ is interior and does not contain a boundary vertex of degree 4. Then $c(\D) > 0$ implies that $\D$ contains an interior vertex of degree 4, which therefore (up to cyclic permutation and inversion) has label $b \mu b \mu$. So (up to inversion) there are two cases for $\D$, namely $\D$ of Figure {{C}}(i) and {{G}}(i). The curvature distributions for these two cases are largely analogous. We first describe the distribution for Figure {{C}}(i), and in Section \ref{sec:symmetricases} we explain how to obtain the distribution for the case of Figure {{G}}(i) from this. Consider Figure {{C}}(i). The region $\D$ has a vertex of degree 4, and interior regions can have at most one vertex of degree 4, so the other two vertices each have degree at least 6. If the vertex $u_3$ of Figure {{C}}(i) is an exceptional vertex then distribute $c(\D)\leq c(4,6,6)=\pi/6$ to $\hD$. Assume from now on that $u_3$ is not exceptional (so, in particular, $d(u_3)\geq 4$ and even).

First let $c(\D)\leq c(4,6,8)= \pi/12$. If $\D$ shares an edge with at least one boundary region then GR applies, so assume otherwise. The four cases for $d(u_1)\geq 8$, $d(u_2)=6$ are shown in Figures {{C}}(ii)--(v); and the six cases for $d(u_1)=6$, $d(u_2)\geq 8$ are shown in Figures {{C}}(vi)--(ix). (For example, in Figure {{C}}(vi) $l(u_1)=a\mu ba^{-1}\lambda b^{-1}$ or $a\mu ba^{-1}b\mu$.) In Figures {{C}}(ii),(iii),(vi),(vii) distribute $c(\D)\leq \pi/12$ to $\hD$ as shown; and in Figures {{C}}(iv),(v),(viii),(ix) distribute $c(\D)/2\leq \pi/24$ to each of $\hD_1$ and $\hD_2$, as shown.
    
Now let $c(\D)>c(4,6,8)=\pi/12$; that is $c(\D)= c(4,6,6)= \pi/6$. If $\D$ shares an edge with at least two boundary regions then GR applies, so assume otherwise. If $u_1$ or $u_2$ is a boundary vertex then $\D$ is given by Figures {{D}}(i)--(iii) and in each case distribute $c(\D)=\pi/6$ to $\hD$ as shown. So let $u_1,u_2$ be interior vertices (of degree 6), in which case $\D$ is given by Figures {{D}}(iv).

The three cases when $d(u_3)\geq 8$ are given by Figures {{D}}(v)--(vii). In Figure {{D}}(v) and (vii) distribute $c(\D)/2=\pi/12$ to each of $\hD_1$ and $\hD_2$ as shown; and in Figure {{D}}(vi) distribute $c(\D)/2=\pi /12$ to $\hD_1$ and $c(\D)/4=\pi /24$ to $\hD_3$ and $\hD_4$ as shown. 
    
Let $d(u_3)=4$, forcing $d(u_3)=d^{-1}ba^{-1}b$. If the vertex $u_4$ of Figure {{D}}(iv) is exceptional or $d(u_4)\geq 8$ then distribute $c(\D)=\pi/6$ to $\hD$ as shown in Figure {{D}}(viii). Assume otherwise, in which case $l(u_4)=b^{-1}ae^{-1}\lambda$ and we have Figure {{D}}(ix). If the vertex $v_2$ of Figure {{D}}(ix) is exceptional then distribute $c(\D)=\pi /6$ to $\hD_2$ as shown; or if $v_2$ is not exceptional but $v_1$ is exceptional then distribute $c(\D)=\pi /6$ to $\hD_1$ as shown in Figure {{D}}(ix). In Figure {{D}}(ix), when $c(\D)$ is distributed to $\hD_1$ we introduce the exceptional rule to GR that none of $c(\hD_2)$ is transferred to $\hD_1$ (in the event that $\hD_1$ is a boundary region) but instead that $ c(\hD_2)/2\leq c(4,6,8)/2=\pi/24$ is transferred to each of the other two neighbouring regions of $\hD_2$.

Assume that neither $v_1$ nor $v_2$ is exceptional.  Then, in particular, $d(v_1)\geq 6$ and, as noted, $d(v_2)\geq 8$. Consider the vertex $v_3$ of Figure {{D}}(ix). If $v_3$ is exceptional or $d(v_3)\geq 6$ then distribute $c(\D)=\pi/6$ to the region $\hD$ and introduce the exceptional rule to GR: distribute $c(\D_1)\leq \pi/12$ to $\hD_1$ as shown in Figure {{D}}(x). 

Let $d(v_3)=4$ forcing $l(v_3)=d^{-1}ba^{-1}b$. If  $d(v_2)\geq 10$ then distribute $c(\D)/2=\pi/12$ to each of the regions $\hD_1$ and $\hD_2$ and introduce the exceptional rule to GR: distribute $c(\D_1)= c(4,6,10)=\pi/30$ to $\hD_3$ as shown in Figure {{D}}(xi). This leaves the case $d(v_3)=4$ and $d(v_2)=8$. Then $l(v_2)\in \{a\mu a \mu ad^{-1}a \mu, a\mu a\mu ba^{-1}\lambda c^{-1}/e^{-1}\}$.
    
If the vertex $u$ of Figure {{D}}(xi) is exceptional then distribute $c(\D)/2=\pi/12$ to each of $\hD_1$ and $\hD_2$ as shown in Figure {{D}}(xii). If the vertex $u$ is not exceptional then distribute $c(\D)=\pi/6$ to the region $\hD$ as shown in Figures {{D}}(xiii),(xiv) according to the label of $v_2$.
    
In Figure {{D}}(xiii) we introduce an exceptional rule to GR: none of $c(\D_1)$ is distributed to $\hD$ but rather, as shown, we follow Figure {{B}}(viii) and distribute $c(\D_1)\leq \pi/4$ to $\hD_1$. In Figure {{D}}(xiv)  we introduce an exceptional rule to GR: none of $c(\D_1)$ is distributed to $\hD$ but, as shown, $c(\D_1)/2=\pi/24$ is distributed to each of the other two neighbouring regions of $\D_1$. This completes the case when $d(u_3)=4$. 
    
Thus we may assume $d(u_3)=6$ and so we are left with the case when both $u_2$ and $u_3$ are interior vertices of degree 6. (Note that, since $u_3$ is not exceptional and has a sublabel $ba^{-1}b$, it must be interior by Lemma \ref{lem:stargraphlemma}.) In this case $u_2,u_3$ are given by Figure {{D}}(xv). Suppose that vertices $u_4,u_5$ of Figure {{D}}(xv) are interior of degree 6. Then $u_4$ and $u_5$ are given by Figure {{D}}(xvi) where, if $u_6,u_7$ are interior of degree 6 then $u_6,u_7$ have the same labels as $u_4,u_5$, respectively, namely $l(u_6)=b^{-1}ab^{-1}\lambda a^{-1}\lambda$,  $l(u_7)=a^{-1}b\mu a\mu b$. So we proceed in this way until we obtain an even $j\geq 2$ such that $u_{j+2},u_{j+3}$ are not both interior of degree 6 but $u_j,u_{j+1}$ are both interior of degree 6. Since the curvature distribution will be exactly the same for each pair $u_j,u_{j+1}$ for $j\geq 6$ we can assume without any loss that $j=6$, as shown in Figure {{E}}(i). The two differences that occur when $j=2$ or $j=4$ will be described in Figures {{E}}(xi),(xii).

In Figure {{E}}(i) if $u_8$ is exceptional then distribute $c(\D)=\pi /6$ to $\hD_1$, as shown. Assume from now on that $u_8$ is not exceptional. Assume until otherwise stated that $d(u_8)\neq 4$; that is, $l(u_8)\neq db^{-1}ab^{-1}$. In Figure {{E}}(i) if now $u_9$ is exceptional then distribute $c(\D)=\pi/6$ to $\hD_2$ as shown. Assume from now on that $u_9$ is not exceptional. If the vertex $v_1$ of Figure {{E}}(i) is exceptional then distribute $c(\D)=\pi/6$ to the region $\hD_4$ as shown; or if $v_1$ is not exceptional and $v_2$ is exceptional then distribute $c(\D)=\pi/6$ to $\hD_3$ as shown. Assume from now on that neither $v_1$ nor $v_2$ is exceptional. 

Let $d(u_8)\geq 8$. Then either $d(u_9)\geq 6$ or $d(u_9)=4$ and the cases are shown in Figures {{E}}(ii)--(iv): in Figure {{E}}(ii) and {{E}}(iv) distribute $c(\D)/2=\pi/12$ to each of $\hD_1,\hD_2$ as shown; in Figure {{E}}(iii) distribute $c(\D)/2=\pi/12$ to $\hD_1$ and $c(\D)/4=\pi/24$ to $\hD_3$ and $\hD_4$ as shown. 

Let $d(u_8)=6$ and $d(u_9)\geq 8$. If $u_{10}$ is exceptional or $d(u_{10})\geq 6$ then distribute $c(\D)/2=\pi /12$ to each of $\hD_1$ and $\hD_2$ as shown in Figure {{E}}(v); if $l(u_{10})=db^{-1}ab^{-1}$ then distribute $c(\D)/2=\pi/12$ to $\hD_1$ and $c(\D)/4=\pi/24$ to each of $\hD_3,\hD_4$ as shown in Figure {{E}}(vi); or if $l(u_{10})=c^{-1}\lambda b^{-1}a$ or $e^{-1}\lambda b^{-1}a$ then distribute $c(\D)/2=\pi/12$ to each of $\hD_1$ and $\hD_2$ as shown in Figure {{E}}(vii).

If $d(u_8)=6$ and $d(u_9)=6$ then their vertex labelling forces both $u_8,u_9$ to be interior vertices, a contradiction to the definition of $j$.

Let $d(u_8)=6$ and $d(u_9)=4$. If $u_{10}$ is exceptional or $d(u_{10})\geq 8$ then distribute $c(\D)=\pi /6$ to $\hD$ as shown in Figure {{E}}(viii); or if $u_{10}$ is non-exceptional then $l(u_{10})$ has a sublabel $b^{-1}ae^{-1}$, then $d(u_{10}) < 8$ forces $l(u_{10}) = b^{-1}ae^{-1}\lambda$ and so distribute $c(\D)/2=\pi/12$ to each of $\hD_1$ and $\hD_2$ as shown in Figure {{E}}(vii). 

Assume that $d(u_8)=4$. Then $l(u_8)=db^{-1}ab^{-1}$. If the vertex $u$ in Figure {{E}}(ix) is exceptional or $d(u)\geq 8$ then distribute $c(\D)=\pi/6$ to $\hD$ as shown; so assume $u$ is not exceptional and $d(u)\leq 6$. Then $d(u)=4$ and $l(u)=ea^{-1}b\mu$, as in Figure {{E}}(x), in which case distribute $c(\D)/2=\pi/12$ to each of $\hD_1,\hD_2$ as shown in Figure {{E}}(x). 

This completes the curvature distribution of curvature for $u_j,u_{j+1}$ when $j\geq 6$.

If $j=2$ then distribute $\pi/6$ to $\hD$, as shown in Figure {{E}}(xi) (noting that $l(u)$ in Figure {{E}}(xi) differs from $l(u)$ in Figure {{E}}(ix)). If $j=4$ then (differing from Figure {{E}}(x)) distribute $c(\D)/4=\pi/24$ to each of $\hD_1$ and $\hD_2$ and $c(\D)/2=\pi/12$ to $\hD_3$, as shown in Figure {{E}}(xii).

\subsubsection*{Exceptions to curvature distribution rules}

In Figure {{D}}(v) or {{E}}(v) suppose that $\hD_1$ receives $\pi/12$ from $\D_1$ according to Figure {{C}}(iii). In these cases, the vertex $u_3$ in Figure {{D}}(v) and the vertex $u_9$ of Figure {{E}}(v) correspond to the vertex $u_1$ in Figure {{C}}(iii) and so the labels of these vertices are $ba^{-1} ba^{-1}w$. Assume further that $d(u_3)=8$ in Figure {{D}}(v) and $d(u_9)=8$ in Figure {{E}}(v). Note that if $\D_1$ of Figure {{D}}(v) or Figure {{E}}(v) shares an edge with a boundary region, then it distributes its curvature to that region, rather than to $\hD_1$, a contradiction. Therefore $\D_1$ does not share an edge with a boundary region. Applying Lemma \ref{lem:stargraphlemma}(iv) this forces $l(u_3),l(u_9) \in \{ ba^{-1} ba^{-1}\lambda a^{-1} da^{-1}, ba^{-1} ba^{-1}\lambda d^{-1} a\mu \}$  and the two cases are given by  Figure {{F}}(i),(ii)  respectively. Note that in Figure {{F}}(i) when the $a$-corner vertex of $\D_0$ has label $ab^{-1}\lambda a^{-1} b\mu$ then this corresponds to Figure {{D}}(v); when it has label $ab^{-1}\lambda a^{-1}\lambda b^{-1}$ then this corresponds to Figure {{E}}(v); and similarly for Figure {{F}}(ii). The fact that  $\hD_1$ receives $\pi/12$ from $\D_1$ forces  $l(v) = ba^{-1}ba^{-1}w$ (see Figures {{F}}(i),(ii)). Moreover, if $\D_1$ shares an edge with a boundary region, then it distributes its curvature to that region, rather than to $\hD_1$, a contradiction. Therefore $\D_1$ does not share an edge with a boundary region. Applying Lemma \ref{lem:stargraphlemma}(iv) this forces $l(v)\in \{ba^{-1}ba^{-1}\lambda a^{-1}da^{-1},  ba^{-1}ba^{-1}\lambda d^{-1}a \mu\}$. If the vertex in Figure {{F}}(i) or (ii) with label $\encircle{\geq 6}$ is exceptional then distribute all of the curvature $c(\D)=\pi/6$ to $\hD_2$ as shown by the dotted arrow in Figure {{F}}(i) and (ii). Otherwise distribute $c(\D)/2=\pi/12$ to each of $\hD_2$ and $\hD_3$, as shown by the solid arrow in Figure {{F}}(i) and (ii) except when $l(u)=d^{-1}ba^{-1}b$ in Figure {{F}}(ii), in which case distribute $c(\D)/2=\pi/12$ to $\hD_2$ and $c(\D)/4=\pi/24$ to each boundary region $\hD_4$ and $\hD_5$, as shown by the dotted arrows. The fact that no curvature is distributed from $\D$ to $\hD_1$ in Figures {{F}}(i),(ii) but instead $c(\D)$ is distributed among other regions means that all of this is an exception to the curvature distribution of Figures {{D}}(v) and {{E}}(v).

Now consider Figure {{E}}(ii) and suppose that $\hD_1$ receives $\pi/12$ from $\D_1$ according to Figure {{C}}(iii) and suppose further that $d(u_8)=8$. Then $l(u_8)=ab^{-1}ab^{-1}w$, where the last letter of $w$ is not $d^{-1}$ or $e^{-1}$ (for otherwise $\D_1$ is adjacent to a boundary region and so, by GR, does not distribute curvature to $\hD_1$), and this forces $l(u_8)\in \{ab^{-1}ab^{-1}ad^{-1}a\mu , ab^{-1}ab^{-1}\lambda a^{-1}d\mu\}$. In this case distribute $c(\D)/2=\pi/12$ to each of $\hD_2$ and $\hD_3$ as shown by the solid arrows in Figure {{F}}(iii),(iv) except when $l(u)=db^{-1}ab^{-1}$ in Figure {{F}}(iv), in which case distribute $c(\D)/2=\pi/12$ to $\hD_2$ and $c(\D)/4=\pi/24$ to each boundary region $\hD_4,\hD_5$ as shown by the dotted arrows. The fact that no curvature is distributed from $\D$ to $\hD_1$ in Figures {{F}}(iii),(iv) means that all of this is an exception to the curvature distribution of Figure {{E}}(ii).

Finally, consider Figure {{E}}(vii) and suppose that $\hD_1$ receives $\pi/12$ from $\D_1$ according to Figure {{C}}(iii) and that $d(v_1)=8$. Then $l(v_1)=a^{-1}ba^{-1}ba^{-1}w$ which forces $l(v_1)=a^{-1}ba^{-1}ba^{-1}\lambda a^{-1}d$ or $a^{-1}ba^{-1}ba^{-1}da^{-1}\lambda$. But in the latter case $\D_1$ is adjacent to a boundary region, so according to GR does not transfer curvature to $\hD_1$, a contradiction. So $l(v_1)=a^{-1}ba^{-1}ba^{-1}\lambda a^{-1}d$. In this case distribute $c(\D)/2=\pi/12$ to each of $\hD_2,\hD_3$ as shown in Figure {{F}}(v). The fact that no curvature is distributed from $\D$ to $\hD_1$ in Figure {{F}}(v) means that all of this is an exception to the curvature distribution of Figure {{E}}(vii).

\subsubsection{The symmetric cases}\label{sec:symmetricases}

We have completed the curvature distribution rules for the case when $\D$ is given by 
Figure {{C}}(i) and we turn now to the case of Figure {{G}}(i). As we explain below, the arguments are exactly the same except in regard to the labelling of degree 8 vertices in Figures {{D}}/{{H}}(xiii), {{F}}/{{J}}(i), {{F}}/{{J}}(iii), {{F}}/{{J}}(v), {{D}}/{{H}}(xiv), when $m=5$ or $15$. Except in those cases, applying the transformation
\begin{equation}
 \lambda \leftrightarrow \hat{b}, \mu\leftrightarrow b, a\leftrightarrow \hat{a}, c\leftrightarrow \hat{c}, d\leftrightarrow \hat{d}, e\leftrightarrow \hat{e}   \label{eq:transformation}
\end{equation}
to Figures {{C}}(i)--(ix), {{D}}(i)--(xvi), {{E}}(i)--(xii), {{F}}(i)--(v) respectively, yields Figures {{G}}(i)--(ix), {{H}}(i)--(xvi), {{I}}(i)--(xii), {{J}}(i)--(v).  

The degree 8 labels $l_1=a\mu a\mu ad^{-1}a\mu$, $ba^{-1}ba^{-1}\lambda a^{-1}da^{-1}$, $ab^{-1}ab^{-1}ad^{-1}a\mu$, $a^{-1}ba^{-1}ba^{-1}\lambda a^{-1}d$ of Figure {{D}}(xiii), {{F}}(i), {{F}}(iii), {{F}}(v) have corresponding labels $l_2=a^{-1}ba^{-1}ba^{-1}da^{-1}\lambda$, $\mu a \mu a \mu a d^{-1} a$, $a^{-1}\lambda a^{-1}\lambda a^{-1}da^{-1}\lambda$, $a\mu a\mu a\mu a d^{-1}$ in Figure {{H}}(xiii), {{J}}(i), {{J}}(iii), {{J}}(v), respectively, as shown and in each case $l_1$ does not transform to $l_2$ by \eqref{eq:transformation}. The degree 8 labels  $a\mu a\mu wd^{-1}/c^{-1}/e^{-1}$ of Figure {{D}}(xiv) where $wd^{-1}=a\mu a d^{-1}$, $wc^{-1}=ba^{-1}\lambda c^{-1}$, $we^{-1}=ba^{-1}\lambda e^{-1}$ have corresponding labels $a^{-1}ba^{-1}b w'd/c/e$ in Figure {{H}}(xiv), as shown, where  $w'd=a^{-1}\lambda a^{-1}d$, $w'c=\mu ab^{-1}c$, $w'e=\mu ab^{-1}e$, respectively. Note that $wc^{-1},we^{-1}$ transform to $w'c,w'e$ by \eqref{eq:transformation} (so these are not exceptions), whereas $wd^{-1}$ does not transform to $w'd$ (so this is an exception).

The reason for the remaining exceptions is as follows. In Lemma \ref{lem:stargraphlemma}(i),(ii),(iv) the exponent sum of $t$ is an even multiple of $m$, and hence is congruent to $0 \bmod 2m$, apart from in (iv) when $w_1\in \{ a\mu ad^{-1}, ad^{-1}a\mu\}$, $w_2\in \{ a^{-1}\lambda a^{-1} d, a^{-1}da^{-1}\lambda \}$, $w_3\in \{ da^{-1}\lambda a^{-1}, \lambda a^{-1}da^{-1}\}$, and \linebreak $w_4\in \{d^{-1}a\mu a, \mu ad^{-1} a\}$. Apart from theses cases in (iv), the exponent sum 0 label transforms, by \eqref{eq:transformation}, to an allowable exponent sum 0 label. However, for the labels $l_1=a\mu a\mu ad^{-1}a\mu$, $ba^{-1}ba^{-1}\lambda a^{-1}da^{-1}$, $ab^{-1}ab^{-1}ad^{-1}a\mu$, $a^{-1}ba^{-1}ba^{-1}\lambda a^{-1}d$, the exponent sum yields $m\equiv 15 \bmod 2m$, and so $m=5$ or $15$, whereas the exponent sum corresponding to their transformations (by \eqref{eq:transformation}) gives $15\equiv 0\bmod 2m$, a contradiction. For this reason we must apply Lemma \ref{lem:stargraphlemma}(iv) again as follows. Observe that the four labels for $l_1$ are of the form $a\mu a\mu w_1$, $ba^{-1}ba^{-1}w_3$, $ab^{-1}ab^{-1}w_3^{-1}$ and $a^{-1}ba^{-1}bw_2$. Applying transformation \eqref{eq:transformation} we get $a^{-1}ba^{-1}bw_2$, $\mu a\mu aw_4$, $a^{-1}\lambda a^{-1}\lambda w_4^{-1}$ and $a\mu a\mu w_1$ where the sublabels $w_2$, $w_4$, $w_4^{-1}$ and $w_1$ are chosen from the pairs listed above. However $w_2$ is not $a^{-1}\lambda a^{-1}d$ as this is already accounted for in Figure {{H}}(xiv) as described in the previous paragraph, so $w_2=a^{-1}da^{-1}\lambda$; $w_4$ cannot be $d^{-1}a\mu a$ as this would contradict the assumption that $\Delta_1$ distributes $\pi/12$ to $\hD_1$ in Figure {{J}}(i), so $w_4=\mu ad^{-1}a$; $w_4^{-1}$ cannot be $(d^{-1}a\mu a)^{-1}$ as this would contradict the assumption that $\Delta_1$ distributes $\pi/12$ to $\hD_1$ in Figure {{J}}(iii), so $w_4^{-1}=(\mu ad^{-1}a)^{-1}$; and $w_1$ cannot be $ad^{-1}a\mu$ as this would contradict the assumption that $\Delta_1$ distributes $\pi/12$ to $\hD_1$ in Figure {{J}}(v), so $w_1=a\mu ad^{-1}$. The four resulting labels are those listed for $l_2$ above and hence we obtain the Figures {{H}}(xiii), J(i), J(iii) and J(v).

\subsection{Implications of Stage I distribution}\label{sec:StageIconsequences}

The following properties hold in Figures {{B}}--{{J}}. It is routine to confirm properties (P1),(P3) by inspection; the confirmation of the remaining properties is lengthier and more involved. We provide some comments regarding this, but for reasons of space we omit the details.
\begin{itemize}
    \item[(P1)] If $\hD$ is an interior region that receives curvature across its edge $e$ then neither of the endpoints of $e$ has label $(b\mu b\mu)^{\pm 1}$.

    \item[(P2)] If $\hD$ is a region that receives curvature across its edge $e$ then neither of the endpoints of $e$ is an exceptional vertex.
    
    \item[(P3)] If curvature is distributed from any given  region to a non-adjacent region  then (in its travels) it never crosses an edge that has a vertex with label $(b\mu b\mu)^{\pm 1}$, except possibly when the non-adjacent region is $\hD$ of Figure {{D}}(xiii) or {{H}}(xiii).
    
    \item[(P4)] No interior region receives curvature from more than one non-adjacent region.
    
    \item[(P5)] No boundary region $\hD$ receives curvature from two non-adjacent regions across the same edge. 
\end{itemize}

Note that in (P5), it is possible that a boundary region $\hD$ simultaneously receives curvature both from an adjacent and a non-adjacent region across the same edge. This will be discussed later, in the proof of Claim \ref{claim:two}. 

\begin{proof}[Confirming Property (P2)]
By applying the transformation \eqref{eq:transformation}, we need only check Figures {{B}},{{C}},{{D}},{{E}} and {{F}}, and the exceptional cases described in Section \ref{sec:symmetricases}. Since we assume that $\D$ does not contain an exceptional vertex and, moreover, that the region $\D_1$ of Figure {{B}}(ix), (x), (xii), and (xii) does not contain an exceptional vertex, we need only check regions that receive curvature from non-adjacent regions. That is, we need only check 
Figures  {{C}}(iv), (v), (viii), (ix), {{D}}(i), (v)--(xiv), {{E}}(i)--(xii), {{F}}(i)--(v), {{H}}(xiii), {{H}}(xiv), {{J}}(i), {{J}}(iii), {{J}}(v). In some cases (P2) is readily verified by inspecting vertex labels of the figures; for example in Figure {{B}}(x). In other cases it is necessary to refer to the text; for example the description of the curvature distribution for Figure {{E}}(i) establishes that none of the vertices $v_1, v_2, u_8$ or $u_9$ are exceptional in the subsequent Figures {{E}}(ii)--(xii).
\end{proof}

\begin{proof}[Confirming Property (P4)]
It suffices to confirm (P4) up to inversion, and thus only positive regions need be considered. There are 54 positive regions to consider: 4 {{D}} regions; 16 {{E}} regions; 7 {{F}} regions; 4 {{H}} regions; 16 {{I}} regions; and 7 {{J}} regions. The {{D}},{{E}},{{F}} regions are the following 27 (where  \textoverline{X} denotes the inverse of Figure X): {{D}}(v) $\hD_2$; \textoverline{{{D}}}(ix) $\hD_2$; {{D}}(ix) $\hD_1$; \textoverline{{{D}}}(xii) $\hD_2$; {{E}}(i) $\hD_1$; \textoverline{{{E}}}(i) $\hD_2$; \textoverline{{{E}}}(i) $\hD_2$; \textoverline{{{E}}}(i) $\hD_3$; \textoverline{{{E}}}(i) $\hD_4$; {{E}}(ii) $\hD_1$; \textoverline{{{E}}}(ii) $\hD_2$; {{E}}(iii) $\hD_1$; {{E}}(iv) $\hD_1$; \textoverline{{{E}}}(v) $\hD_1$; {{E}}(v) $\hD_2$; \textoverline{{{E}}}(vi) $\hD_1$; {{E}}(vii) $\hD_1$;  \textoverline{{{E}}}(vii) $\hD_2$; {{E}}(x) $\hD_1$; \textoverline{{{E}}}(x) $\hD_2$; \textoverline{{{E}}}(xii) $\hD_3$; {{F}}(i) $\hD_2$; {{F}}(ii) $\hD_2$;  \textoverline{{{F}}}(ii) $\hD_3$; \textoverline{{{F}}}(iii) $\hD_2$; \textoverline{{{F}}}(iv) $\hD_2$; {{F}}(iv) $\hD_3$; and {{F}}(v) $\hD_2$. There are 27 further regions from {{H}},{{I}},{{J}} that are obtained from the above list by applying the symmetry \eqref{eq:transformation}.

The comparisons required are as follows: {{D}}(v) $\hD_2$ with (the remaining) 53 regions; \textoverline{{D}}(ix) $\hD_2$ with 52 regions; $\dots$ ; {{F}}(iv) $\hD_3$ with 28 regions; and {{F}}(v) $\hD_2$ with 27 regions (that is, with all the {{H}},{{I}},{{J}} regions referred to above). This gives a total of 1080 comparisons. Note that no further comparisons are needed, as any coincidences among the {{H}},{{I}},{{J}} regions would, by symmetry, have appeared among the {{D}},{{E}},{{F}} regions. Note further that 9 of the {{E}} regions and 6 of the {{F}} regions represent an infinite family of regions (each at a distance $j\geq 6$ from $\D$) and so another 15 comparisons (for example {{E}}(ii) $\hD_1$ with {{E}}(ii) $\hD_1$ for different values of $j$) must be made, bringing the total to 1095.

The first three initial checks are as follows:

\begin{itemize}
    \item Compare the degrees of the vertices. For example, {{D}}(v) $\hD_2$ can be distinguished from \textoverline{{I}}(i) $\hD_1$ by comparing the degrees of the $b$-vertices.

    \item Compare information about whether vertices are exceptional or are non-exceptional. For example, in Figure \textoverline{{{F}}}(iv) $\hD_2$ the $b$-vertex is not exceptional whereas in Figure {{H}}(ix) $\hD_2$  the $b$-vertex is exceptional (in the event that curvature is distributed from $\D$ to $\hD_2$), which distinguishes these regions.
    
    \item Compare the corner labels. For example, comparing the labelling of the $a$-vertices distinguishes {{E}}(ii) $\hD_1$ from \textoverline{II}(ii) $\hD_1$.
    \end{itemize}

These tests distinguish a total of 973 pairs, leaving 122 pairs of regions to compare. We then apply Property (P3), noting that the region $\hD$ in {{D}}(xiii) and {{H}}(xiii) is a boundary region. For example, this shows that {{E}}(v) $\hD_2$ cannot coincide with {{F}}(i) $\hD_2$ (which we recall corresponds either to {{D}}(v) or {{E}}(v)), since this would force the crossing of an edge, one of whose endpoints having label $(b\mu b\mu)^{\pm 1}$. This (P3) test distinguishes a further 57 pairs leaving 65 pairs of regions yet to be considered.
    
The argument that distinguishes the regions in each of the remaining pairs is the same. As an illustration, consider \textoverline{{{E}}}(ii) $\hD_2$ and {{I}}(vii) $\hD_2$. If these two regions were to coincide then {{I}}(vii) $\hD_2$ would have to receive the curvature from $\D$ as in \textoverline{{{E}}}(ii) from directly above. But, as can be seen from {{I}}(vii), there are no interior vertices of degree 4 between $\hD_2$ and the boundary that would allow for such a $\D$, a contradiction.
\end{proof}

\begin{proof}[Confirming Property (P5)]
As for (P4), it suffices to consider positive regions. There is a total of 96 positive regions to be considered. There are 25 regions (such as {{C}}(iv) $\hD_1$), of \emph{Type 1}, say, in which the final edge crossed is from a $(b^{-1},a^{-1})$-edge (relative to the penultimate region) to a $(b,\lambda)$-region (relative to $\hD$) and we denote this crossing by $(b^{-1},a^{-1})\rightarrow (b,\lambda)$; 
there are 25 regions, of \emph{Type 2}, with crossing of type $(a^{-1},\mu)\rightarrow (b,\lambda)$; 
there are 18 regions, of \emph{Type 3}, with crossing of type $(b^{-1},a^{-1})\rightarrow (\lambda,a)$; 
there are 18 regions, of \emph{Type 4}, with crossing of type $(a^{-1},\mu)\rightarrow (a,b)$; 
there are 5 regions, of \emph{Type 5}, with crossing of type $(\mu,b^{-1})\rightarrow (\lambda,a)$; 
and there are 5 regions, of \emph{Type 6}, with crossing of type $(\mu,b^{-1})\rightarrow (a,b)$. Observe that the transformation \eqref{eq:transformation} transforms regions of Type 1 to Type 2;  regions of Type 3 to Type 4; and regions of Type 5 to Type 6. Therefore, it is only necessary to confirm that there are no coincidences among the 25 regions of Type 1, or among the 18 regions of Type 3, or among the 5 regions of Type 5. (Note that there is no need to compare regions of different type, since Property (P5) concerns curvature being transferred across a single edge.)

The (positive) Type 1 regions are: Figure {{C}}(iv) $\hD_1$; {{C}}(v) $\hD_2$; {{C}}(viii) $\hD_1$; {{D}}(v) $\hD_2$; \textoverline{{{D}}}(vi) $\hD_4$; {{D}}(vii) $\hD_2$; {{D}}(viii) $\hD$; {{D}}(xi) $\hD_1$; \textoverline{{{E}}}(i) $\hD_2$; \textoverline{{{E}}}(ii) $\hD_2$; {{E}}(iii) $\hD_4$; \textoverline{{{E}}}(iv) $\hD_2$; {{E}}(v) $\hD_2$; \textoverline{{{E}}}(vi) $\hD_4$; {{E}}(viii) $\hD$;
\textoverline{{{F}}}(ii) $\hD_3$; {{F}}(ii) $\hD_4$; \textoverline{{{F}}}(iii) $\hD_2$; \textoverline{{{F}}}(iv) $\hD_2$; {{F}}(iv) $\hD_3$; {{F}}(iv) $\hD_5$; \textoverline{{{G}}}(ix) $\hD_2$; \textoverline{{{H}}}(i) $\hD$; \textoverline{{{I}}}(vii) $\hD_1$; and \textoverline{{{I}}}(x) $\hD_1$. (Note that the regions {{F}}(ii) $\hD_4$ and {{F}}(iv) $\hD_5$ only receive curvature if the vertex $u$ has label $d^{-1}ba^{-1}b$.)
The (positive) Type 3 regions are: Figure {{C}}(ix) $\hD_1$; \textoverline{{{D}}}(xii) $\hD_2$; \textoverline{{{E}}}(xii) $\hD_2$; \textoverline{{{F}}}(i) $\hD_3$; {{F}}(iii) $\hD_3$; \textoverline{{{F}}}(v) $\hD_3$; \textoverline{{{G}}}(iv) $\hD_2$; \textoverline{{{G}}}(viii) $\hD_2$; {{H}}(vi) $\hD_3$; {{H}}(ix) $\hD_2$; \textoverline{{{H}}}(xiv) $\hD$; \textoverline{{{I}}}(iii) $\hD_3$; {{I}}(vi) $\hD_3$; {{I}}(ix) $\hD$; {{I}}(xi) $\hD$; {{I}}(xii) $\hD_2$;  
\textoverline{{{J}}}(ii) $\hD_5$; \textoverline{{{J}}}(iv) $\hD_4$. 
The (positive) Type 5 regions are: {{D}}(x) $\hD$; {{D}}(xi) $\hD_2$; {{D}}(xii) $\hD_1$; {{D}}(xiii) $\hD$; and \textoverline{{{G}}}(v) $\hD_1$. 

The confirmation that no pair of regions from Type 1, or pair from Type 3, or pair from Type 5 coincide uses the same tests as those for confirming Property (P4), so we omit the details.
\end{proof}

We now prove some other consequences of the Stage I curvature distribution.

\begin{notation*}
For an interior region $\hD$ let $c^*(\hD)$ denote $c(\hD)$ plus all the curvature $\hD$ receives in Stage I minus all the curvature $\hD$ distributes in Stage I. For any interior or boundary region $\hD$ let $\hc(\hD)$ denote the sum $c(\hD)$ plus any curvature $\hD$ received in Stage I.
\end{notation*}

\begin{claim}\label{claim:one}
If $\hD$ is an interior region with no exceptional vertices then $c^*(\hD)\leq 0$. 
\end{claim}

\begin{proof}[Proof of Claim \ref{claim:one}]
    If $\hD$ receives curvature from only a single adjacent or non-adjacent region then inspecting Figures {{C}}--{{J}} shows that $c^*(\hD)\leq 0$. Thus, by Property (P4) we may assume either Case 1: $\hD$ receives curvature from more that one adjacent region, but no non-adjacent regions; or Case 2: $\hD$ receives curvature from a single non-adjacent region and at least one adjacent region.

    \smallskip

\noindent {Case 1.} Suppose that $\hD$ receives curvature from more than one adjacent region, but no non-adjacent regions. In particular, $\hD$ receives curvature across an edge other than its $(b,\lambda)$-edge. We see from Figures {{C}}--{{J}} that the region $\hD$ receives from an adjacent region across other than its $(b,\lambda)$-edge precisely in Figure {{C}}(iii) ($(a,\lambda)$-edge) and its copy $\hD_1$ in Figures {{F}}(i)--(v); {{D}}(xi), $\hD=\hD_3$ ($(a,b)$-edge);  Figure {{G}}(iii) ($(a,b)$-edge) and its copy $\hD_1$ in Figures {{J}}(i)--(v); and {{H}}(xi), $\hD=\hD_3$ ($(a,\lambda )$-edge).

In Figures {{D}}(xi) and {{H}}(xi), $\hD=\hD_3$ can, in addition, only receive curvature from the region $\D_2$ shown. If $c(\D_2)\leq \pi/12$ then $c^*(\hD_3)\leq c(6,6,10)+\pi/12+\pi/30<0$. On the other hand, if $c(\D_2)=c(4,6,6)=\pi/6$ then    according to Figures {{D}}(v)--(vii) or Figures {{H}}(v)--(vii), $c(\D_2)/2=\pi/12$ is transferred to $\hD$ and  $c^*(\hD_3)\leq c(6,6,10)+\pi/12+\pi/30<0$.

If $\hD$ is given by Figure {{C}}(iii) then (since it cannot coincide with $\hD$ of Figure {{G}}(iii)), $\hD$ can, in addition, only receive curvature across its $(b,\lambda)$-edge and if it does so, then inspecting Figures {{B}}--{{J}} shows that the following are the possibilities: $\hD=\hD$ of Figure {{C}}(vi); or $\hD=\hD_1$ of Figures {{D}}(v),(vi) or (vii).  For Figure {{C}}(vi) we see that $c^*(\hD) \leq c(6,8,8)+2\pi/12=0$. Consider Figures {{D}}(v)--(vii). In each case, if $\hD_1$ coincides with $\hD$ of Figure {{C}}(iii) then this forces $l(u_3)=ba^{-1}ba^{-1}w$ (where $u_3$ refers to its position in Figures {{D}}(v)--(vii), rather than in Figure {{C}}(iii)). In Figure {{D}}(vii) this forces $d(u_3)\geq 10$ (see Lemma \ref{lem:stargraphlemma}(iv)) and $c^*(\hD)\leq c(6,6,10)+\pi/12+c(\D)\leq -2\pi/15+\pi/12+\pi/30<0$. For Figure {{D}}(vi) if $d(u_3)\geq 10$ then again $c^*(\hD_1)\leq 0$; but if $d(u_3)=8$ then $l(u_3) \in \{ba^{-1}ba^{-1}cb^{-1}a\mu, \linebreak ba^{-1}ba^{-1}eb^{-1}a\mu\}$, in which case, according to GR, $c(\D_1)=\pi/12$ is distributed to an adjacent boundary region and not to $\hD_1$, a contradiction (to the assumption that $\hD$ receives curvature from more than one adjacent region). For Figure {{D}}(v) if $d(u_3)\geq 10$ then $c^*(\hD_1)\leq 0$ and if $d(u_3)=8$ then $l(u_3)\in \{ ba^{-1}ba^{-1}\lambda a^{-1}da^{-1},$ $ba^{-1}ba^{-1}\lambda d^{-1}a\mu\}$. But then $\hD_1$ does not receive any curvature from $\D$ in Figure {{D}}(v) since Figures {{F}}(i) and (ii) now apply. A similar argument shows that if $\hD$ is given by Figure {{G}}(iii) and receives curvature across its $(a,b)$-edge then $\hD=\hD$ of Figure {{G}}(vi) or $\hD=\hD_1$ of Figure {{H}}(v),(vi),(vii), and again $c^*(\hD)\leq 0$ (where Figures {{J}}(i),(ii) are used in place of Figures {{F}}(i),(ii)).   

\smallskip
    
\noindent {Case 2.} Suppose that $\hD$ receives curvature from a single non-adjacent region and at least one adjacent region. Then since $\hD$ receives from a non-adjacent region, Figures {{B}}--{{J}} show that it must be across the $(b,\lambda)$-edge; moreover, $\hD$ cannot in addition receive from the adjacent region sharing the $(b,\lambda)$-edge. It follows that $\hD$ must be given by either Figure {{C}}(iii) or Figure {{G}}(iii) and that $\hD$ receives from a single adjacent region. 
    
Let $\hD$ be given by Figure {{C}}(iii). Then $\hD$ must coincide with one of the following regions: Figure {{E}}(ii) $\hD_1$; {{E}}(ii) $\hD_2$; {{E}}(iii) $\hD_1$; {{E}}(iv) $\hD_1$; {{E}}(v) $\hD_1$; {{E}}(v) $\hD_2$; {{E}}(vi) $\hD_1$; {{I}}(vii) $\hD_1$; {{I}}(x) $\hD_1$; {{F}}(ii) $\hD_2$; {{F}}(iii) $\hD_2$; {{F}}(iv) $\hD_2$.
    
If $\hD$ is given by Figure {{E}}(ii) $\hD_2$, {{E}}(v) $\hD_2$, {{F}}(ii) $\hD_2$, {{F}}(iii) $\hD_2$ or {{F}}(iv) $\hD_2$ then $c^*(\hD)\leq c(6,8,8)+2\pi/12 =0$. If $\hD$ is given by Figure {{E}}(iii) $\hD_1$ then $l(u_8)=ab^{-1}ab^{-1}\lambda w$; or by Figure {{E}}(iv) $\hD_1$ then $l(u_8)=ab^{-1}ab^{-1}c/d w$ and in both cases $d(u_8)\geq 10$. Or if $\hD$ is given by Figure {{E}}(vi) $\hD_1$ then $l(u_9)=\mu ba^{-1}ba^{-1}w$ and this forces $d(u_9)\geq 10$. In all three cases $c^*(\hD)\leq c(6,6,10)+\pi/12+\pi/30<0$. This leaves Figures {{E}}(ii) $\hD_1$, {{E}}(v) $\hD_1$, {{I}}(vii) $\hD_1$, and {{I}}(x) $\hD_1$. If $d(u_8)\geq 10$ in Figure {{E}}(ii) or $d(u_9)\geq 10$ in Figure {{E}}(v) or $d(v_1)\geq 10$ in Figure {{I}}(vii) or {{I}}(x) then $c^*(\hD)\leq c(6,6,10)+\pi/12+\pi/30<0$. On the other hand, if $d(u_8)=8$, $d(u_9)=8$, or $d(v_1)=8$ then given $\hD=\hD_1$ receives from an adjacent region,  it does not receive from a non-adjacent region since instead we apply the exceptional rules shown in Figure {{F}}(iii) or (iv), Figure {{F}}(i) or (ii), or Figure {{F}}(v) (respectively) and in these figures $c^*(\hD_1)\leq c(6,6,8)+\pi/12=0$.

We have shown that if $\hD$ is given by Figure {{C}}(iii) then $c^*(\hD)\leq 0$. If $\hD$ is given by Figure {{G}}(iii) then we can apply symmetry and a similar argument again shows that $c^*(\hD)\leq 0$.
\end{proof}

\begin{claim}\label{claim:two}
If $\hD$ is a boundary region, none of whose vertices are exceptional, then either $\hc(\hD)\leq \pi/2$ or ($\pi/2<\hc(\hD)\leq 7\pi /12$ and one of the boundary vertices has degree at least 8).
\end{claim}

\begin{proof}[Proof of Claim \ref{claim:two}]
First observe that if all three vertices of $\hD$ have degree 4 then inspecting Figures {{B}}--{{J}} shows that $\hD$ does not receive any curvature and so $\hc (\hD)=c(\hD)=\pi/2$. Assume then that $c(\hD)\leq c(4,4,6)=\pi/3$. Note that if the maximum total curvature $\hD$ receives across an edge is $\pi/12$ then $\hc(\hD)\leq \pi/3+2\pi/12=\pi/2$. Thus we may assume that $\hD$ receives more than $\pi/12$ across at least one of its edges.

Suppose first that $\hD$ receives curvature from adjacent regions only. In Figures {{B}}(iv), (vi) with $\hD=\hD_2,\hD_1$ (respectively), we have $\hc(\hD)\leq c(\hD)+\pi/4+c(4,6,8)=c(4,4,8)+\pi/4+c(4,6,8)=7\pi/12$, but note that, by Lemma \ref{lem:stargraphlemma}, $\hD$ has a boundary vertex of degree at least 8 in both cases; in {{B}}(v),(vii) with $\hD=\hD_2,\hD_1$ (respectively) we have 
 $\hc(\hD)\leq c(4,4,8)+\pi/6+ c(4,4,8)/3=\pi/2$; in {{B}}(viii)--(x) with $\hD=\hD_2$ or {{B}}(xi)--(xiii) with $\hD=\hD_1$ we see that $\hc(\hD) \leq c(4,4,8) + \pi/4+0=\pi/2$; and in {{B}}(xiv),(xv) with $\hD=\hD_1,\hD_2$ (respectively), $\hc(\hD)\leq c(4,6,8)+\pi/3+c(4,6,8)=\pi/2$.

For Figures {{C}}--{{J}} the maximum amount of curvature distributed across an edge is $\pi/6$ and if any other amount is distributed, then that amount is at most $\pi/12$ so we need only consider the cases when $\pi/6$ is involved. These are Figures {{D}}(ii), {{D}}(iii), {{H}}(ii), {{H}}(iii). In Figures {{D}}(ii) and {{H}}(ii) we have $\hc(\hD)\leq c(4,6,6)+\pi/6+\pi/6=\pi/2$. In Figures {{D}}(iii) and {{H}}(iii) if $d(u)>4$ then $\hc(\hD)\leq c(4,6,6)+\pi/6+\pi/6=\pi/2$; or if $d(u)=4$ then, according to Figures {{B}}(xiv) and {{B}}(xv), $\hc(\hD)\leq c(4,4,6)+\pi/6+0=\pi/2$.

From now on, suppose that $\hD$ receives curvature from at least one non-adjacent region. This assumption together with Property (P3) implies that we need not consider the regions $\hD_2$, $\hD_1$ of Figure {{B}}(iv), {{B}}(vi) (respectively). (In what follows we will use Property (P3)  often without explicit mention.) Given Property (P5) the region $\hD$ can receive curvature across an edge from an adjacent region or a non-adjacent region  or possibly both (across the same edge).

\smallskip

\noindent {Case 1.} Suppose $\hD$ receives more than $\pi/6$ across a single edge. The following lists all cases when this can happen. The first amount is from an adjacent region and the second from a non-adjacent region and where the question marks represent a yet to be determined non-negative amount of curvature:
Figure {{B}}(viii), $\hD_2$, $\pi/4+?\geq \pi/4$;
Figure {{B}}(xi), $\hD_1$, $\pi/4+?\geq \pi/4$;
Figure {{B}}(xiv), $\hD_1$, $\pi/3+?\geq \pi/3$;
Figure {{B}}(xv), $\hD_2$, $\pi/3+?\geq \pi/3$;
Figure {{D}}(viii), $\hD$, $\pi/12+\pi/6$;
Figure {{E}}(viii), $\hD$, $\pi/12+\pi/6$;
Figure {{E}}(ix), $\hD$, $\pi/12+\pi/6$;
Figure {{E}}(xi), $\hD$, $\pi/12+\pi/6$;
Figure {{H}}(viii), $\hD$, $\pi/12+\pi/6$;
Figure {{I}}(viii), $\hD$, $\pi/12+\pi/6$;
Figure {{I}}(ix), $\hD$, $\pi/12+\pi/6$; or
Figure {{I}}(xi), $\hD$, $\pi/12+\pi/6$.

Let $\hD$ be $\hD_2$ of Figure {{B}}(viii). According to {{B}}(ix) and {{B}}(x) the region $\hD_2$ of {{B}}(viii) does not receive any curvature from the adjacent region $\D_1$ and (by Property (P3)) we see that $\hD_2$ does not receive from a non-adjacent region across its $(b^{-1},a^{-1})$-edge. If the vertex marked $\geq 8$ of $\hD_2$ has degree exactly 8 then inspecting the figures shows that $\hD_2$ does not receive from a non-adjacent region across its $(a^{-1},\mu)$-edge, and so we can assume that this third vertex has degree $\geq 10$ (since we are assuming that curvature is being received from a non-adjacent region). In this case, if $\hD_2$ receives from a non-adjacent region then it must coincide with $\hD_2$ of Figure {{D}}(xi). Thus we get
\begin{alignat*}{1}
\hc(\hD)=\hc(\hD_2) &\leq c(\hD_2)+c(\D)+\pi/12\\
                    &\leq c(4,4,10)+c(4,4,10)+\pi/12\\
                    &=29\pi/60<\pi/2
\end{alignat*}
(where the $\pi/12$ is the curvature distributed from a non-adjacent region, as in Figure {{D}}(xi)).

A similar argument holds for $\hD_1$ of {{B}}(xi) with {{B}}(xii),(xiii), {{H}}(xi) playing the roles of {{B}}(ix),(x), {{D}}(xi) respectively. That is, if $\hD$ is $\hD_1$ of Figure {{B}} (xi) then $\hc(\hD)<\pi/2$.

For the remaining 10 cases of $\hD$ in the above list  we note that no two can coincide. This means that $\hD$ receives at most $\pi/6$ across its other edge. Thus if $\hD$ is $\hD_1$ of Figure {{B}}(xiv) or $\hD_2$ of Figure {{B}}(xv) then (noting that the vertex of $\hD$ that is not degree 4 has exterior label $d$ and so has degree at least 8, by Lemma \ref{lem:stargraphlemma}(ii)) $\hc(\hD)\leq c(4,6,8)+\pi/3+\pi/6=7\pi/12$, but $\hD$ has a boundary vertex of degree at least 8.
In all 8 remaining cases observe that $c(\hD)\leq c(4,6,8)=\pi/12$ and so $\hc(\hD)\leq \pi/12+(\pi/12+\pi/6)+\pi/6=\pi/2$.

\smallskip

\noindent {Case 2.} Suppose $\hD$ receives at most $\pi/6$ across a single edge. Recall that $\hD$ receives from a non-adjacent region. If $c(\hD)\leq c(4,6,6)$ then $\hc(\hD) \leq c(4,6,6)+\pi/6+\pi/6=\pi/2$, so we can assume otherwise; that is, that $\hD$ has two vertices of degree 4. Recall that $\hD$ must receive more than $\pi/12$ across at least one edge. It follows that $\hD$ is one of the following regions: Figure {{D}}(xi) $\hD_1$; {{D}}(xii) $\hD_1$; {{H}}(xi) $\hD_1$; {{H}}(xii) $\hD_1$. 

In Figures {{D}}(xi) and {{H}}(xi) the exceptional rule means that $\hD_1$ does not receive curvature from the adjacent region $\D_1$ and, moreover, inspecting the figures shows that $\hD_1$ does not receive curvature across its $(a,\lambda)$-edge from a non-adjacent region. Therefore $\hc(\hD_1)\leq c(4,4,6)+\pi/12+c(4,6,6)/2=\pi /2$.
In Figures {{D}}(xii) and {{H}}(xii) $\hD_1$ does not receive curvature from the adjacent (boundary) region $\hD_2$ and, moreover, inspecting the figures shows that $\hD_1$ does not receive curvature across its $(a,b)$-edge from a non-adjacent region. Therefore $\hc(\hD_1)\leq c(4,4,8)+\pi/12+c(4,6,8)/2<\pi/2$. 
\end{proof}

\begin{claim}\label{claim:three}
Let $\hD$ be an interior region with exactly one vertex  $u^*$ of $\hD$ an exceptional vertex as shown in Figure {{K}}(i) and let $d(u^*)=k$. Then $\hc(\hD)\leq 2\pi/k-\pi/6$, except possibly when $\hD$ is given by $\hD$ of Figure {{K}}(ii)--(v) in which case $\hc(\hD)\leq 2\pi/k$.
\end{claim}

\begin{proof}[Proof of Claim \ref{claim:three}]
By (P2), $\hD$ can receive only across the edge $e$, say, with endpoints $u_1$ and $u_2$. If $u_1$ is interior of degree 4 and $d(u_2)\geq 6$ or $u_2$ is interior of degree 4 and $d(u_1)\geq 6$ then by (P1) we have $\hc(\hD)=c(\hD)\leq c(k,4,6)=2\pi/k-\pi/6$, so assume otherwise.

We claim that the maximum amount of curvature that $\hD$ can receive across the edge $e$ is $\pi/6$. Checking Figures {{C}}--{{J}} confirms this if $\hD$ receives only from an adjacent region. Suppose that $\hD$ receives curvature from at least one non-adjacent region, and possibly also from an adjacent region. Then $\hD$ is one of the following regions: Figure {{D}}(v), {{H}}(v) $\hD_2$; {{D}}(ix), {{H}}(ix) $\hD_1$ or $\hD_2$; {{D}}(xii), {{H}}(xii) $\hD_2$; {{E}}(i), {{I}}(i) $\hD_i$ ($1\leq i\leq 4$); {{E}}(v), {{I}}(v) $\hD_2$; {{F}}(i), {{J}}(i) $\hD_2$; {{F}}(ii), {{J}}(ii) $\hD_3$; or {{F}}(iv), {{J}}(iv) $\hD_3$. By (P4) none of these regions coincide and, moreover, in each case $\hD$ does not receive any curvature from the corresponding adjacent region and again $\pi/6$ is the maximum. If $d(u_1)\geq 6$ and $d(u_2)\geq 6$ then $\hc(\hD)\leq c(k,6,6)+\pi/6 =2\pi/k-\pi/6$. Noting that no region can have two interior vertices of degree 4, this leaves  the case when at least one of $u_1$ or $u_2$ is a boundary vertex of degree 4. Suppose that $u_1$ is such a vertex. Then the two cases are shown in Figure {{K}}(ii),(iii) in which the vertex $u_3$ is not exceptional, is exceptional (respectively). If $d(u_2)\geq 6$ then $\hc(\hD)\leq c(k,4,6)+\pi/6 =2\pi/k$; if $d(u_2)=4$ and is interior then by (P1) we have $\hc(\hD)=c(\hD)=c(k,4,4)=2\pi/k$; or if $u_2$ is a boundary vertex and $d(u_2)=4$ then the adjacent region sharing the edge $e$ is forced to be a boundary region, so again $\hc(\hD)=c(\hD)=c(k,4,4)=2\pi/k$. If, instead, $u_2$ is a boundary vertex of degree 4  then the two cases are shown in Figures {{K}}(iv),(v), in which the vertex $u_3$ is not exceptional, is exceptional (respectively); and a similar argument shows that $\hc(\hD)\leq 2\pi/k$.
\end{proof}

\begin{claim}\label{claim:four}
Let $\hD$ be the boundary region shown in Figure {{K}}(vi) in which $u^*$ is an exceptional vertex and $d(u^*)>2$. Then the maximum amount of curvature $\hD$ can receive across its edge having endpoints $v_1$ and $v_2$ is $\pi/3$.
\end{claim}

\begin{proof}[Proof of Claim \ref{claim:four}]
If $\hD$ receives more than $\pi/6$ from either an adjacent or a non-adjacent region then $\hD$ is one of {{B}}(iv) $\hD_2$, {{B}}(vi) $\hD_1$, {{B}}(viii) $\hD_2$, {{B}}(x) $\hD$, {{B}}(xi) $\hD_1$, {{B}}(xii) $\hD_1$, {{B}}(xiii) $\hD$, {{B}}(xiv) $\hD_1$, {{B}}(xv) $\hD_2$, in which case $\hD$ receives at most $\pi/3$, and in these cases inspecting the figures shows that $\hD$ does not receive any curvature from a non-adjacent region, and so the claim holds for $\hD$. Thus we may assume that the maximum amount of curvature that $\hD$ can receive from an adjacent region or non-adjacent region is $\pi/6$ and it follows from Property (P5) that the most $\hD$ can receive is $\pi/6+\pi/6=\pi/3$, as required.
\end{proof}

\subsection{Stage II: $\hD$ either is interior and contains at least one exceptional vertex or $\hD$ is a boundary region}\label{sec:StageII}

We now define the curvature distribution scheme for Stage II. Let $\hc (\hD)>0$ and assume until otherwise stated that $\hD$ is not given by $\hD$ of Figure {{K}}(ii)--(v) nor is $\hD$ given by $\hD_1$ of Figure {{K}}(ii)--(v) and so $\hD$ is given in one of Figures {{B}}--{{J}}, so Property (P2) holds for $\hD$. Further, we recall that $\D$ has at most two exceptional vertices. Suppose that $\hD$ has no exceptional vertices and hence is a boundary region. Distribute $\hc (\hD)$ from $\hD$ to $\D^*$ as shown in Figure {{K}}(vii). Suppose then that $\hD$ has exactly one exceptional vertex $u^*$, say. Distribute $\hc (\hD)$ from $\hD$ to $\D^*$ through $u^*$ as shown in Figure {{K}}(viii). An exception to this rule is if the vertex $u$ of Figure {{K}}(vii) is an exceptional vertex. In this case we still distribute $\hc (\hD)$ to $\D^*$ across the shared boundary edge, as stated before (as in {{K}}(vii)). Suppose that $\hD$ has two exceptional vertices $u^*$ and $v^*$, say, with $d(u^*)=k_1$ and $d(v^*)=k_2$. Then from (P2) $\hD$ does not receive any curvature across any of its edges so $\hc(\hD)=c(\hD)\leq c(k_1,k_2,4)=2\pi/k_1+2\pi/k_2 -\pi/2$. Distribute, when positive, $2\pi/k_1-\pi/4, 2\pi/k_2-\pi/4$ from $\hD$ to $\D^*$  through $u^*,v^*$ respectively, as shown in Figure {{K}}(viii).

Now consider Figures {{K}}(ii)--(v) where $d(u^*)=k$.

\begin{itemize}
    \item If $\hD$ is $\hD$ of Figure {{K}}(ii)--(v) then by Claim \ref{claim:three} we have $\hc(\hD)\leq 2\pi/k$. Distribute $\pi/6$ from $\hD$ to $\hD_1$ and, if positive, distribute the remaining $\hc(\hD) -\pi/6 \leq 2\pi/k-\pi/6$ to $\hD^*$ through the exceptional vertex $u^*$, as shown. 
    
    \item Let $\hD$ be given by $\hD_1$ of Figure {{K}}(ii) or (iv) (in which it is assumed  -- see the proof of Claim \ref{claim:three} -- that the vertex $u_3$ is not exceptional). In Stage I no curvature is distributed to regions with an exceptional vertex so $\hc(\hD_1)=c(\hD_1)$. But from the two figures we see that $\hD_1$ can receive two amounts of $\pi/6$ in Stage II, and so there is at most $c(\hD_1)+2\pi/6\leq  c(k,4,4)+2\pi/6=2\pi/k+\pi/3$ to distribute. Distribute $\pi/2$ from $\hD_1$ to $\D^*$ across the boundary edge as shown; and distribute, if positive, the remaining $2\pi/k-\pi/6$ to $\D^*$ through $u^*$, as shown.

    \item Let $\hD$ be given by $\hD_1$ of Figure {{K}}(iii) or (v) (in which it is assumed  -- see the proof of Claim \ref{claim:three} -- that the vertex $u_3$ is an exceptional vertex) and let $d(u_3)=k_3$. Again, from Stage I we have $\hc(\hD_1)=c(\hD_1)$. Then $\hD_1$ does not receive curvature across the $(u^*,u_3)$-edge (and one of the other edges of $\hD_1$ is a boundary edge) so the maximum amount of curvature to be distributed from $\hD_1$ is $c(\hD_1)+\pi/6 \leq c(k,k_3,4)+\pi/6 =(2\pi/k -\pi/6)+(2\pi/k_3-\pi/6)$. Distribute, when positive, $2\pi/k -\pi/6$, $2\pi/k_3-\pi/6$ from $\hD_1$ to $\D^*$ through $u^*,u_3$, respectively, as shown in both figures.
\end{itemize}

This completes Stage II and there is no further curvature distribution.

\begin{remark*}
It is clear from Figure {{K}} that although Property (P2) no longer holds in general, it still holds when $\hD$ receives curvature across its edge $e$ from a non-adjacent region; and it is also clear that properties (P1), (P3), (P4), (P5) and Claim \ref{claim:four} all still hold. (The remaining claims are discussed below.)
\end{remark*}

\subsection{Curvature after Stage II: concluding the proof}\label{sec:concludingtheproof}

\begin{notation*}   
From now on, for any given region $\hD$, we use the previous notation $c^*(\hD)$ to now denote $c(\hD)$ plus all the curvature $\hD$ receives minus all the curvature distributed from $\hD$ in Stage I or II.
\end{notation*}

If $\hD$ is interior having no exceptional vertices, then clearly Claim \ref{claim:one} still holds, that is $c^*(\hD)\leq 0$. Suppose otherwise and 
$\hD\neq \D^*$. If $\hD$ is not $\hD$ or $\hD_1$ of Figure {{K}}(ii)--(v) then, according to Stage II, $\hc(\hD)$ is distributed to $\D^*$ so $c^*(\hD)\leq 0$. If $\hD$ is $\hD$ of Figure {{K}}(ii)--(v) then $\hc(\hD)$ is distributed to $\hD_1$ and $\D^*$; or if $\hD$ is $\hD_1$ of Figure {{K}}(ii)--(v) then $\hc(\hD_1)=c(\hD_1)$ plus the $\pi/6$ or $2\pi/6$ that $\hD_1$ receives and this is distributed to $\D^*$. So in all cases $c^*(\hD)\leq 0$. This implies that $c^*(\D^*)\geq 4\pi$ and we show that this cannot happen. 

Observe that if curvature is distributed to $\D^*$ across an edge shared with a boundary region $\hD$ then it follows from Claim \ref{claim:two} and Figures {{K}}(ii),(iv),(vii) that the maximum amount distributed is $\pi/2$ except when $\hD$ has a boundary vertex $v$ of degree at least 8, and the maximum is then $7\pi/12$. Thus $\D^*$ may receive $\pi/12$ more than $\pi/2$ across each of the two edges of $\D^*$ that share $v$. In the case where $v$ is a boundary vertex of degree 8, at most $\pi/2+\pi/12=7\pi/12$ can be distributed to $\D^*$ across each of two edges, and so the most that can be distributed is $2\pi/8 + 7\pi/12+7\pi/12= 17\pi/12$. On the other hand, if $d(v)=4$ then the most that can be distributed is $2\pi/4+\pi/2+\pi/2>17\pi/12$. Thus, in order to maximise $c^*(\D^*)$, we may assume that $d(v)=4$ and that the two amounts distributed across the two edges are $\pi/2$.

\begin{notation*}
For an exceptional vertex $u^*$ of degree $k$, let $\tau(u^*)$ denote the sum $2\pi/d(u^*)=2\pi/k$ plus the total amount of curvature $\D^*$ receives through $u^*$. 
\end{notation*}

Let $u^*$ be an exceptional vertex of degree $k$.

\begin{itemize}
\item  Suppose $k\geq 4$ and let $\hD$ be an interior region, one of whose vertices is $u^*$. If $\hD$ has a second exceptional vertex then, as described in Stage II, $\D^*$ receives at most $2\pi/k-\pi/4\leq 0$ from $\hD$ through $u^*$, so suppose otherwise.

If $\hD$ is not given by $\hD$ or $\hD_1$ of Figure {{K}}(ii)--(v) then Claim \ref{claim:three} asserts that $\D^*$ receives at most $2\pi/k-\pi/6$ from $\hD$ through $u^*$ as in Figure {{K}}(viii). But checking the description of Stage II for Figures {{K}}(ii)--(v) shows that in each case $\D^*$ again receives at most $2\pi/k-\pi/6$ from $\hD$ through $u^*$.

Now let $\hD$ be a boundary region. If $\hD$ is given by Figure {{K}}(vii) with $u^*=u$ then $\D^*$ does not receive any curvature from $\hD$ through $u^*$. On the other hand, if $\hD$ is given by Figure {{K}}(vi) then, applying Claim \ref{claim:four}, the maximum amount of curvature that $\D^*$ can receive from $\hD$ through $u^*$ is $c(k,4,4)+\pi/3
=2\pi/k+\pi/3$.

Let $\tau(u^*)$ denote the sum $2\pi/d(u^*)=2\pi/k$ plus the total amount of curvature $\D^*$ receives through $u^*$. It follows from the above that (because there are $(k-3)$ interior regions and 2 boundary regions that can distribute curvature to $\D^*$ through $u^*$) we have
\begin{alignat}{1}
    \tau(u^*) &\leq (k-3)\left( 2\pi/k - \pi/6\right) +2(2\pi/k+\pi/3)+2\pi/d(u^*)\nonumber\\
    &= (19-k)\pi/6.\label{eq:tauu*}
\end{alignat} 
Therefore if $k\geq 5$ then $\tau(u^*)\leq 7\pi/3$ or if $k=4$ then $\tau(u^*)\leq 5\pi/2$. 

\item Suppose $k=3$, i.e. $u^*$ is an exceptional vertex of degree 3, as shown in Figure {{K}}(ix) in which it is assumed that neither $v_1$ nor $v_2$ is exceptional (since we assume that $\hD$ has exactly one exceptional vertex).

If $d(v_i)=4$ ($1\leq i\leq 3$) then neither $\hD_1$ nor $\hD_2$ receives any curvature from adjacent or non-adjacent regions and so $\tau(u^*)= 2\pi/d(u^*) + \hc(\hD_1)+ \hc(\hD_2)= 2\pi/d(u^*) + c(\hD_1)+ c(\hD_2) \leq   2\pi/3 +2c(3,4,4)=2\pi$. (This can be seen as follows. If $v_3$ is a boundary vertex, then, given that $v_1$ is a boundary vertex, as explained after the GR, the vertex $u^*$ has degree 2, a contradiction; therefore $v_3$ is interior. Consider the region $\hD_1$.  An inspection of Figures {{D}}(xiii) and {{H}}(xiii) shows that Property (P3) implies that $\hD_1$ does not receive any curvature from a non-adjacent region. Note also that regions adjacent to $\hD_1, \hD_2$ must be interior, since $d(v_1) = d(v_2) = 4$, and also that neither $v_1$ nor $v_2$ are exceptional, so only  Stage I(a) applies. The figures with a region in which curvature is distributed across an edge with one endpoint being boundary of degree 4 and the other endpoint being interior of degree 4 are Figures {{B}}(iv),(v),(vi),(vii). In Figures {{B}}(iv),(vi) the vertex corresponding to $v_2$ is exceptional or of degree at least 6, a contradiction; and in Figures {{B}}(v),(vii) the third vertex has degree greater than 3, a contradiction to $d(u^*)=3$; therefore $\hD_1$ does not receive curvature from an adjacent region. Similarly  $\hD_2$ does not receive curvature from an adjacent region.)

If $d(v_1)\geq 6$, $d(v_2)=4$, $d(v_3)=4$ or $d(v_1)=4$, $d(v_2)\geq 6$, $d(v_3)=4$ then $c(\hD_1)\leq c(3,4,6)$, $c(\hD_2)\leq c(3,4,4)$, so by GR the most $\hD_1$ can receive from its adjacent region is $c(4,4,6)=\pi/3$ and by Stage I, the most $\hD_2$ can receive from its adjacent region is $\pi/4$ as in Figures {{B}}(iv), (vi). Thus $\tau(u^*)=2\pi/d(u^*)+ \hc(\hD_1)+ \hc(\hD_2)\leq 2\pi/3+ (c(\hD_1) +\pi/3) + (c(\hD_2)+\pi/4) \leq 29\pi/12<5\pi/2$.
    
If $d(v_1)\geq 6$, $d(v_2)\geq 6$, $d(v_3)=4$ then $\tau(u^*) =2\pi/d(u^*)+\hc(\hD_1)+\hc(\hD_2)\leq 2\pi/3 + (c(\hD_1)+\pi/3) + (c(\hD_2)+\pi/3) = 2\pi/3+2c(3,4,6)+2\pi/3=7\pi/3$; or if $d(v_1)\geq 4$, $d(v_2)\geq 4$, $d(v_3)\geq 6$ then $\tau(u^*) =2\pi/d(u^*)+\hc(\hD_1)+\hc(\hD_2)\leq 2\pi/3 + (c(\hD_1)+\pi/3) + (c(\hD_2)+\pi/3) =2\pi/3 + 2c(3,4,6) +2\pi/3=7\pi /3$.

In conclusion $\tau(u^*)\leq 7\pi/3$.

    \item Suppose $k=2$, i.e. $u^*$ is an exceptional vertex of degree 2 as shown in Figure {{K}}(x) in which it is assumed that neither $v_1$ nor $v_2$ is exceptional (since we assume that $\hD$ has exactly one exceptional vertex). According to the curvature distribution scheme, the maximum amount that a region can receive from a non-adjacent region is $\pi/6$. Noting this, and Property (P5), the following hold:

    \begin{itemize}
        \item If $d(v_1)=d(v_2)=4$ then $\hD$ does not receive any curvature from non-adjacent regions and, as shown in Figure {{B}}(ii), $\hD$ can receive at most $c(4,4,4)=\pi/2$ from $\D$. Thus 
        $\tau(u^*)=2\pi/d(u^*)+\hc(\hD)\leq 2\pi/2 + (c(\hD)+\pi/2)\leq 2\pi/2+c(2,4,4)+\pi/2= 5\pi/2$.

        \item On the other hand, if at least one of $v_1$ or $v_2$ has degree at least $6$ then
        $\tau(u^*)=2\pi/d(u^*)+\hc(\hD)\leq 2\pi/2 + (c(\hD)+c(\D)+\pi/6) \leq 7\pi/3$.
    \end{itemize}
\end{itemize}

Let $d(\D^*)=k^*$. If $\D^*$ has exactly one exceptional vertex $u^*$ then, writing $d_1,\ldots , d_{k^*-1}$ for the degrees of the $k^*-1$ boundary vertices other than $u^*$, we have
\begin{alignat*}{1}
    c(\D^*) 
    &=(2-k^*)\pi +2\pi \sum_{i=1}^{k^{*}-1} \frac{1}{d_i} + \frac{2\pi}{d(u^*)}\\
    &\leq (2-k^*)\pi +2\pi \frac{(k^*-1)}{4} + \frac{2\pi}{d(u^*)}.
\end{alignat*}
Thus, since at most $\pi/2$ is distributed to $\D^*$ over $k^*-2$ boundary edges, zero curvature is distributed over the remaining two boundary edges by (P2), and $\tau(u^*)-2\pi/d(u^*)$ is distributed to $\D^*$ through $u^*$,  we have
\begin{alignat*}{1}
    c^*(\D^*) 
    &\leq c(\D^*) + (k^*-2)\frac{\pi}{2} + \left( \tau(u^*) -\frac{2\pi}{d(u^*)} \right)\\
    &\leq (2-k^*)\pi +2\pi \frac{(k^*-1)}{4} + \frac{2\pi}{d(u^*)} + (k^*-2)\frac{\pi}{2} + \left( \tau(u^*) -\frac{2\pi}{d(u^*)} \right)\\
    &= \pi/2 +\tau(u^*)\\
    &\leq \pi/2+ 5\pi/2=3\pi.
\end{alignat*}
So to attain $c^*(\D^*) \geq 4\pi$ we require $\hD$ to have two exceptional vertices. Before proceeding with this case, note that having obtained a contradiction in the case where there is at most one exceptional vertex, we have dealt with the cases in \eqref{eq:starandstarstar} of Section \ref{sec:generalsetup} where either $\beta=0$ or $\alpha=0$ and this shows that the elements $A$ and $B$ have infinite order in $E_{2m}$. Now assume that there are exactly two exceptional vertices $u^*$ and $v^*$ of degrees $k_1$ and $k_2$ respectively.
First assume that $u^*$ and $v^*$ are not adjacent on the boundary of $\D^*$. Then there are $(k^*-2)$ vertices that are not $u^*,v^*$ of degree at least 4 and there are $k^*-4$ boundary edges that are not incident either to $u^*$ or $v^*$, and so at most $\pi/2$ can be distributed over those edges and $0$ over the remaining 4 edges, and hence
\begin{alignat*}{1}
c^* (\D^*)&\leq (2-k^*)\pi + (k^*-2)\cdot 2\pi/4 + (k^*-4)\pi/2 + \tau(u^*)+\tau(v^*)\\
&=-\pi + \tau(u^*)+\tau(v^*)
\end{alignat*}
so if $c(\D^*)\geq 4\pi$ we must have $\tau(u^*)+\tau(v^*)\geq 5\pi$. Thus, by \eqref{eq:tauu*}, $d(u^*), d(v^*)\leq 4$. But, as shown above, if $d(u^*)=3$ or $d(v^*)=3$ then $\tau(u^*)\leq 7\pi/3$, $\tau(v^*)\leq 7\pi/3$, respectively. Thus $\{d(u^*), d(v^*)\}\in \{2,4\}$. But at most one of $u^*,v^*$ is an exceptional vertex $u_*$ created in the detaching process, so by Lemma \ref{lem:stargraphlemma}(iii), without any loss, we may assume $l(u^*)=fa^{-1}$.
    
For the remainder of the proof, we use the fact that the maximum amount of curvature a region can receive from a non-adjacent region is $\pi/6$. Suppose that $u^*$ is given by Figure {{K}}(x).
    
If $d(v_1)\geq 6$, $d(v_2)\geq 6$, $d(v_3)\geq 6$ then $\tau(u^*)= 2\pi/d(u^*)+ c(\hD)+\hc(\D)\leq 2\pi/2+ c(2,6,6) +(c(6,6,6)+2\pi/6)= 2\pi$; if $d(v_1)\geq 6$, $d(v_2)\geq 6$, $d(v_3)=4$ then $\tau(u^*)= 2\pi/d(u^*)+ c(\hD)+\hc(\D)\leq 2\pi/2+ c(2,6,6) +(c(4,6,6)+2\pi/6)= 13\pi/6$; if ($d(v_1)\geq 6$, $d(v_2)=4$, $d(v_3)\geq 6$) or ($d(v_1)=4$, $d(v_2)\geq 6$, $d(v_3)\geq 6$) then $\tau(u^*)= 2\pi/d(u^*)+ c(\hD)+\hc(\D)\leq 2\pi/2+ c(2,4,6) +(c(4,6,6)+2\pi/6)= 7\pi/3$; if ($d(v_1)=4$, $d(v_2)\geq 6$, $d(v_3)=4$) or ($d(v_1)\geq 6$, $d(v_2)=4$, $d(v_3)=4$) then $\tau(u^*)= 2\pi/d(u^*)+ c(\hD)+\hc(\D)\leq 2\pi/2+ (c(2,4,6)+0) +(c(4,6,6)+\pi/6)=7\pi/3$ (where, for the first of the two cases for example, we used the fact that at most $\pi/6$ can be distributed from the adjacent region to $\D$ across the $(v_2,v_3)$-edge, and 0 can be distributed from the adjacent region to $\D$ across the $(v_1,v_3)$-edge, as both vertices have degree 4). This leaves $d(v_1)=d(v_2)=4$, as shown in Figure {{K}}(xi). But $d(v_3)=4$ in Figure {{K}}(xi) forces $l(v_3)=ca^{-1}b\mu$ or $ea^{-1}b\mu$ and this would force two further exceptional vertices. (To see this note that with this labelling, Figure {{K}}(xi) then gives a diagram with 4 regions plus the exterior region, with two boundary vertices of degree 2, with label $\hat{b}w$. These are exceptional vertices because only exceptional vertices have degree 2 by Lemma \ref{lem:stargraphlemma}.) Therefore $d(v_3)\geq 6$ and $\tau(u^*)\leq 2\pi/2 + c(2,4,4)+c(4,4,6)=7\pi/3$ (where we use the fact that no curvature is distributed to $\hD$ from a non-adjacent region to $\hD$, since it would involve distributing across an edge with two degree 4 vertices, namely $v_1,v_2$). In conclusion, since $\tau(u^*)\leq 7\pi/3$ and $\tau(v^*)\leq 5\pi/2$, regardless of its degree, we have $\tau(u^*)+\tau(v^*)\leq 7\pi/3+5\pi/2<5\pi$, a contradiction.
    
Now assume that $u^*$ and $v^*$ are adjacent on the boundary of $\D^*$. Then since there are $(k-3)$ boundary edges not incident to either $u^*$ or $v^*$, and at most $\pi/2$ is distributed across each of these edges, and zero over the remaining 3 edges, and the degrees of the boundary vertices other than $u^*,v^*$ are at least 4, we have
\begin{alignat*}{1}
c^*(\D^*) &\leq (2-k^*)\pi + (k^*-2)\cdot 2\pi/4 +(k^*-3)\pi/2 + \tau(u^*)+\tau(v^*)\\
&= -\pi/2 +\tau (u^*)+\tau (v^*)
\end{alignat*}
so if $c^*(\D^*)\geq 4\pi$ we must have $\tau (u^*)+\tau (v^*)\geq 9\pi/2$. To complete the proof we show that \linebreak $\max\{\tau (u^*), \tau (v^*)\}\leq 2\pi$, a contradiction. Consider Figure {{K}}(xii).

Let $k_1=d(u^*)\geq 4$. Given $d(u^*)\geq 4$, suppose that $k_2=d(v^*)\geq 3$. Then $\hD$ of Figure {{K}}(xii) cannot coincide with $\hD_1$ of Figure {{K}}(iii),(v). (To see this, observe that if $\D$ of {{K}}(xii) coincides with $\hD_1$ of {{K}}(ii) then vertex $u_3$ must coincide with $v^*$, because $(u^*,u_1)$ is an interior edge. But then $\hD_1$ of {{K}}(ii) has two boundary edges $(u^*,u_3)$ and $(u_3,u_1)$ and hence $d(v^*)=d(u_3)=2$, a contradiction.) Moreover, Figure {{K}}(xii) cannot coincide with $\hD_1$ of Figure {{K}}(ii),(iv), since in those cases $u_3$ is not exceptional.  Therefore, using (P2), and checking the remaining figures in Figure {{K}} we see that $\hD$ does not receive curvature across its edges. Observe then that $\hc(\hD)=c(\hD)\leq c(k_1,k_2,4)$ and so, as in Figure {{K}}(viii), $2\pi/k_1-\pi/4$, $2\pi/k_2-\pi/4$ is distributed from $\hD$ to $\D^*$ through $u^*,v^*$, respectively. There are at most $k_1-3$ interior regions incident to $u^*$, each of which distribute at most $2\pi/k_1-\pi/6$ to $\D^*$ through $u^*$ (by Claim \ref{claim:three}); the boundary region incident to $u^*$ has curvature at most $c(k_1,4,4)$, at most $(2\pi/k_1-\pi/4)$ is distributed to $\D^*$ through $u^*$, and at most $\pi/3$ is distributed from the adjacent region to $\hD$ that is not incident to $u^*$ (by Claim \ref{claim:four}). It follows that 
\begin{alignat*}{1}
\tau(u^*)&\leq (k_1-3)\left(\frac{2\pi}{k_1}-\frac{\pi}{6}\right) + c(k_1,4,4)+\pi/3 +\left( \frac{2\pi}{k_1}-\frac{\pi}{4} \right) + \frac{2\pi}{k_1}\\
&= \frac{15-k_1}{6}\pi\leq \frac{15-4}{6}\pi=11\pi/6.
\end{alignat*}
Given $d(u^*)\geq 4$, suppose now that $d(v^*)=2$. Then $\hD$ can coincide with $\hD_1$ of Figure {{K}}(iii) or (v) and so $\hD$ has at most $c(k_1,k_2,4)+\pi/6=c(k_1,2,4)+\pi/6=(2\pi/k_1-\pi/6)+(2\pi/2-\pi/6)$ to distribute to $\hD^*$ through $u^*$ and $v^*$. Therefore
\begin{alignat*}{1}
\tau(u^*)&\leq (k_1-3)\left(\frac{2\pi}{k_1}-\frac{\pi}{6}\right) + c(k_1,4,4)+\pi/3 +\left( \frac{2\pi}{k_1}-\frac{\pi}{6} \right) + \frac{2\pi}{k_1}\\
&= \frac{16-k_1}{6}\pi\leq \frac{16-4}{6}\pi=2\pi.
\end{alignat*}
Thus if $d(u^*)\geq 4$ then $\tau(u^*)\leq 2\pi$. Similarly, if $d(v^*)\geq 4$ then $\tau(v^*)\leq 2\pi$.
We have also shown that if $d(v^*)=2$ and $d(u^*)\geq 4$ then $\tau(v^*)\leq 2\pi$ and so again, by symmetry, if $d(u^*)=2$ and $d(v^*)\geq 4$ then $\tau(u^*)\leq 2\pi$.

Now let $d(u^*)=3$. Given this, suppose that $d(v^*)\geq 3$. We may take $u^*$ to be the vertex in Figure {{K}}(ix) and $v^*=v_2$. Thus (since $v_1,v_3$ are not exceptional) in Figure {{K}}(ix) we have $d(v_1)\geq 4$ and $d(v_3)\geq 4$. The region $\hD_2$ of Figure {{K}}(ix) cannot coincide with $\hD_1$ of Figure {{K}}(iii) or (v) and so there is $\hc(\hD_2)\leq c(3,4,k_2)=(2\pi/3-\pi/4)+(2\pi/k_2-\pi/4)$ to distribute from $\hD_2$ to $\D^*$ through $u^*$ and $v^*$. We have $d(v_1)=d(v_3)=4$ or $d(v_1)\geq 4, d(v_3)\geq 6$ or $d(v_1)\geq 6, d(v_3)\geq 4$. In the first of these cases, since no curvature can pass over the edge with two endpoints of degree 4, we have
\[\tau(u^*)\leq c(3,4,4)+(2\pi/3-\pi/4)+2\pi/3= 7\pi/4.\]
In each of the latter two cases the region $\hD_1$ of Figure {{K}}(ix) can receive at most $c(4,4,6)=\pi/3$ from its adjacent region, so
\[\tau(u^*)\leq (c(3,4,6)+\pi/3)+(2\pi/3-\pi/4)+2\pi/3= 23\pi/12.\]
Therefore $\tau(u^*)<2\pi$ and the same holds for $v^*$ if $d(v^*)=3$ and $d(u^*)\geq 3$.

Finally, this leaves ($d(u^*)=2$ and $d(v^*)=3$) or ($d(u^*)=3$ and $d(v^*)=2$) and, without loss of generality, we may assume the former. It can be assumed without loss that $u^*=u^*$ of Figure {{K}}(x) with $v^*=v_2$. Thus the $(v_2,v_3)$-edge, i.e.\,the $(v^*,v_3)$-edge, is an edge of $\D^*$. Then $\hD$ of Figure {{K}}(x) can coincide with $\hD_1$ of Figure {{K}}(iii) or (vi) and so $\hD$ has at most $c(2,3,4)+\pi/6=(2\pi/2-\pi/6)+(2\pi/3-\pi/6)$ to distribute to $\D^*$ through $u^*$ and $v^*$. It follows that
\[ \tau(u^*)\leq (2\pi/2-\pi/6)+2\pi/2=11\pi/6\]
and (using Claim \ref{claim:four} and the fact that non-exceptional boundary vertices have degree at least 6) that
\[ \tau(v^*)\leq (c(3,4,6)+\pi/3) + (2\pi/3-\pi/6)+2\pi/3=2\pi.\]

Thus $\max\{\tau (u^*), \tau (v^*)\}\leq 2\pi$,  and the proof of Main Lemma is complete.

\section*{Funding}

This work was partially supported by the Leverhulme Trust Research Project Grant RPG-2017-334.

\newpage
\subsection*{Figure {{A}}}
\psset{arrowscale=0.7,linewidth=1.0pt,dotsize=3pt}
\psset{xunit=1.4cm,yunit=1.4cm,}
\begin{center}

\begin{tabular}{@{}c@{}c}

\begin{pspicture*}(2.,-4.)(6.,0.5)
\psset{dotscale=2}

\psdots[dotstyle=*](4.,0.)
\psdots[dotstyle=*](3.,-2.)
\psdots[dotstyle=*](5.,-2.)

\psline[ArrowInside=->,arrowsize=0.3](5.,-2.)(3.,-2.)
\psline[ArrowInside=->,arrowsize=0.3](3.,-2.)(4.,0.)
\psline[ArrowInside=->,arrowsize=0.3](4.,0.)(5.,-2.)

\psarc[arrowsize=0.3]{->}(3,-2){0.8}{0}{42}
\psarc(3,-2){0.8}{0}{63}
\psarc[arrowsize=0.3]{->}(5,-2){0.8}{117}{159}
\psarc(5,-2){0.8}{117}{180}
\psarc[arrowsize=0.3]{->}(4,0){0.8}{243}{285}
\psarc(4,0){0.8}{243}{298}

\rput(4.0,-1.2){$\Delta$}
\rput(4.9,-1){$x$}
\rput(3.1,-1){$x$}
\rput(4,-2.3){$x$}

\rput(3.25,-1.8){$\lambda$}
\rput(4.7,-1.8){$b$}
\rput(4,-0.35){$a$}

\rput(4.0,-2.8){A (i)}
\end{pspicture*}

&
\begin{pspicture*}(2.,-4.)(6.,0.5)
\psset{dotscale=2}

\psdots[dotstyle=*](4.,0.)
\psdots[dotstyle=*](3.,-2.)
\psdots[dotstyle=*](5.,-2.)

\psline[ArrowInside=->,arrowsize=0.3](3.,-2.)(5.,-2.)
\psline[ArrowInside=->,arrowsize=0.3](4.,0.)(3.,-2.)
\psline[ArrowInside=->,arrowsize=0.3](5.,-2.)(4.,0.)

\psarc[arrowsize=0.3]{<-}(3,-2){0.8}{21}{63}
\psarc(3,-2){0.8}{0}{63}
\psarc[arrowsize=0.3]{<-}(5,-2){0.8}{138}{180}
\psarc(5,-2){0.8}{117}{180}
\psarc[arrowsize=0.3]{<-}(4,0){0.8}{256}{298}
\psarc(4,0){0.8}{243}{298}

\rput(4.0,-1.2){$\Delta^{-1}$}
\rput(4.9,-1){$x$}
\rput(3.1,-1){$x$}
\rput(4,-2.3){$x$}

\rput(3.25,-1.8){$b$}
\rput(4.7,-1.8){$\lambda$}
\rput(4,-0.35){$a$}
\rput(4.0,-2.8){A (ii)}
\end{pspicture*}
\\
\end{tabular}%
\end{center}

\newpage
\subsection*{Figure {{B}}}

\psset{arrowscale=0.7,linewidth=1.0pt,dotsize=3pt}
\psset{xunit=.8cm,yunit=1.0cm,}
\begin{center}

\begin{tabular}{@{}c@{}c@{}c}
\begin{pspicture*}(0.,-4.)(6.,0.5)
\psdots[dotstyle=*](4.,0.)
\psdots[dotstyle=*](3.,-2.)
\psdots[dotstyle=*](5.,-2.)
\psline(3.,-2.)(5.,-2.)

\psline(3.,-2.)(4.,0.)
\psline(4.,0.)(5.,-2.)
\rput(4.0,-1.0){$\Delta$}
\rput(5.0,-0.5){$\hat{\Delta}$}
%
\psline{->}(4.3,-1.1)(4.8,-.9)
\rput(5.0,-1.0){$\frac{\pi}{12}$}

\rput(4.0,-3.5){B (i)}
\end{pspicture*}
&
\begin{pspicture*}(0.,-4.)(7.5,0.5)
\psdots[dotstyle=*](4.,0.)
\psdots[dotstyle=*](3.,-2.)
\psdots[dotstyle=*](5.,-2.)

\psdots[dotstyle=*](6.1,-0.5)
\rput(6.5,-0.4){$u^*$}
\psbezier
(4,0)(6.5,0.5)
(7,-1.5)(5,-2)
\psline(3.,-2.)(5.,-2.)

\psline(3.,-2.)(4.,0.)
\psline(4.,0.)(5.,-2.)
\rput(4.0,-1.0){$\Delta$}
\rput(5.0,-0.5){$\hat{\Delta}$}
%
\rput(4.0,0.3){$v_1$}
\rput(5.0,-2.3){$v_2$}
\psline{->}(4.3,-1.1)(4.8,-.9)
\rput(5.5,-1.0){$c(\Delta)$}

\rput(4.0,-3.5){B (ii)}
\end{pspicture*}
 &

 \begin{pspicture*}(1.,-4.)(6.,0.5)
\psdots[dotstyle=*](4.,0.)
\psdots[dotstyle=*](3.,-2.)
\psdots[dotstyle=*](5.,-2.)
\psline(2.,0.)(6.,0.) \psline(3.,-2.)(5.,-2.)

\psline(3.,-2.)(4.,0.)
\psline(4.,0.)(5.,-2.)
\rput(4.0,0.3){$\Delta^*$}



\rput(2.5,-2.3){$\geq 6$}
\rput(3.8,-2.5){$\hat{\Delta}_3$}
\rput(5.5,-2.3){$\geq 6$}
%
\rput(2.5,-.5){$\hat{\Delta}_2$}
\rput(4.0,-1.0){$\Delta$}
\rput(5.5,-.5){$\hat{\Delta}_1$}
%

\psline{->}(3.7,-1.1)(3.2,-0.9)
\psline{->}(4.3,-1.1)(4.8,-.9)
\psline[linestyle=dotted]{->}(4.0,-1.7)(4.0,-2.3)
\rput(3.0,-1.0){$\frac{\pi}{12}$}
\rput(5.0,-1.0){$\frac{\pi}{12}$}
%
\rput(4.0,-3.5){B (iii)}
\end{pspicture*}

\end{tabular}

\begin{tabular}{@{}c@{}c}

\begin{pspicture*}(-2.,-4.)(8.,0.5)
\psdots[dotstyle=*](2.,0.)
\psdots[dotstyle=*](4.,0.)
\psdots[dotstyle=*](6.,0.)
\psdots[dotstyle=*](1.,-2.)
\psdots[dotstyle=*](3.,-2.)
\psdots[dotstyle=*](5.,-2.)
\psline(2.,0.)(6.,0.) \psline(1.,-2.)(5.,-2.)

\psline(1.,-2.)(2.,0.)
\psline(2.,0.)(3.,-2.)
\psline(3.,-2.)(4.,0.)
\psline(4.,0.)(5.,-2.)
\psline(5.,-2.)(6.,0.)

\rput(2.0,0.3){}
\rput(4.0,0.3){$c$}
\rput(6.0,0.3){}

\rput(1.5,-0.3){}
\rput(2.0,-0.5){$a$}
\rput(2.5,-0.3){$\hat{b}$}
\rput(3.5,-0.3){$\hat{a}$}
\rput(4.0,-0.5){$\lambda$}
\rput(4.5,-0.3){$\hat{b}$}
\rput(5.5,-0.3){$\hat{a}$}
\rput(6.0,-0.5){}

\rput(0.5,-1.7){}
\rput(1.0,-1.5){}
\rput(1.5,-1.7){$\lambda$}
\rput(2.5,-1.7){$b$}
\rput(3.0,-1.5){$\mu$}
\rput(3.5,-1.7){$b$}
\rput(4.5,-1.7){$a$}
\rput(5.0,-1.5){$\mu$}
\rput(5.5,-1.7){}
%
\rput(0.5,-2.3){$\encircle{\geq 6}$}
\rput(1.5,-2.3){$\hat{a}$}
\rput(3.0,-2.3){$\mu$}
\rput(4.5,-2.3){$\hat{b}$}
\rput(5.5,-2.3){$\geq 6$}
%
\rput(4.0,-1.0){$\Delta$}
%
\psbezier
(1.,-2.)(1.3,-3.3)
(4.7,-3.3)(5,-2.)
%
\psline{->}(3.7,-1.1)(3.2,-0.9)
\psline{->}(4.3,-1.1)(4.8,-.9)
\rput(3.0,-1.0){$\frac{\pi}{4}$}
\rput(5.0,-1.0){$\frac{\pi}{12}$}
%
\rput(4.0,-3.5){B (iv)}
\end{pspicture*}


&

\begin{pspicture*}(-2.,-4.)(8.,0.5)
\psdots[dotstyle=*](1.,0.)
\psdots[dotstyle=*](4.,0.)
\psdots[dotstyle=*](6.,0.)
\psdots[dotstyle=*](-1.,-2.)
\psdots[dotstyle=*](1.,-2.)
\psdots[dotstyle=*](3.,-2.)
\psdots[dotstyle=*](5.,-2.)
\psline(0.,0.)(7.,0.) 
\psline(-1.,-2.)(5.,-2.)

\psline(1.,0.)(-1.,-2.)
\psline(1.,-2.)(1.,0.)
\psline(1.,0.)(3.,-2.)
\psline(3.,-2.)(4.,0.)
\psline(4.,0.)(5.,-2.)
\psline(5.,-2.)(6.,0.)
\rput(2.0,0.3){}
\rput(4.0,0.3){$c$}
\rput(6.0,0.3){}

\rput(0.7,-0.6){$\mu$}
\rput(1.2,-0.5){$a$}
\rput(1.7,-0.3){$\hat{b}$}
\rput(3.5,-0.3){$\hat{a}$}
\rput(4.0,-0.5){$\lambda$}
\rput(4.5,-0.3){$\hat{b}$}
\rput(5.5,-0.3){$\hat{a}$}
\rput(6.0,-0.5){}

\rput(-.3,-1.7){$\hat{a}$}
\rput(0.7,-1.6){$\hat{b}$}
\rput(1.5,-1.7){$\lambda$}
\rput(2.4,-1.7){$b$}
\rput(3.0,-1.5){$\mu$}
\rput(3.5,-1.7){$b$}
\rput(4.5,-1.7){$a$}
\rput(5.0,-1.5){$\mu$}
\rput(5.5,-1.7){}
%
\rput(0.5,-2.3){$c/e$}
\rput(1.5,-2.3){$\hat{a}$}
\rput(3.0,-2.3){$\mu$}
\rput(4.5,-2.3){$\hat{b}$}
\rput(5.5,-2.3){$\geq 6$}
%
\rput(4.0,-1.0){$\Delta$}
%
\psbezier
(1.,-2.)(1.3,-3.3)
(4.7,-3.3)(5,-2.)
%
\psline{->}(3.7,-1.1)(3.2,-0.9)
\psline{->}(4.3,-1.1)(4.8,-.9)
\rput(3.0,-1.0){$\frac{\pi}{6}$}
\rput(5.0,-1.0){$\frac{\pi}{12}$}
\psline{->}(4.0,-1.8)(4.0,-2.2)
\rput(4.0,-2.5){$\frac{\pi}{12}$}
%
\rput(4.0,-3.5){B (v)}

\end{pspicture*}

\\
\end{tabular}
\end{center}

\newpage

\begin{center}

\begin{tabular}{@{}c@{}c}

%
\begin{pspicture*}(-2.,-4.)(8.,0.5)
\psdots[dotstyle=*](0.,0.)
\psdots[dotstyle=*](2.,0.)
\psdots[dotstyle=*](4.,0.)
\psdots[dotstyle=*](6.,0.)
\psdots[dotstyle=*](-1.,-2.)
\psdots[dotstyle=*](1.,-2.)
\psdots[dotstyle=*](3.,-2.)
\psdots[dotstyle=*](5.,-2.)
\psline(-1.,0.)(7.,0.) \psline(-1.,-2.)(5.,-2.)

\psline(-1.,-2.)(0.,0.)
\psline(0.,0.)(1.,-2.)
\psline(1.,-2.)(2.,0.)
\psline(2.,0.)(3.,-2.)
\psline(3.,-2.)(4.,0.)
\psline(4.,0.)(5.,-2.)
\psline(5.,-2.)(6.,0.)
\rput(2.0,0.3){$d$}
\rput(4.0,0.3){$e$}
\rput(6.0,0.3){}

\rput(0.,-0.5){$b$}
\rput(0.5,-0.3){$\hat{a}$}
\rput(1.5,-0.3){$\mu$}
\rput(2.0,-0.5){$a$}
\rput(2.5,-0.3){$\mu$}
\rput(3.5,-0.3){$\hat{b}$}
\rput(4.0,-0.5){$\lambda$}
\rput(4.5,-0.3){$\hat{a}$}
\rput(5.5,-0.3){$\mu$}
\rput(6.0,-0.5){}

\rput(-.5,-1.7){$a$}
\rput(0.5,-1.7){$\lambda$}
\rput(1.0,-1.5){$\hat{b}$}
\rput(1.5,-1.7){$\lambda$}
\rput(2.5,-1.7){$b$}
\rput(3.0,-1.5){$\hat{a}$}
\rput(3.5,-1.7){$b$}
\rput(4.5,-1.7){$a$}
\rput(5.0,-1.5){$\hat{b}$}
\rput(5.5,-1.7){}
%
\rput(-.5,-2.3){$\mu$}
\rput(1.,-2.3){$\hat{b}$}
\rput(2.3,-2.3){$\hat{a}$}
\rput(3.5,-2.3){$\geq 8$}
\rput(4.5,-2.3){}
\rput(5.,-2.3){$\geq 6$}
%
\rput(2.0,-1.0){$\Delta$}
\rput(4.0,-1.0){$\Delta_1$}
\rput(5.0,-.5){$\hat{\Delta}$}
%
\psbezier
(-1.,-2.)(-0.7,-3.3)
(2.7,-3.3)(3,-2.)
%
\psline{->}(2.3,-1.1)(2.8,-.9)
\rput(3.0,-1.0){$\frac{\pi}{4}$}
\psline{->}(4.3,-1.1)(4.8,-.9)
\rput(5.0,-1.0){$\frac{\pi}{12}$}
%
\rput(4.0,-3.5){B (xii)}
\end{pspicture*}

&
%
%
\begin{pspicture*}(-2.,-4.)(8.,0.5)
\psdots[dotstyle=*](2.,0.)
\psdots[dotstyle=*](4.,0.)
\psdots[dotstyle=*](6.,0.)
\psdots[dotstyle=*](1.,-2.)
\psdots[dotstyle=*](3.,-2.)
\psdots[dotstyle=*](5.,-2.)
\psline(-.5,0.)(6.0,0.) \psline(0.,-2.)(5.,-2.)

\psline(0.,0.)(1.,-2.)
\psline(1.,-2.)(2.,0.)
\psline(2.,0.)(3.,-2.)
\psline(3.,-2.)(4.,0.)
\psline(4.,0.)(5.,-2.)
\psline(5.,-2.)(6.,0.)
\psline(5,-2.0)(5.,-2.5)
\rput(2.0,0.3){$d$}
\rput(4.0,0.3){$e$}
\rput(6.0,0.3){$u^*$}

\rput(0.,-0.5){}
\rput(0.5,-0.3){$\hat{a}$}
\rput(1.5,-0.3){$\mu$}
\rput(2.0,-0.5){$a$}
\rput(2.5,-0.3){$\mu$}
\rput(3.5,-0.3){$\hat{b}$}
\rput(4.0,-0.5){$\lambda$}
\rput(4.5,-0.3){$\hat{a}$}
\rput(5.5,-0.3){$\mu$}
\rput(6.0,-0.5){}

%
\rput(0.5,-1.7){$\lambda$}
\rput(1.0,-1.5){$\hat{b}$}
\rput(1.5,-1.7){$\lambda$}
\rput(2.5,-1.7){$b$}
\rput(3.0,-1.5){$\hat{a}$}
\rput(3.5,-1.7){$b$}
\rput(4.5,-1.7){$a$}
\rput(5.0,-1.5){$\hat{b}$}
\rput(5.5,-1.7){}
%
\rput(1.,-2.3){$\hat{b}$}
\rput(2.3,-2.3){$\hat{a}$}
\rput(3.5,-2.3){$\geq 8$}
\rput(4.5,-2.3){}
\rput(4.7,-2.3){$\hat{b}$}
\rput(5.3,-2.1){$d$}
%
\rput(2.0,-1.0){$\Delta$}
\rput(4.0,-1.0){$\Delta_1$}
\rput(5.0,-.5){$\hat{\Delta}$}
%
\psbezier
(2.0,-2.7)(2.5,-2.7)
(2.7,-2.5)(3,-2.)
%
\psline{->}(4.3,-1.1)(4.8,-.9)
\rput(5.0,-1.0){$\frac{\pi}{4}$}
%
\rput(4.0,-3.5){B (xiii)}
\end{pspicture*}
\\
\end{tabular}
\end{center}

\newpage
\subsection*{Figure {{C}}}
\psset{arrowscale=0.7,linewidth=1.0pt,dotsize=3pt}
\psset{xunit=.8cm,yunit=1.0cm,}
\begin{center}

\begin{tabular}{@{}c@{}c@{}c}

\begin{pspicture*}(2.,-4.)(7.5,1.5)
\psdots[dotstyle=*](4.,0.)
\psdots[dotstyle=*](6.,0.)
\psdots[dotstyle=*](4.,-2.)
\psdots[dotstyle=*](6.,-2.)
\psline(4.,0)(6.,0.)
\psline(4.,-2.)(6.,-2.)
\psline(4.,0)(4.,-2.)
\psline(6.,0)(6.,-2.)
\psline(4.,0)(6.,-2.)
\psline(3.3,0.)(4.,0.)
\psline(3.3,-2.)(4.,-2.)
\psline(4.,-2.)(4.,-2.7)
\psline(6.,-2.)(6.,-2.7)
%
%
\rput(4.0,0.2){$u_2$}
\rput(6.0,0.2){$u_3$}
\rput(3.5,-0.2){$\mu$}
\rput(4.2,-0.5){$a$}
\rput(3.5,-1.7){$\hat{b}$}
\rput(4.3,-1.7){$\lambda$}
\rput(5.5,-1.8){$b$}
\rput(6.5,-2.){$u_1$}
\rput(3.5,-2.3){$\lambda$}
\rput(4.3,-2.3){$\hat{b}$}
\rput(5.5,-2.3){$\hat{a}$}
%
\rput(4.8,-1.4){$\Delta$}
\rput(5.2,-0.5){$\hat{\Delta}$}
%
\psline[linestyle=dotted]{->}(4.5,-1.3)(5.5,-.7)
\rput(5.5,-1.0){$\frac{\pi}{6}$}

\rput(5.0,-3.3){C (i)}
\end{pspicture*}

&

\begin{pspicture*}(2.,-4.)(7.5,1.5)
\psdots[dotstyle=*](4.,0.)
\psdots[dotstyle=*](6.,0.)
\psdots[dotstyle=*](4.,-2.)
\psdots[dotstyle=*](6.,-2.)
\psline(4.,0)(6.,0.)
\psline(4.,-2.)(6.,-2.)
\psline(4.,0)(4.,-2.)
\psline(6.,0)(6.,-2.)
\psline(4.,0)(6.,-2.)
\psline(3.3,0.)(4.,0.)
\psline(3.3,-2.)(4.,-2.)
\psline(4.,0.)(4.,0.7)
\psline(6.,0.)(6.,0.7)
\psline(4.,-2.)(4.,-2.7)
\psline(6.,-2.)(6.,-2.7)
\psline(4.,0.)(3.3,0.7)
%
\rput(3.5,0.2){$b$}
\rput(3.8,0.5){$\hat{a}$}
\rput(4.3,0.2){$\lambda$}
\rput(5.6,0.2){$b$}
\rput(6.6,0.0){$\geq 6$}
\rput(3.5,-0.2){$\mu$}
\rput(4.2,-0.5){$a$}
\rput(4.6,-0.2){$\hat{b}$}
\rput(5.6,-0.2){$\hat{a}$}
\rput(3.5,-1.7){$\hat{b}$}
\rput(4.3,-1.7){$\lambda$}
\rput(5.5,-1.8){$b$}
\rput(5.8,-1.5){$\mu$}
\rput(6.6,-2.){$\geq 8$}
\rput(3.5,-2.3){$\lambda$}
\rput(4.3,-2.3){$\hat{b}$}
\rput(5.5,-2.3){$\hat{a}$}
%
\rput(4.8,-1.4){$\Delta$}
\rput(5.2,-0.5){$\hat{\Delta}$}
%
\psline{->}(4.5,-1.3)(5.5,-.7)
\rput(5.5,-1.0){$\frac{\pi}{12}$}

\rput(5.0,-3.3){C (ii)}
\end{pspicture*}

&

\begin{pspicture*}(2.,-4.)(7.5,1.5)
\psdots[dotstyle=*](4.,0.)
\psdots[dotstyle=*](6.,0.)
\psdots[dotstyle=*](4.,-2.)
\psdots[dotstyle=*](6.,-2.)
\psline(4.,0)(6.,0.)
\psline(4.,-2.)(6.,-2.)
\psline(4.,0)(4.,-2.)
\psline(6.,0)(6.,-2.)
\psline(4.,0)(6.,-2.)
\psline(3.3,0.)(4.,0.)
\psline(3.3,-2.)(4.,-2.)
\psline(4.,0.)(4.,0.7)
\psline(6.,0.)(6.,0.7)
\psline(4.,-2.)(4.,-2.7)
\psline(6.,-2.)(6.,-2.7)
\psline(4.,0.)(3.3,0.7)
%
\rput(3.5,0.2){$b$}
\rput(3.8,0.5){$\hat{a}$}
\rput(4.3,0.2){$b$}
\rput(5.6,0.2){$a$}
\rput(6.6,0.0){$\geq 6$}
\rput(3.5,-0.2){$\mu$}
\rput(4.2,-0.5){$a$}
\rput(4.6,-0.2){$\mu$}
\rput(5.6,-0.2){$\hat{b}$}
\rput(3.5,-1.7){$\hat{b}$}
\rput(4.3,-1.7){$\lambda$}
\rput(5.5,-1.8){$b$}
\rput(5.8,-1.5){$\hat{a}$}
\rput(6.6,-2.){$\geq 8$}
\rput(3.5,-2.3){$\lambda$}
\rput(4.3,-2.3){$\hat{b}$}
\rput(5.5,-2.3){$\hat{a}$}
%
\rput(4.8,-1.4){$\Delta$}
\rput(5.2,-0.5){$\hat{\Delta}$}
%
\psline{->}(4.5,-1.3)(5.5,-.7)
\rput(5.5,-1.0){$\frac{\pi}{12}$}

\rput(5.0,-3.3){C (iii)}
\end{pspicture*}
\\
\end{tabular}

\end{center}

\newpage
\subsection*{Figure {{D}}}
\psset{arrowscale=0.7,linewidth=1.0pt,dotsize=3pt}
\psset{xunit=1.0cm,yunit=1.0cm,}

\begin{center}
\begin{tabular}{c@{}c}

\begin{pspicture*}(1.,2.5)(9.,7.)
\psdots[dotstyle=*](4.,4.)
\rput(3.7,4.3){$\hat{b}$}
\rput(4.2,4.2){$\lambda$}
\rput(4.2,3.8){$\hat{b}$}
\rput(3.7,3.8){$\lambda$}

\psdots[dotstyle=*](4.,6.)
\rput(3.5,6.2){$b$}
\rput(3.8,6.5){$\hat{a}$}
\rput(4.2,6.2){$\lambda$}
\rput(4.5,5.8){$\hat{b}$}
\rput(4.2,5.6){$a$}
\rput(3.7,5.8){$\mu$}

\psdots[dotstyle=*](6.,4.)
\rput(5.6,4.2){$b$}
\rput(5.8,4.5){$\mu$}
\rput(6.2,4.2){$b$}
\rput(6.5,3.8){$\hat{d}$}
\rput(6.2,3.6){$\lambda$}
\rput(5.7,3.8){$\hat{a}$}

\psdots[dotstyle=*](6.,6.)
\rput(5.7,6.2){$b$}
\rput(6.2,5.6){$\lambda$}
\rput(5.7,5.8){$\hat{a}$}

\psdots[dotstyle=*](8.,4.)
\rput(7.6,4.2){$a$}

\psline(6.5,3.5)(3.5,6.5)
\psline(6.,6.)(8.,4.)
\psline(3.5,6.)(6.,6.)
\psline(3.5,4.)(8.,4.)
\psline(4.,3.5)(4.,6.5)
\psline(6.,3.5)(6.,6.5)
%
\rput(4.7,4.7){$\Delta$}
\rput(6.9,4.7){$\hat{\Delta}$}
%
\psline{->}(4.7,5.)(6.4,5.)
\rput(6.6,5.){$\frac{\pi}{6}$}

\rput(6.0,2.8){D (i)}
\end{pspicture*}

&

\begin{pspicture*}(1.,2.5)(9.,7.)
\psdots[dotstyle=*](4.,4.)
\rput(3.7,4.3){$\hat{b}$}
\rput(4.2,4.2){$\lambda$}
\rput(4.2,3.8){$\hat{b}$}
\rput(3.7,3.8){$\lambda$}

\psdots[dotstyle=*](4.,6.)
\rput(3.5,6.2){$b$}
\rput(3.8,6.4){$\mu$}
\rput(4.2,6.2){$d$}
\rput(4.5,5.8){$\hat{b}$}
\rput(4.2,5.6){$a$}
\rput(3.7,5.8){$\mu$}

\psdots[dotstyle=*](6.,4.)
\rput(5.6,4.2){$b$}
\rput(5.8,4.5){$\mu$}
\rput(6.2,4.2){$a$}
\rput(5.7,3.8){$\hat{a}$}

\psdots[dotstyle=*](6.,6.)
\rput(6.2,5.7){$b$}
\rput(5.7,5.8){$\hat{a}$}

\psline(6.5,3.5)(3.5,6.5)
\psline(3.5,6.)(6.5,6.)
\psline(3.5,4.)(6.5,4.)
\psline(4.,3.5)(4.,6.5)
\psline(6.,3.5)(6.,6.)
%
\rput(4.7,4.7){$\Delta$}
\rput(5.3,5.3){$\hat{\Delta}$}
%
\psline{->}(5.0,4.6)(5.4,5.0)
\rput(5.6,4.9){$\frac{\pi}{6}$}

\rput(6.0,2.8){D (ii)}
\end{pspicture*}
\\
\end{tabular}
\end{center}

\pagebreak

\begin{center}
\begin{tabular}{c@{}c}

\begin{pspicture*}(1.,0.5)(9.,11.)
\psdots[dotstyle=*](4.,4.)
\rput(3.7,4.){$v_1$}
\rput(4.2,4.2){$\lambda$}
\rput(4.4,3.8){$\hat{a}$}
\rput(4.2,3.6){$b$}

\psdots[dotstyle=*](4.,6.)
\rput(4.2,6.2){$a$}
\rput(4.4,5.8){$\mu$}
\rput(4.2,5.6){$a$}
\rput(3.7,6.){$v_2$}

\psdots[dotstyle=*](4.,8.)
\rput(3.7,8.){$v_3$}
\rput(4.2,7.6){$b$}

\psdots[dotstyle=*](6.,2.)
\rput(5.8,2.0){$\lambda$}
\rput(5.8,2.5){$\hat{b}$}
\rput(6.2,2.2){$\lambda$}
\rput(6.2,1.8){$\hat{b}$}

\psdots[dotstyle=*](6.,4.)
\rput(5.6,4.2){$b$}
\rput(5.8,4.5){$\hat{a}$}
\rput(6.2,4.2){$\lambda$}
\rput(6.5,3.8){$\hat{b}$}
\rput(6.2,3.6){$a$}
\rput(5.7,3.8){$\mu$}

\psdots[dotstyle=*](6.,6.)
\rput(5.6,6.2){$\lambda$}
\rput(6.2,6.1){$\hat{e}$}
\rput(6.2,5.6){$a$}
\rput(5.7,5.8){$\hat{b}$}

\psdots[dotstyle=*](8.,2.)
\rput(7.6,2.2){$b$}
\rput(7.8,2.4){$\mu$}
\rput(8.2,2.2){$a$}
\rput(7.7,1.8){$\hat{a}$}

\psdots[dotstyle=*](8.,4.)
\rput(7.6,4.2){$b$}
\rput(8.2,4.1){$\hat{d}$}
\rput(8.2,3.6){$b$}
\rput(7.7,3.8){$\hat{a}$}

\psline(6.,2.)(4.,4.)
\psline(8.5,1.5)(4.,6.)
\psline(8.5,3.5)(4.,8.)
\psline(4.,6.)(6.,6.)
\psline(4.,4.)(8.,4.)
\psline(6.,2.)(8.5,2.)
\psline(4.,3.5)(4.,8.)
\psline(6.,1.5)(6.,6.)
\psline(8.,1.5)(8.,4.)
%
\rput(6.7,2.7){$\Delta$}
\rput(5.3,5.7){$\hat{\Delta}_2$}
\rput(4.3,5.2){$\hat{\Delta}_1$}
%
\psbezier{->}
(6.7,2.9)(6.7,4.5)
(7.,5.)(5.6,5.0)
\rput(5.5,5.2){$\frac{\pi}{6}$}

\psline[linestyle=dotted]{->}(5.6,5.0)(4.8,4.5)
\rput(4.7,4.6){$\frac{\pi}{6}$}

\rput(6.0,1.0){D (ix)}
\end{pspicture*}

&

\begin{pspicture*}(1.,0.5)(9.,11.)
\psdots[dotstyle=*](4.,4.)
\rput(3.6,4.){$\geq 6$}
\rput(4.2,4.2){$\lambda$}
\rput(4.4,3.8){$\hat{a}$}
\rput(4.2,3.6){$b$}

\psdots[dotstyle=*](4.,6.)
\rput(4.2,6.2){$a$}
\rput(4.4,5.8){$\mu$}
\rput(4.2,5.6){$a$}
\rput(3.6,6.){$\geq 8$}

\psdots[dotstyle=*](4.,8.)
\rput(3.5,8.){$\encircle{\geq 6}$}
\rput(4.2,7.6){$b$}

\psdots[dotstyle=*](6.,2.)
\rput(5.8,2.0){$\lambda$}
\rput(5.8,2.5){$\hat{b}$}
\rput(6.2,2.2){$\lambda$}
\rput(6.2,1.8){$\hat{b}$}

\psdots[dotstyle=*](6.,4.)
\rput(5.6,4.2){$b$}
\rput(5.8,4.5){$\hat{a}$}
\rput(6.2,4.2){$\lambda$}
\rput(6.5,3.8){$\hat{b}$}
\rput(6.2,3.6){$a$}
\rput(5.7,3.8){$\mu$}

\psdots[dotstyle=*](6.,6.)
\rput(5.6,6.2){$\lambda$}
\rput(6.2,6.1){$\hat{e}$}
\rput(6.2,5.6){$a$}
\rput(5.7,5.8){$\hat{b}$}

\psdots[dotstyle=*](8.,2.)
\rput(7.6,2.2){$b$}
\rput(7.8,2.4){$\mu$}
\rput(8.2,2.2){$a$}
\rput(7.7,1.8){$\hat{a}$}

\psdots[dotstyle=*](8.,4.)
\rput(7.6,4.2){$b$}
\rput(8.2,4.1){$\hat{d}$}
\rput(8.2,3.6){$b$}
\rput(7.7,3.8){$\hat{a}$}

\psline(6.,2.)(4.,4.)
\psline(8.5,1.5)(4.,6.)
\psline(8.5,3.5)(4.,8.)
\psline(4.,6.)(6.,6.)
\psline(4.,4.)(8.,4.)
\psline(6.,2.)(8.5,2.)
\psline(4.,3.5)(4.,8.)
\psline(6.,1.5)(6.,6.)
\psline(8.,1.5)(8.,4.)
%
\rput(6.7,2.7){$\Delta$}
\rput(4.5,7.){$\hat{\Delta}$}
\rput(5.3,5.7){$\Delta_1$}
\rput(7.,4.7){$\hat{\Delta}_1$}
%
\psbezier{->}
(6.7,2.9)(6.7,7.)
(5.,4)(4.7,6.3)
\rput(4.9,6.3){$\frac{\pi}{6}$}

\psline{->}(5.7,4.9)(6.3,4.9)
\rput(6.2,4.6){$\frac{\pi}{12}$}

\rput(6.0,1.0){D (x)}
\end{pspicture*}
\\
\end{tabular}
\end{center}

\pagebreak

\begin{center}
\begin{tabular}{c@{}c}

\begin{pspicture*}(1.,0.)(9.,11.)
\psdots[dotstyle=*](2.,6.)
\rput(2.2,6.2){$b$}
\psdots[dotstyle=*](2.,8.)
\rput(2.2,7.6){$\lambda$}
\rput(2.5,7.8){$\hat{b}$}
\rput(2.2,8.2){$\lambda$}
\psdots[dotstyle=*](4.,4.)
\rput(3.6,4.){$\geq{6}$}
\rput(4.2,4.2){$\lambda$}
\rput(4.5,3.8){$\hat{a}$}
\rput(4.2,3.6){$b$}

\psdots[dotstyle=*](4.,6.)
\rput(3.5,6.2){$a$}
\rput(3.8,6.4){$\mu$}
\rput(4.2,6.2){$a$}
\rput(4.4,5.8){$\mu$}
\rput(4.2,5.6){$a$}
\psline(4.,6.)(3.1,5.3)
\rput(3.4,5.8){$\hat{d}$}
\psline(4.,6.)(3.6,5.3)
\rput(3.5,5.4){$a$}
\rput(3.85,5.4){$\mu$}
\rput(4.4,5.8){$\mu$}

\psdots[dotstyle=*](4.,8.)
\rput(3.6,8.2){$b$}
\rput(4.2,8.2){$\hat{d}$}
\rput(4.2,7.6){$b$}
\rput(3.7,7.8){$\hat{a}$}

\psdots[dotstyle=*](6.,2.)
\rput(5.8,2.0){$\lambda$}
\rput(5.8,2.5){$\hat{b}$}
\rput(6.2,2.2){$\lambda$}
\rput(6.2,1.8){$\hat{b}$}

\psdots[dotstyle=*](6.,4.)
\rput(5.6,4.2){$b$}
\rput(5.8,4.5){$\hat{a}$}
\rput(6.2,4.2){$\lambda$}
\rput(6.45,3.8){$\hat{b}$}
\rput(6.2,3.6){$a$}
\rput(5.7,3.8){$\mu$}

\psdots[dotstyle=*](6.,6.)
\rput(5.6,6.2){$\lambda$}
\rput(6.2,6.1){$\hat{e}$}
\rput(6.2,5.6){$a$}
\rput(5.7,5.8){$\hat{b}$}

\psdots[dotstyle=*](8.,2.)
\rput(7.6,2.2){$b$}
\rput(7.8,2.4){$\mu$}
\rput(8.2,2.2){$a$}
\rput(7.7,1.8){$\hat{a}$}

\psdots[dotstyle=*](8.,4.)
\rput(7.6,4.2){$b$}
\rput(8.2,4.1){$\hat{d}$}
\rput(8.2,3.6){$b$}
\rput(7.7,3.8){$\hat{a}$}

\psline(6.,2.)(4.,4.)
\psline(8.5,1.5)(2.,8.)
\psline(8.5,3.5)(3.5,8.5)
\psline(2.,8.)(4.,8.)
\psline(2.,6.)(6.,6.)
\psline(4.,4.)(8.,4.)
\psline(6.,2.)(8.5,2.)
\psline(2.,6.)(2.,8.5)
\psline(4.,3.5)(4.,8.)
\psline(6.,1.5)(6.,6.)
\psline(8.,1.5)(8.,4.)
%
\rput(6.7,2.7){$\Delta$}
\rput(4.7,6.7){$\hat{\Delta}_1$}
\rput(3.3,7.3){$\Delta_1$}
\rput(2.7,6.7){$\hat{\Delta}$}
%
\psbezier{->}
(6.9,2.9)(6.0,5.8)
(5.,6.7)(2.9,6.7)
\rput(3.1,6.4){$\frac{\pi}{6}$}
\psline{->}(3.7,7.0)(4.4,6.7)
\rput(4.3,7.1){$\frac{\pi}{4}$}

\rput(6.0,1.0){D (xiii)}
\end{pspicture*}
&
\begin{pspicture*}(1.,0.)(9.,11.)
\psdots[dotstyle=*](2.,8.)
\rput(1.8,8.){$u$}
\rput(2.5,7.8){$\hat{b}$}
\psdots[dotstyle=*](4.,4.)
\rput(3.6,3.5){$\geq{6}$}
\rput(4.2,4.2){$\lambda$}
\rput(4.5,3.8){$\hat{a}$}
\rput(3.4,4.1){$\hat{e}/\hat{c}/\hat{d}$}
\rput(4.2,3.6){$b$}
\psline(4.,4.)(3.0,3.5)

\psdots[dotstyle=*](4.,6.)
\rput(3.8,6.4){$\mu$}
\rput(4.2,6.2){$a$}
\rput(4.4,5.8){$\mu$}
\rput(4.2,5.6){$a$}
\rput(3.3,6.1){$=8$}
\rput(3.5,5.3){$\hat{d}/\hat{c}/\hat{e}$}
\psline(4.,6.)(3.0,5.5)

\psdots[dotstyle=*](4.,8.)
\rput(3.6,8.2){$b$}
\rput(4.2,8.2){$\hat{d}$}
\rput(4.2,7.6){$b$}
\rput(3.7,7.8){$\hat{a}$}

\psdots[dotstyle=*](6.,2.)
\rput(5.8,2.0){$\lambda$}
\rput(5.8,2.5){$\hat{b}$}
\rput(6.2,2.2){$\lambda$}
\rput(6.2,1.8){$\hat{b}$}

\psdots[dotstyle=*](6.,4.)
\rput(5.6,4.2){$b$}
\rput(5.8,4.5){$\hat{a}$}
\rput(6.2,4.2){$\lambda$}
\rput(6.5,3.8){$\hat{b}$}
\rput(6.2,3.6){$a$}
\rput(5.7,3.8){$\mu$}

\psdots[dotstyle=*](6.,6.)
\rput(5.6,6.2){$\lambda$}
\rput(6.2,6.1){$\hat{e}$}
\rput(6.2,5.6){$a$}
\rput(5.7,5.8){$\hat{b}$}

\psdots[dotstyle=*](8.,2.)
\rput(7.6,2.2){$b$}
\rput(7.8,2.4){$\mu$}
\rput(8.2,2.2){$a$}
\rput(7.7,1.8){$\hat{a}$}

\psdots[dotstyle=*](8.,4.)
\rput(7.6,4.2){$b$}
\rput(8.2,4.1){$\hat{d}$}
\rput(8.2,3.6){$b$}
\rput(7.7,3.8){$\hat{a}$}

\psline(6.,2.)(4.,4.)
\psline(8.5,1.5)(2.,8.)
\psline(8.5,3.5)(3.5,8.5)
\psline(2.,8.)(4.,8.)
\psline(4.,6.)(6.,6.)
\psline(4.,4.)(8.,4.)
\psline(6.,2.)(8.5,2.)
\psline(4.,3.5)(4.,8.)
\psline(6.,1.5)(6.,6.)
\psline(8.,1.5)(8.,4.)
%
\rput(6.7,2.7){$\Delta$}
\rput(5.3,5.3){$\Delta_1$}
\rput(4.7,4.7){$\hat{\Delta}$}
%
\psbezier{->}
(6.7,2.9)(7.,4.5)
(6.,5.)(5.,4.5)
\rput(5.,4.25){$\frac{\pi}{6}$}

\psline{->}(5.5,5.5)(6.4,5.2)
\rput(6.5,5.){$\frac{\pi}{24}$}
\psline{->}(5.5,5.5)(5.2,6.3)
\rput(5.0,6.4){$\frac{\pi}{24}$}

\rput(6.0,1.0){D (xiv)}
\end{pspicture*}
\\
\end{tabular}

\end{center}

\newpage
\subsection*{Figure {{E}}}
\psset{arrowscale=0.7,linewidth=1.0pt,dotsize=3pt}
\psset{xunit=1.0cm,yunit=1.0cm,}

\begin{center}
\begin{tabular}{c@{}c}

\begin{pspicture*}(2.5,0.5)(10.,11.)

\psline(4.,4.)(4.,10.)
\psline(6.,1.5)(6.,10.)
\psline(8.,1.5)(8.,10.)
\psline(6.,2.)(4.,4.)
\psline(8.5,1.5)(4.,6.)
\psline(8.5,3.5)(4.,8.)
\psline(8.5,5.5)(4.,10.)
\psline(8.5,7.5)(6.,10.)
\psline(8.5,9.5)(8.,10.)

\psline(6.,2.)(8.5,2.)

\psdots[dotstyle=*](6.,2.)
\rput(5.8,2.0){$\lambda$}
\rput(5.8,2.5){$\hat{b}$}
\rput(6.2,2.2){$\lambda$}
\rput(6.2,1.8){$\hat{b}$}

\psdots[dotstyle=*](8.,2.)
\rput(7.6,2.2){$b$}
\rput(7.8,2.4){$\mu$}
\rput(8.2,2.2){$a$}
\rput(7.7,1.8){$\hat{a}$}

%
\psline(4.,4.)(8.5,4.)
\psdots[dotstyle=*](4.,4.)
\rput(4.2,4.2){$\lambda$}
\rput(4.45,3.8){$\hat{a}$}

\psdots[dotstyle=*](6.,4.)
\rput(5.6,4.2){$b$}
\rput(5.8,4.5){$\hat{a}$}
\rput(6.2,4.2){$\lambda$}
\rput(6.5,3.8){$\hat{b}$}
\rput(6.2,3.6){$a$}
\rput(5.7,3.8){$\mu$}

\psdots[dotstyle=*](8.,4.)
\rput(7.6,4.2){$b$}
\rput(7.85,4.4){$\mu$}
\rput(8.2,4.2){$a$}
\rput(8.4,3.85){$\mu$}
\rput(8.2,3.6){$b$}
\rput(7.7,3.8){$\hat{a}$}

\psline(4.,6.)(8.5,6.)

\psdots[dotstyle=*](4.,6.)
\rput(3.7,6.){$v_1$}
\rput(4.2,6.2){$a$}
\rput(4.45,5.8){$\mu$}
\rput(4.2,5.6){$a$}

\psdots[dotstyle=*](6.,6.)
\rput(5.6,6.2){$\lambda$}
\rput(5.8,6.5){$\hat{a}$}
\rput(6.2,6.2){$\lambda$}
\rput(6.5,5.8){$\hat{b}$}
\rput(6.2,5.6){$a$}
\rput(5.7,5.8){$\hat{b}$}

\psdots[dotstyle=*](8.,6.)
\rput(7.6,6.2){$b$}
\rput(7.85,6.4){$\mu$}
\rput(8.2,6.2){$a$}
\rput(8.4,5.85){$\mu$}
\rput(8.2,5.6){$b$}
\rput(7.7,5.8){$\hat{a}$}

\psline(4.,8.)(8.5,8.)

\psdots[dotstyle=*](4.,8.)
\rput(3.7,8.){$v_2$}
\rput(4.2,8.2){$a$}
\rput(4.45,7.8){$\mu$}
\rput(4.2,7.6){$b$}

\psdots[dotstyle=*](6.,8.)
\rput(5.6,8.2){$\lambda$}
\rput(5.8,8.5){$\hat{a}$}
\rput(6.2,8.2){$\lambda$}
\rput(6.5,7.8){$\hat{b}$}
\rput(6.2,7.6){$a$}
\rput(5.7,7.8){$\hat{b}$}

\psdots[dotstyle=*](8.,8.)
\rput(7.6,8.2){$b$}
\rput(7.85,8.4){$\mu$}
\rput(8.2,8.2){$a$}
\rput(8.4,7.85){$\mu$}
\rput(8.2,7.6){$b$}
\rput(7.7,7.8){$\hat{a}$}

\psline(4.,10.)(8.,10.)

\psdots[dotstyle=*](4.,10.)
\rput(4.45,9.8){$\mu$}
\rput(4.2,9.6){$b$}

\psdots[dotstyle=*](6.,10.)
\rput(6.,10.3){$u_8$}
\rput(6.5,9.8){$\hat{b}$}
\rput(6.2,9.6){$a$}
\rput(5.7,9.8){$\hat{b}$}

\psdots[dotstyle=*](8.,10.)
\rput(8.,10.3){$u_9$}
\rput(8.2,9.6){$b$}
\rput(7.7,9.8){$\hat{a}$}

\rput(6.5,2.5){$\Delta$}
\rput(6.5,8.5){$\hat{\Delta}_1$}
\rput(7.5,9.5){$\hat{\Delta}_2$}

\rput(5.2,7.5){$\hat{\Delta}_3$}
\rput(5.2,5.5){$\hat{\Delta}_4$}

\psline{->}(7.,2.5)(7.,8.7)
\rput(6.8,8.8){$\frac{\pi}{6}$}
\psline[linestyle=dotted]{->}(7.,8.2)(7.5,9.)
\rput(7.7,9.){$\frac{\pi}{6}$}

\psline[linestyle=dotted]{->}(7.,6.7)(5.5,7.)
\rput(5.7,7.3){$\frac{\pi}{6}$}
\psline[linestyle=dotted]{->}(7.,4.7)(5.5,5.)
\rput(5.7,5.3){$\frac{\pi}{6}$}

\rput(6.0,1.0){E (i)}
\end{pspicture*}

&

\begin{pspicture*}(2.5,0.5)(10.,11.)

\psline(4.,4.)(4.,10.)
\psline(6.,1.5)(6.,10.)
\psline(8.,1.5)(8.,10.)
\psline(6.,2.)(4.,4.)
\psline(8.5,1.5)(4.,6.)
\psline(8.5,3.5)(4.,8.)
\psline(8.5,5.5)(4.,10.)
\psline(8.5,7.5)(6.,10.)
\psline(8.5,9.5)(8.,10.)

\psline(6.,2.)(8.5,2.)

\psdots[dotstyle=*](6.,2.)
\rput(5.8,2.0){$\lambda$}
\rput(5.8,2.5){$\hat{b}$}
\rput(6.2,2.2){$\lambda$}
\rput(6.2,1.8){$\hat{b}$}

\psdots[dotstyle=*](8.,2.)
\rput(7.6,2.2){$b$}
\rput(7.8,2.4){$\mu$}
\rput(8.2,2.2){$a$}
\rput(7.7,1.8){$\hat{a}$}

%
\psline(4.,4.)(8.5,4.)
\psdots[dotstyle=*](4.,4.)
\rput(4.2,4.2){$\lambda$}
\rput(4.45,3.8){$\hat{a}$}

\psdots[dotstyle=*](6.,4.)
\rput(5.6,4.2){$b$}
\rput(5.8,4.5){$\hat{a}$}
\rput(6.2,4.2){$\lambda$}
\rput(6.5,3.8){$\hat{b}$}
\rput(6.2,3.6){$a$}
\rput(5.7,3.8){$\mu$}

\psdots[dotstyle=*](8.,4.)
\rput(7.6,4.2){$b$}
\rput(7.85,4.4){$\mu$}
\rput(8.2,4.2){$a$}
\rput(8.4,3.85){$\mu$}
\rput(8.2,3.6){$b$}
\rput(7.7,3.8){$\hat{a}$}

\psline(4.,6.)(8.5,6.)

\psdots[dotstyle=*](4.,6.)
\rput(4.2,6.2){$a$}
\rput(4.45,5.8){$\mu$}
\rput(4.2,5.6){$a$}

\psdots[dotstyle=*](6.,6.)
\rput(5.6,6.2){$\lambda$}
\rput(5.8,6.5){$\hat{a}$}
\rput(6.2,6.2){$\lambda$}
\rput(6.5,5.8){$\hat{b}$}
\rput(6.2,5.6){$a$}
\rput(5.7,5.8){$\hat{b}$}

\psdots[dotstyle=*](8.,6.)
\rput(7.6,6.2){$b$}
\rput(7.85,6.4){$\mu$}
\rput(8.2,6.2){$a$}
\rput(8.4,5.85){$\mu$}
\rput(8.2,5.6){$b$}
\rput(7.7,5.8){$\hat{a}$}

\psline(4.,8.)(8.5,8.)

\psdots[dotstyle=*](4.,8.)
\rput(4.2,8.2){$a$}
\rput(4.45,7.8){$\mu$}
\rput(4.2,7.6){$b$}

\psdots[dotstyle=*](6.,8.)
\rput(5.6,8.2){$\lambda$}
\rput(5.8,8.5){$\hat{a}$}
\rput(6.2,8.2){$\lambda$}
\rput(6.5,7.8){$\hat{b}$}
\rput(6.2,7.6){$a$}
\rput(5.7,7.8){$\hat{b}$}

\psdots[dotstyle=*](8.,8.)
\rput(7.6,8.2){$b$}
\rput(7.85,8.4){$\mu$}
\rput(8.2,8.2){$a$}
\rput(8.4,7.85){$\mu$}
\rput(8.2,7.6){$b$}
\rput(7.7,7.8){$\hat{a}$}

\psline(4.,10.)(8.,10.)

\psdots[dotstyle=*](4.,10.)
\rput(4.45,9.8){$\mu$}
\rput(4.2,9.6){$b$}

\psdots[dotstyle=*](6.,10.)
\rput(6.,10.3){$\geq 8$}
\rput(6.5,9.8){$\hat{b}$}
\rput(6.2,9.6){$a$}
\rput(5.7,9.8){$\hat{b}$}

\psdots[dotstyle=*](8.,10.)
\rput(8.,10.3){$\geq 6$}
\rput(8.2,9.6){$b$}
\rput(7.7,9.8){$\hat{a}$}

\rput(6.5,2.5){$\Delta$}
\rput(6.5,6.5){$\Delta_0$}
\rput(6.5,8.5){$\hat{\Delta}_1$}
\rput(5.5,9.3){$\Delta_1$}
\rput(7.5,9.5){$\hat{\Delta}_2$}
\rput(8.5,8.7){$\Delta_2$}

\psline{->}(7.,2.5)(7.,8.7)
\rput(6.8,8.8){$\frac{\pi}{12}$}
\psline{->}(7.,8.2)(7.5,9.)
\rput(7.7,9.){$\frac{\pi}{12}$}

\rput(6.0,1.0){E (ii)}
\end{pspicture*}
\\
\end{tabular}
\end{center}

\pagebreak

\begin{center}
\begin{tabular}{c@{}c}

\begin{pspicture*}(2.5,0.5)(10.,13.)

\psline(4.,4.)(4.,10.)
\psline(6.,1.5)(6.,10.)
\psline(8.,1.5)(8.,10.)
\psline(6.,2.)(4.,4.)
\psline(8.5,1.5)(4.,6.)
\psline(8.5,3.5)(4.,8.)
\psline(8.5,5.5)(4.,10.)
\psline(8.5,7.5)(6.,10.)
\psline(8.5,9.5)(8.,10.)

\psline(6.,2.)(8.5,2.)

\psdots[dotstyle=*](6.,2.)
\rput(5.8,2.0){$\lambda$}
\rput(5.8,2.5){$\hat{b}$}
\rput(6.2,2.2){$\lambda$}
\rput(6.2,1.8){$\hat{b}$}

\psdots[dotstyle=*](8.,2.)
\rput(7.6,2.2){$b$}
\rput(7.8,2.4){$\mu$}
\rput(8.2,2.2){$a$}
\rput(7.7,1.8){$\hat{a}$}

%
\psline(4.,4.)(8.5,4.)
\psdots[dotstyle=*](4.,4.)
\rput(4.2,4.2){$\lambda$}
\rput(4.45,3.8){$\hat{a}$}

\psdots[dotstyle=*](6.,4.)
\rput(5.6,4.2){$b$}
\rput(5.8,4.5){$\hat{a}$}
\rput(6.2,4.2){$\lambda$}
\rput(6.5,3.8){$\hat{b}$}
\rput(6.2,3.6){$a$}
\rput(5.7,3.8){$\mu$}

\psdots[dotstyle=*](8.,4.)
\rput(7.6,4.2){$b$}
\rput(7.85,4.4){$\mu$}
\rput(8.2,4.2){$a$}
\rput(8.4,3.85){$\mu$}
\rput(8.2,3.6){$b$}
\rput(7.7,3.8){$\hat{a}$}

\psline(4.,6.)(8.5,6.)

\psdots[dotstyle=*](4.,6.)
\rput(4.2,6.2){$a$}
\rput(4.45,5.8){$\mu$}
\rput(4.2,5.6){$a$}

\psdots[dotstyle=*](6.,6.)
\rput(5.6,6.2){$\lambda$}
\rput(5.8,6.5){$\hat{a}$}
\rput(6.2,6.2){$\lambda$}
\rput(6.5,5.8){$\hat{b}$}
\rput(6.2,5.6){$a$}
\rput(5.7,5.8){$\hat{b}$}

\psdots[dotstyle=*](8.,6.)
\rput(7.6,6.2){$b$}
\rput(7.85,6.4){$\mu$}
\rput(8.2,6.2){$a$}
\rput(8.4,5.85){$\mu$}
\rput(8.2,5.6){$b$}
\rput(7.7,5.8){$\hat{a}$}

\psline(4.,8.)(8.5,8.)

\psdots[dotstyle=*](4.,8.)
\rput(4.2,8.2){$a$}
\rput(4.45,7.8){$\mu$}
\rput(4.2,7.6){$b$}

\psdots[dotstyle=*](6.,8.)
\rput(5.6,8.2){$\lambda$}
\rput(5.8,8.5){$\hat{a}$}
\rput(6.2,8.2){$\lambda$}
\rput(6.5,7.8){$\hat{b}$}
\rput(6.2,7.6){$a$}
\rput(5.7,7.8){$\hat{b}$}

\psdots[dotstyle=*](8.,8.)
\rput(7.6,8.2){$b$}
\rput(7.85,8.4){$\mu$}
\rput(8.2,8.2){$a$}
\rput(8.4,7.85){$\mu$}
\rput(8.2,7.6){$b$}
\rput(7.7,7.8){$\hat{a}$}

\psline(4.,10.)(8.,10.)

\psdots[dotstyle=*](4.,10.)
\rput(3.7,10.){$u$}
\rput(4.0,10.5){$\encircle{\geq 8}$}
\rput(4.45,9.8){$\mu$}
\rput(4.2,9.6){$b$}

\psdots[dotstyle=*](6.,10.)
\rput(6.,10.3){$d$}
\rput(6.5,9.8){$\hat{b}$}
\rput(6.2,9.6){$a$}
\rput(5.7,9.8){$\hat{b}$}

\psdots[dotstyle=*](8.,10.)
\rput(8.2,9.6){$b$}
\rput(7.7,9.8){$\hat{a}$}

\rput(6.5,2.5){$\Delta$}
\rput(5.5,9.5){$\hat{\Delta}$}

\psbezier{->}%
(6.7,2.7)(7.5,4.5)
(7.5,9.)(5.7,9.)
\rput(5.5,9.){$\frac{\pi}{6}$}

\rput(6.0,1.0){E (ix)}
\end{pspicture*}

&

\begin{pspicture*}(1.8,0.5)(10.,13.)

\psline(4.,4.)(4.,10.)
\psline(6.,1.5)(6.,10.)
\psline(8.,1.5)(8.,10.)
\psline(6.,2.)(4.,4.)
\psline(8.5,1.5)(4.,6.)
\psline(8.5,3.5)(2.,10.)
\psline(8.5,5.5)(4.,10.)
\psline(8.5,7.5)(6.,10.)
\psline(8.5,9.5)(8.,10.)

\psline(6.,2.)(8.5,2.)

\psdots[dotstyle=*](6.,2.)
\rput(5.8,2.0){$\lambda$}
\rput(5.8,2.5){$\hat{b}$}
\rput(6.2,2.2){$\lambda$}
\rput(6.2,1.8){$\hat{b}$}

\psdots[dotstyle=*](8.,2.)
\rput(7.6,2.2){$b$}
\rput(7.8,2.4){$\mu$}
\rput(8.2,2.2){$a$}
\rput(7.7,1.8){$\hat{a}$}

%
\psline(4.,4.)(8.5,4.)
\psdots[dotstyle=*](4.,4.)
\rput(4.2,4.2){$\lambda$}
\rput(4.45,3.8){$\hat{a}$}

\psdots[dotstyle=*](6.,4.)
\rput(5.6,4.2){$b$}
\rput(5.8,4.5){$\hat{a}$}
\rput(6.2,4.2){$\lambda$}
\rput(6.5,3.8){$\hat{b}$}
\rput(6.2,3.6){$a$}
\rput(5.7,3.8){$\mu$}

\psdots[dotstyle=*](8.,4.)
\rput(7.6,4.2){$b$}
\rput(7.85,4.4){$\mu$}
\rput(8.2,4.2){$a$}
\rput(8.4,3.85){$\mu$}
\rput(8.2,3.6){$b$}
\rput(7.7,3.8){$\hat{a}$}

\psline(4.,6.)(8.5,6.)

\psdots[dotstyle=*](4.,6.)
\rput(3.7,6.){$v_1$}
\rput(4.2,6.2){$a$}
\rput(4.45,5.8){$\mu$}
\rput(4.2,5.6){$a$}

\psdots[dotstyle=*](6.,6.)
\rput(5.6,6.2){$\lambda$}
\rput(5.8,6.5){$\hat{a}$}
\rput(6.2,6.2){$\lambda$}
\rput(6.5,5.8){$\hat{b}$}
\rput(6.2,5.6){$a$}
\rput(5.7,5.8){$\hat{b}$}

\psdots[dotstyle=*](8.,6.)
\rput(7.6,6.2){$b$}
\rput(7.85,6.4){$\mu$}
\rput(8.2,6.2){$a$}
\rput(8.4,5.85){$\mu$}
\rput(8.2,5.6){$b$}
\rput(7.7,5.8){$\hat{a}$}

\psline(4.,8.)(8.5,8.)

\psdots[dotstyle=*](4.,8.)
\rput(3.7,8.){$v_2$}
\rput(3.8,8.5){$\mu$}
\rput(4.2,8.2){$a$}
\rput(4.45,7.8){$\mu$}
\rput(4.2,7.6){$b$}

\psdots[dotstyle=*](6.,8.)
\rput(5.6,8.2){$\lambda$}
\rput(5.8,8.5){$\hat{a}$}
\rput(6.2,8.2){$\lambda$}
\rput(6.5,7.8){$\hat{b}$}
\rput(6.2,7.6){$a$}
\rput(5.7,7.8){$\hat{b}$}

\psdots[dotstyle=*](8.,8.)
\rput(7.6,8.2){$b$}
\rput(7.85,8.4){$\mu$}
\rput(8.2,8.2){$a$}
\rput(8.4,7.85){$\mu$}
\rput(8.2,7.6){$b$}
\rput(7.7,7.8){$\hat{a}$}

\psline(2.,10.)(8.,10.)

\psdots[dotstyle=*](2.,10.)
\rput(2.45,9.8){$\hat{b}$}

\psdots[dotstyle=*](4.,10.)
\rput(3.8,9.8){$\hat{a}$}
\rput(4.0,10.3){$e$}
\rput(4.45,9.8){$\mu$}
\rput(4.2,9.6){$b$}

\psdots[dotstyle=*](6.,10.)
\rput(6.,10.3){$d$}
\rput(6.5,9.8){$\hat{b}$}
\rput(6.2,9.6){$a$}
\rput(5.7,9.8){$\hat{b}$}

\psdots[dotstyle=*](8.,10.)
\rput(8.2,9.6){$b$}
\rput(7.7,9.8){$\hat{a}$}

\rput(6.5,2.5){$\Delta$}
\rput(4.5,7.){$\hat{\Delta}_1$}
\rput(5.7,5.3){$\hat{\Delta}_2$}

\psbezier{->}%
(7.,2.5)(7.5,10.5)
(4.8,10.5)(4.8,6.2)
\rput(4.5,6.5){$\frac{\pi}{12}$}
\psline{->}(4.8,6.7)(5.3,5.6)
\rput(5.2,5.3){$\frac{\pi}{12}$}

\rput(6.0,1.0){E (x)}
\end{pspicture*}
\\
\end{tabular}
\end{center}

\newpage
\subsection*{Figure {{F}}}
\psset{arrowscale=0.7,linewidth=1.0pt,dotsize=3pt}
\psset{xunit=1.0cm,yunit=1.0cm,}
\begin{center}
\begin{tabular}{c@{}c}
\begin{pspicture*}(1.,2.5)(9.,11.)
\psdots[dotstyle=*](2.,8.)
\rput(2.5,7.8){$\mu$}

\psdots[dotstyle=*](4.,4.)
\rput(4.2,4.2){$\lambda$}
\rput(4.2,3.8){$\hat{b}$}

\psdots[dotstyle=*](4.,6.)
\rput(3.85,6.4){$\hat{a}$}
\rput(4.2,6.2){$\lambda$}
\rput(4.5,5.8){$\hat{b}$}
\rput(4.2,5.6){$a$}
\rput(3.3,6.2){$b/\lambda$}
\rput(3.6,5.7){$\mu/\hat{b}$}

\psdots[dotstyle=*](4.,8.)
\rput(3.7,8.5){$\encircle{\geq 6}$}
\rput(3.8,7.7){$\hat{b}$}
\rput(4.2,7.6){$a$}
\rput(4.3,8.0){$\mu$}

\psdots[dotstyle=*](6.,4.)
\rput(5.6,4.2){$b$}
\rput(5.8,4.5){$\mu$}
\rput(6.2,4.2){$a$}
\rput(6.45,3.8){$\mu$}
\rput(6.2,3.6){$b$}
\rput(5.7,3.8){$\hat{a}$}

\psdots[dotstyle=*](6.,6.)
\rput(5.6,6.2){$b$}
\rput(5.65,6.6){$\hat{a}$}
\rput(6.,6.8){$d$}
\rput(6.35,6.6){$\hat{a}$}
\rput(6.5,6.2){$\lambda$}
\rput(5.7,5.8){$\hat{a}$}
\rput(6.2,5.6){$b$}
\rput(6.5,5.8){$\hat{a}$}

\psline(6.4,7.3)(6.,6.)

\psdots[dotstyle=*](5.,9.)
\rput(4.9,8.6){$\hat{b}$}
\psline(5.,9.)(4.,8.)
\psline(5.,9.)(6.,6.)

\psdots[dotstyle=*](8.,6.)
\rput(7.8,5.8){$\mu$}

\psdots[dotstyle=*](8.,4.)
\rput(7.6,4.2){$\lambda$}
\rput(7.85,4.4){$\hat{b}$}
\rput(7.8,3.8){$\hat{b}$}
\rput(8.2,4.1){$\lambda$}

\psline(6.5,3.5)(2.,8.)
\psline(8.5,3.5)(4.,8.)
\psline(7.,7.)(6.,6.)
\psline(2.,8.)(4.,8.)
\psline(3.5,6.)(8.,6.)
\psline(4.,4.)(8.,4.)
\psline(4.,3.5)(4.,8.)
\psline(6.,3.5)(6.,6.)
\psline(8.,3.5)(8.,6.)
%
\rput(4.4,4.7){${\Delta}_0$}
\rput(5.5,5.5){$\hat{\Delta}_1$}
\rput(6.7,4.7){${\Delta}_1$}
\rput(4.3,7.2){$\hat{\Delta}_2$}
\rput(5.0,8.0){$\hat{\Delta}_3$}

\psline{->}(6.3,5.2)(5.7,5.2)
\rput(5.7,4.9){$\frac{\pi}{12}$}
\psline{->}(5.0,3.5)(5.0,6.5)
\rput(4.9,6.7){$\frac{\pi}{12}$}
\psline{->}(5.0,5.5)(5.5,7.)
\rput(5.3,7.2){$\frac{\pi}{12}$}
\psline[linestyle=dotted]{->}(4.7,3.5)(4.7,6.5)
\rput(4.5,6.7){$\frac{\pi}{6}$}

\rput(6.0,2.8){F (i)}
\end{pspicture*}
&

\begin{pspicture*}(1.,2.5)(9.,11.)
\psdots[dotstyle=*](2.,8.)
\rput(2.5,7.8){$\mu$}

\psdots[dotstyle=*](4.,4.)
\rput(4.2,4.2){$\lambda$}
\rput(4.2,3.8){$\hat{b}$}

\psdots[dotstyle=*](4.,6.)
\rput(3.85,6.4){$\hat{a}$}
\rput(4.2,6.2){$\lambda$}
\rput(4.5,5.8){$\hat{b}$}
\rput(4.2,5.6){$a$}
\rput(3.3,6.2){$b/\lambda$}
\rput(3.6,5.7){$\mu/\hat{b}$}

\psdots[dotstyle=*](4.,8.)
\rput(3.7,8.5){$\encircle{\geq 6}$}
\rput(3.8,7.7){$\hat{b}$}
\rput(4.2,7.6){$a$}
\rput(4.3,8.0){$\hat{b}$}

\psdots[dotstyle=*](6.,4.)
\rput(5.6,4.2){$b$}
\rput(5.8,4.5){$\mu$}
\rput(6.2,4.2){$a$}
\rput(6.45,3.8){$\mu$}
\rput(6.2,3.6){$b$}
\rput(5.7,3.8){$\hat{a}$}

\psdots[dotstyle=*](6.,6.)
\rput(5.6,6.2){$b$}
\rput(5.65,6.6){$\mu$}
\rput(6.,6.8){$a$}
\rput(6.35,6.65){$\hat{d}$}
\rput(6.5,6.2){$\lambda$}
\rput(5.7,5.8){$\hat{a}$}
\rput(6.2,5.6){$b$}
\rput(6.5,5.8){$\hat{a}$}

\psline(6.4,7.3)(6.,6.)

\psdots[dotstyle=*](5.,9.)
\rput(4.9,8.6){$\hat{a}$}
\rput(5.3,8.8){$b$}
\rput(5.,9.2){$u$}
\psline(5.,9.)(4.,8.)
\psline(5.,9.)(6.,6.)

\psdots[dotstyle=*](8.,6.)
\rput(7.8,5.8){$\mu$}

\psdots[dotstyle=*](8.,4.)
\rput(7.6,4.2){$\lambda$}
\rput(7.85,4.4){$\hat{b}$}
\rput(8.,3.7){$\hat{b}$}
\rput(8.2,4.1){$\lambda$}

\psline(6.5,3.5)(2.,8.)
\psline(8.5,3.5)(4.,8.)
\psline(7.,7.)(6.,6.)
\psline(2.,8.)(4.,8.)
\psline(3.5,6.)(8.,6.)
\psline(4.,4.)(8.,4.)
\psline(4.,3.5)(4.,8.)
\psline(6.,3.5)(6.,6.)
\psline(8.,4.)(8.,6.)
%
\rput(4.4,4.7){${\Delta}_0$}
\rput(5.5,5.5){$\hat{\Delta}_1$}
\rput(6.7,4.7){${\Delta}_1$}
\rput(4.3,7.2){$\hat{\Delta}_2$}
\rput(4.8,7.8){$\hat{\Delta}_3$}
\rput(4.3,8.8){$\hat{\Delta}_4$}
\rput(5.8,8.2){$\hat{\Delta}_5$}

\psline{->}(6.3,5.2)(5.7,5.2)
\rput(5.7,4.9){$\frac{\pi}{12}$}
\psline{->}(5.0,3.5)(5.0,6.5)
\rput(4.9,6.7){$\frac{\pi}{12}$}
\psline{->}(5.0,5.5)(5.5,7.)
\rput(5.3,7.2){$\frac{\pi}{12}$}

\psline[linestyle=dotted]{->}(5.5,7.1)(5.9,7.7)
\rput(6.,7.5){$\frac{\pi}{24}$}
\psline[linestyle=dotted]{->}(5.3,7.4)(4.6,8.9)
\rput(4.5,9.1){$\frac{\pi}{24}$}
\psline[linestyle=dotted]{->}(4.7,3.5)(4.7,6.5)
\rput(4.5,6.7){$\frac{\pi}{6}$}

\rput(6.0,2.8){F (ii)}
\end{pspicture*}
\\
\end{tabular}%
\end{center}

\newpage
\subsection*{Figure {{G}}}
\psset{arrowscale=0.7,linewidth=1.0pt,dotsize=3pt}
\psset{xunit=.8cm,yunit=1.0cm,}
\begin{center}

\begin{tabular}{@{}c@{}c@{}c}

\begin{pspicture*}(2.,-4.)(7.5,1.5)
\psdots[dotstyle=*](4.,0.)
\psdots[dotstyle=*](6.,0.)
\psdots[dotstyle=*](4.,-2.)
\psdots[dotstyle=*](6.,-2.)
\psline(4.,0)(6.,0.)
\psline(4.,-2.)(6.,-2.)
\psline(4.,0)(4.,-2.)
\psline(6.,0)(6.,-2.)
\psline(4.,0)(6.,-2.)
\psline(3.3,0.)(4.,0.)
\psline(3.3,-2.)(4.,-2.)
\psline(4.,-2.)(4.,-2.7)
\psline(6.,-2.)(6.,-2.7)
%
%
\rput(4.0,0.2){$u_2$}
\rput(6.0,0.2){$u_3$}
\rput(3.5,-0.2){$b$}
\rput(4.2,-0.5){$\hat{a}$}
\rput(3.5,-1.7){$\lambda$}
\rput(4.3,-1.7){$\hat{b}$}
\rput(5.5,-1.8){$\mu$}
\rput(6.5,-2.){$u_1$}
\rput(3.5,-2.3){$\hat{b}$}
\rput(4.3,-2.3){$\lambda$}
\rput(5.5,-2.3){$a$}
%
\rput(4.8,-1.4){$\Delta$}
\rput(5.2,-0.5){$\hat{\Delta}$}
%
\psline[linestyle=dotted]{->}(4.5,-1.3)(5.5,-.7)
\rput(5.5,-1.0){$\frac{\pi}{6}$}

\rput(5.0,-3.3){G (i)}
\end{pspicture*}

&

\begin{pspicture*}(2.,-4.)(7.5,1.5)
\psdots[dotstyle=*](4.,0.)
\psdots[dotstyle=*](6.,0.)
\psdots[dotstyle=*](4.,-2.)
\psdots[dotstyle=*](6.,-2.)
\psline(4.,0)(6.,0.)
\psline(4.,-2.)(6.,-2.)
\psline(4.,0)(4.,-2.)
\psline(6.,0)(6.,-2.)
\psline(4.,0)(6.,-2.)
\psline(3.3,0.)(4.,0.)
\psline(3.3,-2.)(4.,-2.)
\psline(4.,0.)(4.,0.7)
\psline(6.,0.)(6.,0.7)
\psline(4.,-2.)(4.,-2.7)
\psline(6.,-2.)(6.,-2.7)
\psline(4.,0.)(3.3,0.7)
%
\rput(3.5,0.2){$\mu$}
\rput(3.8,0.5){$a$}
\rput(4.3,0.2){$\hat{b}$}
\rput(5.6,0.2){$\mu$}
\rput(6.6,0.0){$\geq 6$}
\rput(3.5,-0.2){$b$}
\rput(4.2,-0.5){$\hat{a}$}
\rput(4.6,-0.2){$\lambda$}
\rput(5.6,-0.2){$a$}
\rput(3.5,-1.7){$\lambda$}
\rput(4.3,-1.7){$\hat{b}$}
\rput(5.5,-1.8){$\mu$}
\rput(5.8,-1.5){$b$}
\rput(6.6,-2.){$\geq 8$}
\rput(3.5,-2.3){$\hat{b}$}
\rput(4.3,-2.3){$\lambda$}
\rput(5.5,-2.3){$a$}
%
\rput(4.8,-1.4){$\Delta$}
\rput(5.2,-0.5){$\hat{\Delta}$}
%
\psline{->}(4.5,-1.3)(5.5,-.7)
\rput(5.5,-1.0){$\frac{\pi}{12}$}

\rput(5.0,-3.3){G (ii)}
\end{pspicture*}

&

\begin{pspicture*}(2.,-4.)(7.5,1.5)
\psdots[dotstyle=*](4.,0.)
\psdots[dotstyle=*](6.,0.)
\psdots[dotstyle=*](4.,-2.)
\psdots[dotstyle=*](6.,-2.)
\psline(4.,0)(6.,0.)
\psline(4.,-2.)(6.,-2.)
\psline(4.,0)(4.,-2.)
\psline(6.,0)(6.,-2.)
\psline(4.,0)(6.,-2.)
\psline(3.3,0.)(4.,0.)
\psline(3.3,-2.)(4.,-2.)
\psline(4.,0.)(4.,0.7)
\psline(6.,0.)(6.,0.7)
\psline(4.,-2.)(4.,-2.7)
\psline(6.,-2.)(6.,-2.7)
\psline(4.,0.)(3.3,0.7)
%
\rput(3.5,0.2){$\mu$}
\rput(3.8,0.5){$a$}
\rput(4.3,0.2){$\mu$}
\rput(5.6,0.2){$\hat{a}$}
\rput(6.6,0.0){$\geq 6$}
\rput(3.5,-0.2){$b$}
\rput(4.2,-0.5){$\hat{a}$}
\rput(4.6,-0.2){$b$}
\rput(5.6,-0.2){$\lambda$}
\rput(3.5,-1.7){$\lambda$}
\rput(4.3,-1.7){$\hat{b}$}
\rput(5.5,-1.8){$\mu$}
\rput(5.8,-1.5){$a$}
\rput(6.6,-2.){$\geq 8$}
\rput(3.5,-2.3){$\hat{b}$}
\rput(4.3,-2.3){$\lambda$}
\rput(5.5,-2.3){$a$}
%
\rput(4.8,-1.4){$\Delta$}
\rput(5.2,-0.5){$\hat{\Delta}$}
%
\psline{->}(4.5,-1.3)(5.5,-.7)
\rput(5.5,-1.0){$\frac{\pi}{12}$}

\rput(5.0,-3.3){G (iii)}
\end{pspicture*}
\\
\end{tabular}
\end{center}

\newpage
\subsection*{Figure {{H}}}
\psset{arrowscale=0.7,linewidth=1.0pt,dotsize=3pt}
\psset{xunit=1.0cm,yunit=1.0cm,}

\begin{center}
\begin{tabular}{c@{}c}

\begin{pspicture*}(1.,2.5)(9.,7.)
\psdots[dotstyle=*](4.,4.)
\rput(3.7,4.3){$\lambda$}
\rput(4.2,4.3){$\hat{b}$}
\rput(4.2,3.8){$\lambda$}
\rput(3.7,3.7){$\hat{b}$}

\psdots[dotstyle=*](4.,6.)
\rput(3.5,6.2){$\mu$}
\rput(3.8,6.5){$a$}
\rput(4.2,6.3){$\hat{b}$}
\rput(4.5,5.8){$\lambda$}
\rput(4.2,5.6){$\hat{a}$}
\rput(3.7,5.8){$b$}

\psdots[dotstyle=*](6.,4.)
\rput(5.5,4.2){$\mu$}
\rput(5.8,4.5){$b$}
\rput(6.2,4.2){$\mu$}
\rput(6.5,3.8){$d$}
\rput(6.2,3.5){$\hat{b}$}
\rput(5.7,3.8){$a$}

\psdots[dotstyle=*](6.,6.)
\rput(5.7,6.2){$\mu$}
\rput(6.2,5.5){$\hat{b}$}
\rput(5.7,5.8){$a$}

\psdots[dotstyle=*](8.,4.)
\rput(7.6,4.2){$\hat{a}$}
\psline(6.5,3.5)(3.5,6.5)
\psline(6.,6.)(8.,4.)

\psline(3.5,6.)(6.,6.)
\psline(3.5,4.)(8.,4.)
\psline(4.,3.5)(4.,6.5)
\psline(6.,3.5)(6.,6.5)
%
\rput(4.7,4.7){$\Delta$}
\rput(6.9,4.7){$\hat{\Delta}$}
%
\psline{->}(4.7,5.)(6.4,5.)
\rput(6.6,5.){$\frac{\pi}{6}$}

\rput(6.0,2.8){H (i)}
\end{pspicture*}

&

\begin{pspicture*}(1.,2.5)(9.,7.)
\psdots[dotstyle=*](4.,4.)
\rput(3.7,4.3){$\lambda$}
\rput(4.2,4.3){$\hat{b}$}
\rput(4.2,3.8){$\lambda$}
\rput(3.7,3.7){$\hat{b}$}

\psdots[dotstyle=*](4.,6.)
\rput(3.5,6.2){$\mu$}
\rput(3.8,6.5){$b$}
\rput(4.2,6.3){$\hat{d}$}
\rput(4.5,5.8){$\lambda$}
\rput(4.2,5.6){$\hat{a}$}
\rput(3.7,5.8){$b$}

\psdots[dotstyle=*](6.,4.)
\rput(5.5,4.2){$\mu$}
\rput(5.8,4.5){$b$}
\rput(6.2,4.2){$\hat{a}$}
\rput(5.7,3.8){$a$}

\psdots[dotstyle=*](6.,6.)
\rput(6.2,5.7){$\mu$}
\rput(5.7,5.8){$a$}

\psline(6.5,3.5)(3.5,6.5)
\psline(3.5,6.)(6.5,6.)
\psline(3.5,4.)(6.5,4.)
\psline(4.,3.5)(4.,6.5)
\psline(6.,3.5)(6.,6.)
%
\rput(4.7,4.7){$\Delta$}
\rput(5.3,5.3){$\hat{\Delta}$}
%

\psline{->}(5.0,4.6)(5.4,5.0)
\rput(5.6,4.9){$\frac{\pi}{6}$}

\rput(6.0,2.8){H (ii)}
\end{pspicture*}
\\
\end{tabular}
\end{center}

\pagebreak

\begin{center}
\begin{tabular}{c@{}c}

\begin{pspicture*}(1.,0.5)(9.,11.)
\psdots[dotstyle=*](4.,4.)
\rput(3.7,4.){$v_1$}
\rput(4.2,4.3){$\hat{b}$}
\rput(4.4,3.8){$a$}
\rput(4.2,3.5){$\mu$}

\psdots[dotstyle=*](4.,6.)
\rput(4.2,6.2){$\hat{a}$}
\rput(4.5,5.8){$b$}
\rput(4.2,5.6){$\hat{a}$}
\rput(3.7,6.){$v_2$}

\psdots[dotstyle=*](4.,8.)
\rput(3.7,8.){$v_3$}
\rput(4.2,7.5){$\mu$}
\psdots[dotstyle=*](6.,2.)
\rput(5.8,1.9){$\hat{b}$}
\rput(5.8,2.5){$\lambda$}
\rput(6.2,2.2){$\hat{b}$}
\rput(6.2,1.8){$\lambda$}

\psdots[dotstyle=*](6.,4.)
\rput(5.5,4.2){$\mu$}
\rput(5.8,4.5){$a$}
\rput(6.2,4.3){$\hat{b}$}
\rput(6.5,3.8){$\lambda$}
\rput(6.2,3.6){$\hat{a}$}
\rput(5.7,3.8){$b$}

\psdots[dotstyle=*](6.,6.)
\rput(5.5,6.2){$\hat{b}$}
\rput(6.2,6.1){$e$}
\rput(6.2,5.6){$\hat{a}$}
\rput(5.7,5.8){$\lambda$}

\psdots[dotstyle=*](8.,2.)
\rput(7.5,2.2){$\mu$}
\rput(7.8,2.45){$b$}
\rput(8.2,2.2){$\hat{a}$}
\rput(7.7,1.8){$a$}

\psdots[dotstyle=*](8.,4.)
\rput(7.5,4.2){$\mu$}
\rput(8.2,4.1){$d$}
\rput(8.2,3.5){$\mu$}
\rput(7.7,3.8){$a$}

\psline(6.,2.)(4.,4.)
\psline(8.5,1.5)(4.,6.)
\psline(8.5,3.5)(4.,8.)
\psline(4.,6.)(6.,6.)
\psline(4.,4.)(8.,4.)
\psline(6.,2.)(8.5,2.)
\psline(4.,3.5)(4.,8.)
\psline(6.,1.5)(6.,6.)
\psline(8.,1.5)(8.,4.)
%
\rput(6.7,2.7){$\Delta$}
\rput(5.3,5.7){$\hat{\Delta}_2$}
\rput(4.3,5.2){$\hat{\Delta}_1$}
%
\psbezier{->}
(6.7,2.9)(6.7,4.5)
(7.,5.)(5.6,5.0)
\rput(5.5,5.2){$\frac{\pi}{6}$}

\psline[linestyle=dotted]{->}(5.6,5.0)(4.8,4.5)
\rput(4.7,4.6){$\frac{\pi}{6}$}

\rput(6.0,1.0){H (ix)}
\end{pspicture*}

&

\begin{pspicture*}(1.,0.5)(9.,11.)
\psdots[dotstyle=*](4.,4.)
\rput(3.6,4.){$\geq 6$}
\rput(4.2,4.3){$\hat{b}$}
\rput(4.4,3.8){$a$}
\rput(4.2,3.5){$\mu$}

\psdots[dotstyle=*](4.,6.)
\rput(4.2,6.2){$\hat{a}$}
\rput(4.5,5.8){$b$}
\rput(4.2,5.6){$\hat{a}$}
\rput(3.6,6.){$\geq 8$}

\psdots[dotstyle=*](4.,8.)
\rput(3.5,8.){$\encircle{\geq 6}$}
\rput(4.2,7.5){$\mu$}

\psdots[dotstyle=*](6.,2.)
\rput(5.8,1.9){$\hat{b}$}
\rput(5.8,2.5){$\lambda$}
\rput(6.2,2.2){$\hat{b}$}
\rput(6.2,1.8){$\lambda$}

\psdots[dotstyle=*](6.,4.)
\rput(5.5,4.2){$\mu$}
\rput(5.8,4.5){$a$}
\rput(6.2,4.25){$\hat{b}$}
\rput(6.5,3.8){$\lambda$}
\rput(6.2,3.6){$\hat{a}$}
\rput(5.7,3.8){$b$}

\psdots[dotstyle=*](6.,6.)
\rput(5.5,6.2){$\hat{b}$}
\rput(6.2,6.1){$e$}
\rput(6.2,5.6){$\hat{a}$}
\rput(5.7,5.8){$\lambda$}

\psdots[dotstyle=*](8.,2.)
\rput(7.5,2.2){$\mu$}
\rput(7.8,2.45){$b$}
\rput(8.2,2.2){$\hat{a}$}
\rput(7.7,1.8){$a$}

\psdots[dotstyle=*](8.,4.)
\rput(7.5,4.2){$\mu$}
\rput(8.2,4.1){$d$}
\rput(8.2,3.5){$\mu$}
\rput(7.7,3.8){$a$}
%
\psline(6.,2.)(4.,4.)
\psline(8.5,1.5)(4.,6.)
\psline(8.5,3.5)(4.,8.)
\psline(4.,6.)(6.,6.)
\psline(4.,4.)(8.,4.)
\psline(6.,2.)(8.5,2.)
\psline(4.,3.5)(4.,8.)
\psline(6.,1.5)(6.,6.)
\psline(8.,1.5)(8.,4.)
%
\rput(6.7,2.7){$\Delta$}
\rput(4.5,7.){$\hat{\Delta}$}
\rput(5.3,5.7){$\Delta_1$}
\rput(7.,4.7){$\hat{\Delta}_1$}
%
\psbezier{->}
(6.7,2.9)(6.7,7.)
(5.,4)(4.7,6.3)
\rput(4.9,6.3){$\frac{\pi}{6}$}

\psline{->}(5.7,4.9)(6.3,4.9)
\rput(6.2,4.6){$\frac{\pi}{12}$}

\rput(6.0,1.0){H (x)}
\end{pspicture*}
\\
\end{tabular}%
\end{center}

\newpage
\subsection*{Figure {{I}}}
\psset{arrowscale=0.7,linewidth=1.0pt,dotsize=3pt}
\psset{xunit=1.0cm,yunit=1.0cm,}

\begin{center}
%
\end{center}

\newpage
\subsection*{Figure {{J}}}
\psset{arrowscale=0.7,linewidth=1.0pt,dotsize=3pt}
\psset{xunit=1.0cm,yunit=1.0cm,}
\begin{center}
%
\end{center}

\newpage
\subsection*{Figure {{K}}}
\psset{arrowscale=0.7,linewidth=1.0pt,dotsize=3pt}
\psset{xunit=.8cm,yunit=1.0cm,}
\begin{center}

\begin{tabular}{@{}c@{}c@{}c}

\begin{pspicture*}(0.5,-5.)(7.0,1.5)
\psdots[dotstyle=*](4.,0.)
\rput(4.,0.3){$u^*$}
\psdots[dotstyle=*](3.,-2.)
\rput(3.,-2.3){$u_1$}
\psdots[dotstyle=*](5.,-2.)
\rput(5.,-2.3){$u_2$}
\psline(2.,0.)(6.,0.) 
\psline(3.,-2.)(5.,-2.)

\psline(3.,-2.)(4.,0.)
\psline(4.,0.)(5.,-2.)
\rput(4.0,-1.0){$\hat{\Delta}$}
%

\rput(4.0,-3.5){K (i)}
\end{pspicture*}
&
\begin{pspicture*}(0.5,-5.)(7.0,1.5)
\psdots[dotstyle=*](4.,0.)
\rput(3.5,0.3){$u^*$}
\psdots[dotstyle=*](2.,-2.)
\rput(2.,-2.3){$u_3$}
\psdots[dotstyle=*](4.,-2.)
\rput(4.,-2.3){$u_1$}
\psdots[dotstyle=*](6.,-2.)
\rput(6.,-2.3){$u_2$}

\psline(2.,0.)(6.,0.) 
\psline(1.5,-2.)(7.,-2.)

\psline(2.,-2.)(4.,0.)
\psline(4.,-2.)(4.,0.)
\psline(6.,-2.)(4.,0.)
%
\psline(4.,-2.)(5.,-2.3)
\rput(3.5,-1.0){$\hat{\Delta}_1$}
\rput(4.5,-1.0){$\hat{\Delta}$}
\rput(5.,0.5){$\Delta^*$}
\rput(2.5,-3.){$\Delta^*$}
%
\psline{->}(3.0,-1.8)(3.0,-2.2)
\rput(3.0,-2.5){$\frac{\pi}{2}$}
\psline{->}(4.2,-1.5)(3.8,-1.5)
\rput(3.6,-1.5){$\frac{\pi}{6}$}
\psline{->}(3.8,-0.5)(4.2,0.6)
\psline{->}(4.2,-0.5)(3.8,0.6)

\rput(4.0,-3.5){K (ii)}
\end{pspicture*}
 &

\begin{pspicture*}(0.5,-5.)(7.0,1.5)
\psdots[dotstyle=*](4.,0.)
\rput(3.5,0.3){$u^*$}
\psdots[dotstyle=*](2.,-2.)
\rput(2.,-2.3){$u_3$}
\psdots[dotstyle=*](4.,-2.)
\rput(4.,-2.3){$u_1$}
\psdots[dotstyle=*](6.,-2.)
\rput(6.,-2.3){$u_2$}

\psline(2.,0.)(6.,0.) 
\psline(1.5,-2.)(7.,-2.)

\psline(2.,-2.)(4.,0.)
\psline(4.,-2.)(4.,0.)
\psline(6.,-2.)(4.,0.)
%
\psline(4.,-2.)(5.,-2.3)

\rput(3.5,-1.0){$\hat{\Delta}_1$}
\rput(4.5,-1.0){$\hat{\Delta}$}
\rput(5.,0.5){$\Delta^*$}
\rput(2.5,-3.){$\Delta^*$}
%
\psline{->}(4.2,-1.5)(3.8,-1.5)
\rput(3.6,-1.5){$\frac{\pi}{6}$}
\psline{->}(2.6,-1.8)(1.4,-2.2)

\psline{->}(3.8,-0.5)(4.2,0.6)
\psline{->}(4.2,-0.5)(3.8,0.6)

\rput(4.0,-3.5){K (iii)}
\end{pspicture*}
\\
\end{tabular}
\end{center}


\begin{thebibliography}{10}

\bibitem{BaumslagGerstenShapiroShort}
G.~Baumslag, S.~M.~Gersten, M.~Shapiro, and H.~Short.
\newblock Automatic groups and amalgams.
\newblock {\em J. Pure Appl. Algebra}, 76(3):229--316, 1991.
\newblock {\url{https://doi.org/10.1016/0022-4049(91)90139-S}}

\bibitem{BogleyPride}
W.~A.~Bogley and S.~J.~Pride.
\newblock Aspherical relative presentations.
\newblock {\em Proc. Edinb. Math. Soc. (2)}, 35(1):1--39, 1992.
\newblock {\url{https://doi.org/10.1017/S0013091500005290}}

\bibitem{BogleyShift}
W.~A.~Bogley.
\newblock On shift dynamics for cyclically presented groups.
\newblock {\em J. Algebra}, 418:154--173, 2014.
\newblock {\url{https://doi.org/10.1016/j.jalgebra.2014.07.009}}

\bibitem{BrownHuebschmann}
R.~Brown and J.~Huebschmann,
\newblock Identities among relations.
\newblock {\em Low-dimensional topology (Bangor, 1979), }
\newblock Lond. Math. Soc. Lecture Note Ser. 48, pages 153--202, 1982.\newline 
\newblock {\url{https://doi.org/10.1017/CBO9780511758935.010}}

\bibitem{BW2}
W.~A.~Bogley and G.~Williams.
\newblock Coherence, subgroup separability, and metacyclic structures for a class of cyclically presented groups.
\newblock {\em J. Algebra}, 480:266--297, 2017.\newline
\newblock {\url{https://doi.org/10.1016/j.jalgebra.2017.02.002}}

\bibitem{CRS05}
A.~Cavicchioli, D.~Repov\v{s}, and F.~Spaggiari.
\newblock Families of group presentations related to topology.
\newblock {\em J. Algebra}, 286(1):41--56, 2005.
\newblock {\url{https://doi.org/10.1016/j.jalgebra.2005.01.008}}

\bibitem{Chalk}
C.~P.~Chalk.
\newblock Fibonacci groups {$F(2,n)$} are hyperbolic for {$n$} odd and {$n\geq 11$}.
\newblock {\em J. Group Theory}, 24(2):359--384, 2021.
\newblock {\url{https://doi.org/10.1515/jgth-2020-0068}}

\bibitem{ChalkEdjvet}
C.~P.~Chalk and M.~Edjvet.
\newblock Curvature distribution and hyperbolicity.
\newblock {\em J. Group Theory}, 26(4):645--663, 2023.
\newblock {\url{https://doi.org/10.1515/jgth-2022-0106}}

\bibitem{ChalkEdjvetJuhasz}
C.~P.~Chalk, M.~Edjvet, and A.~Juh\'{a}sz.
\newblock Curvature distribution, relative presentations and hyperbolicity with an application to {F}ibonacci groups.
\newblock {\em J. Algebra}, 644:796--834, 2024.\newline
\newblock {\url{https://doi.org/10.1016/j.jalgebra.2024.01.023}}

\bibitem{ChinyereWilliamsT5}
I.~Chinyere and G.~Williams.
\newblock Hyperbolic groups of {F}ibonacci type and {T}(5) cyclically presented groups.
\newblock {\em J. Algebra}, 580:104--126, 2021.\newline
\newblock {\url{https://doi.org/10.1016/j.jalgebra.2021.04.003}}

\bibitem{ChinyereWilliamsT6}
I.~Chinyere and G.~Williams.
\newblock Hyperbolicity of {$T(6)$} cyclically presented groups.
\newblock {\em Groups Geom. Dyn.}, 16(1):341--361, 2022.
\newblock {\url{https://doi.org/10.4171/GGD/651}}

\bibitem{EdjvetIrreducible}
M.~Edjvet.
\newblock On irreducible cyclic presentations.
\newblock {\em J. Group Theory}, 6(2):261--270, 2003.\newline 
\newblock {\url{https://doi.org/10.1515/jgth.2003.019}}

\bibitem{EdjvetHowie91}
M.~Edjvet and J.~Howie.
\newblock The solution of length four equations over groups.
\newblock {\em Trans. Amer. Math. Soc.}, 326(1):345--369, 1991.
\newblock {\url{https://doi.org/10.1090/S0002-9947-1991-1002920-5}}

\bibitem{EdjvetWilliams}
M.~Edjvet and G.~Williams.
\newblock The cyclically presented groups with relators {$x_ix_{i+k}x_{i+l}$}.
\newblock {\em Groups Geom. Dyn.}, 4(4):759--775, 2010.
\newblock {\url{https://doi.org/10.4171/GGD/104}}

\bibitem{Howie83}
J.~Howie.
\newblock The solution of length three equations over groups.
\newblock {\em Proc. Edinb. Math. Soc. (2)}, 26(1):89--96, 1983.
\newblock {\url{https://doi.org/10.1017/S0013091500028108}}

\bibitem{HowieWilliamsTadpole}
J.~Howie and G.~Williams.
\newblock Tadpole labelled oriented graph groups and cyclically presented groups.
\newblock {\em J. Algebra}, 371:521--535, 2012.
\newblock {\url{https://doi.org/10.1016/j.jalgebra.2012.09.001}}

\bibitem{HuebschmannPersonalCommunication}
J.~Huebschmann.
\newblock Personal Communication.
\newblock 4/11/2025.

\bibitem{Huebschmann23}
J.~Huebschmann.
\newblock Crossed modules.
\newblock {\em Notices Amer. Math. Soc.}, 70(11):1802--1813, 2023.\newline
\newblock {\url{https://doi.org/10.1090/noti2831}}

\bibitem{McCammondWise}
J.~P.~McCammond and D.~T.~Wise.
\newblock Fans and ladders in small cancellation theory.
\newblock {\em Proc. Lond. Math. Soc. (3)}, 84(3):599--644, 2002.
\newblock {\url{https://doi.org/10.1112/S0024611502013424}}

\bibitem{MohamedWilliams}
E.~Mohamed and G.~Williams.
\newblock An investigation into the cyclically presented groups with length three positive relators.
\newblock {\em Exp. Math.}, 31(2):537--551, 2022.\newline 
\newblock {\url{https://doi.org/10.1080/10586458.2019.1655817}}

\bibitem{NoferiniWilliams3}
V.~Noferini and G.~Williams.
\newblock Smith forms of matrices in {C}ompanion {R}ings, with group theoretic and topological applications.
\newblock {\em Linear Algebra Appl.}, 708:372--404, 2025.\newline 
\newblock {\url{https://doi.org/10.1016/j.laa.2024.12.003}}

\bibitem{Peiffer49}
R.~Peiffer,
\newblock \"{U}ber Identit\"{a}ten zwischen Relationen,
\newblock {\em Math. Ann.}, 121:67--99, 1949. \newline
\newblock {\url{https://doi.org/10.1007/BF01329617}}
\end{thebibliography}
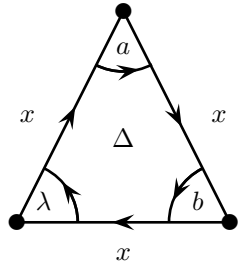
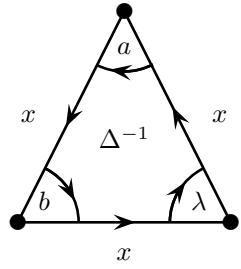
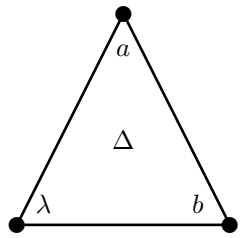
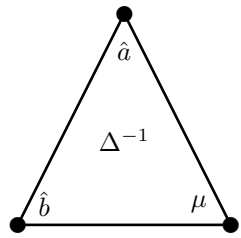
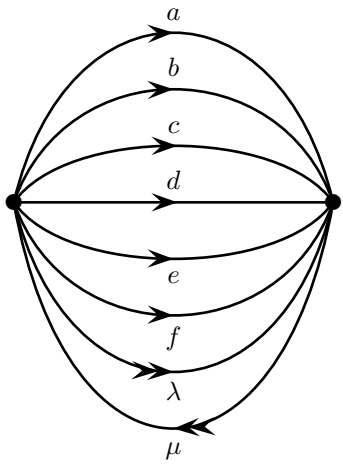
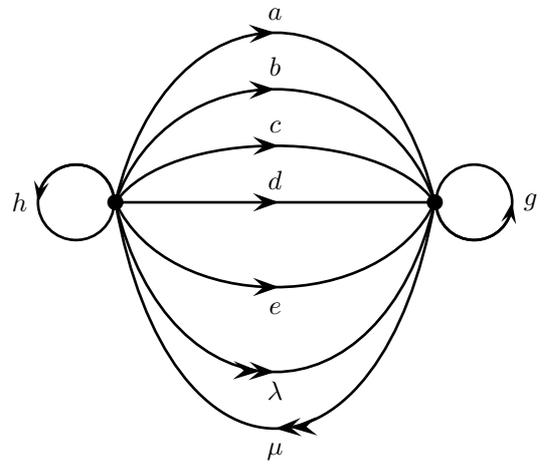
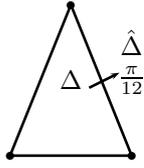
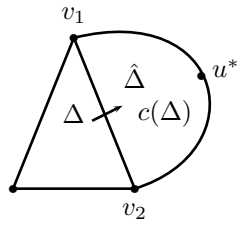
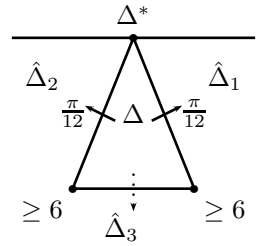
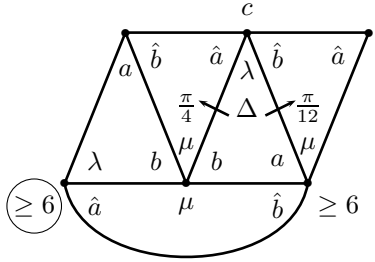
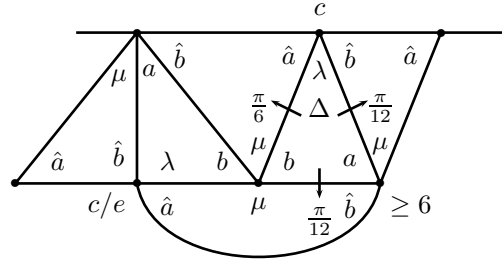
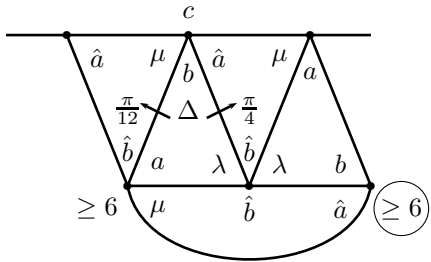
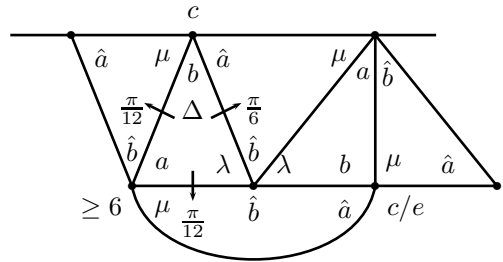
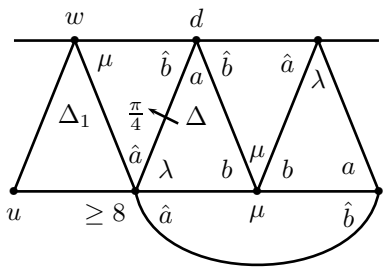
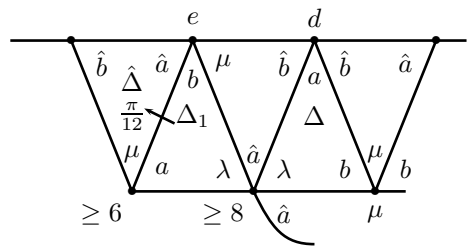
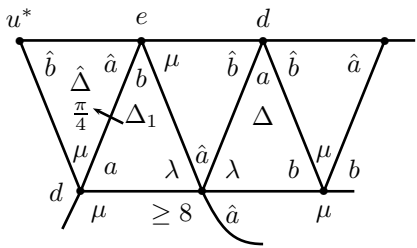
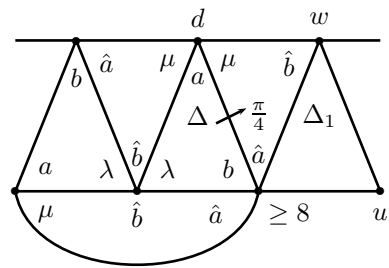
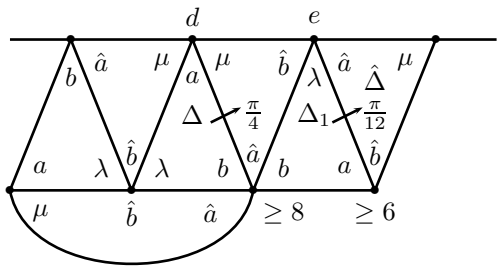
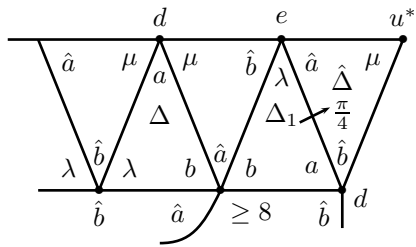
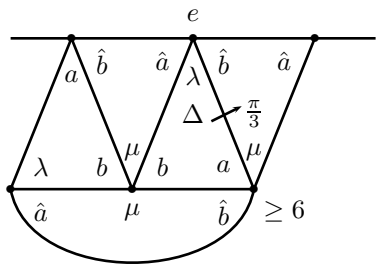
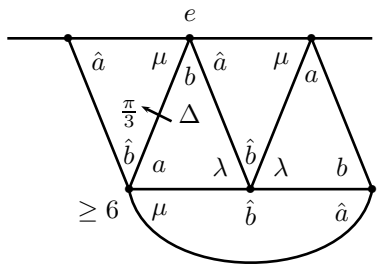
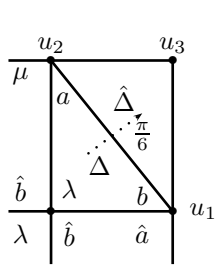
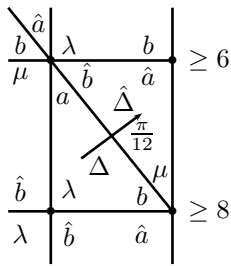
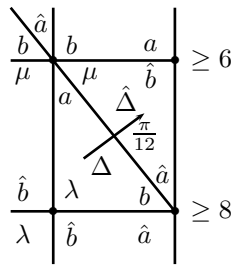
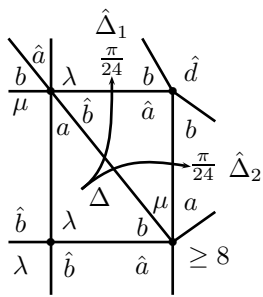
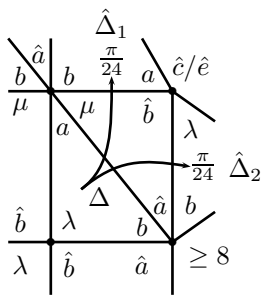
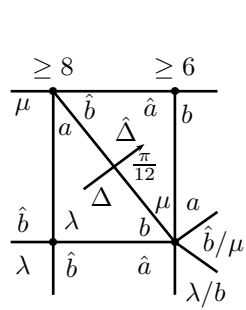
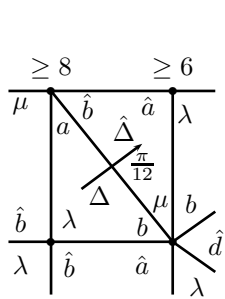
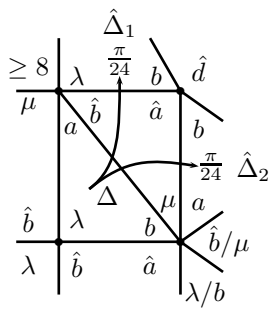
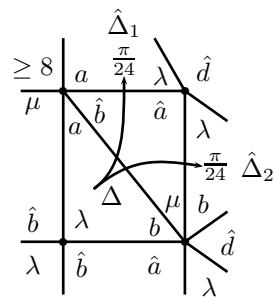
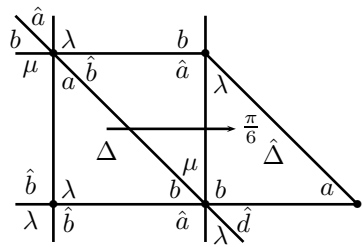
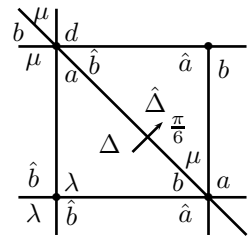
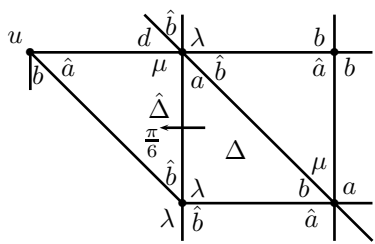
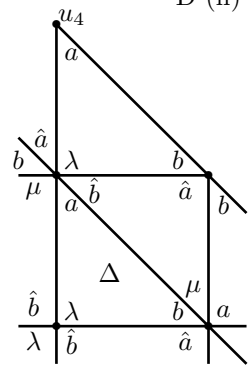
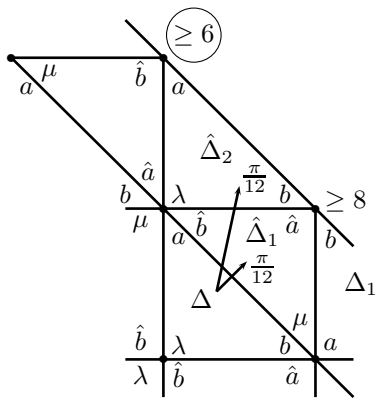
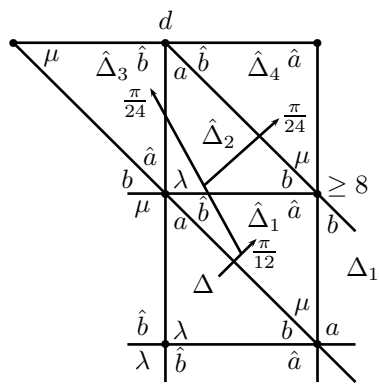
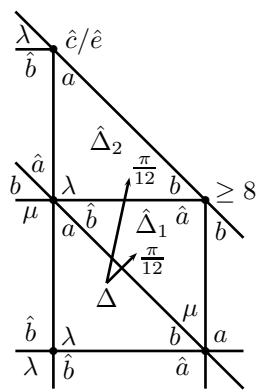
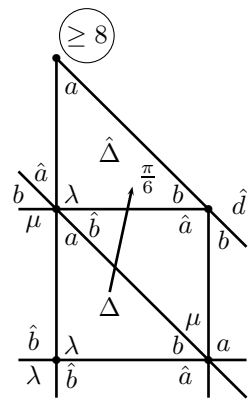
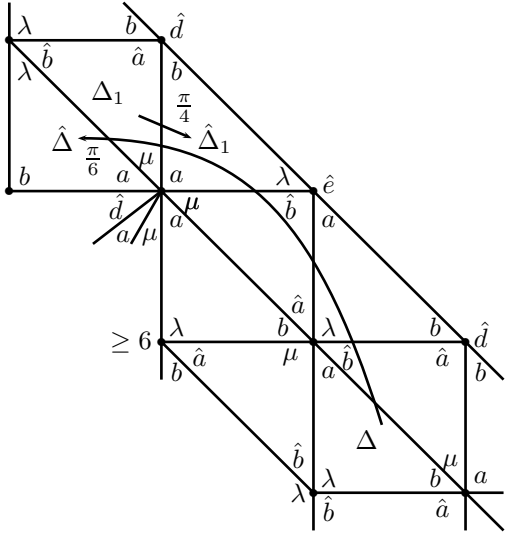
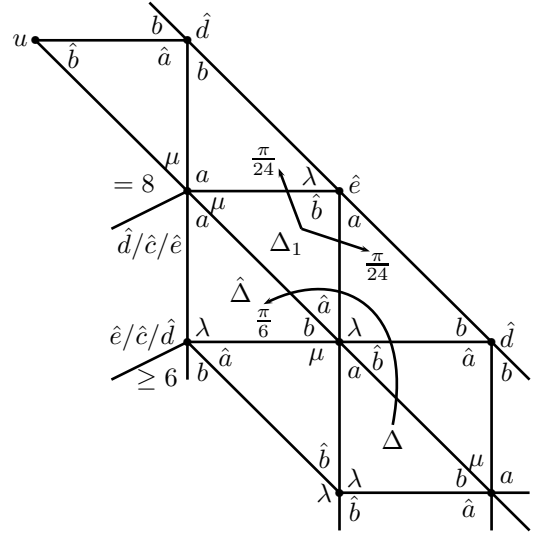
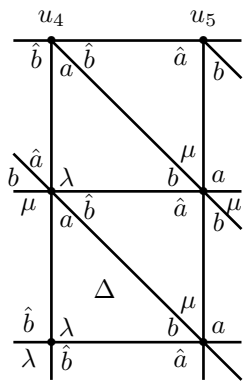
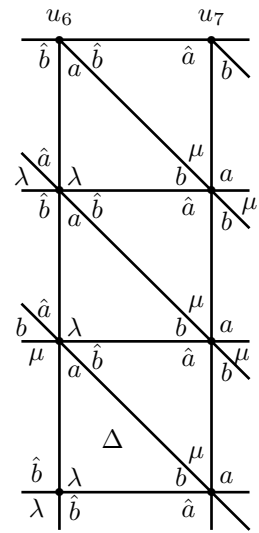
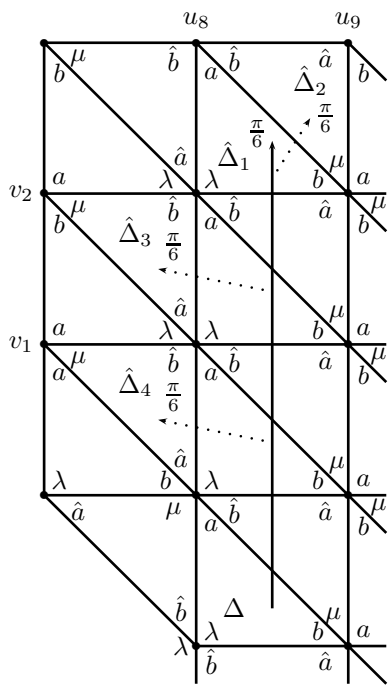
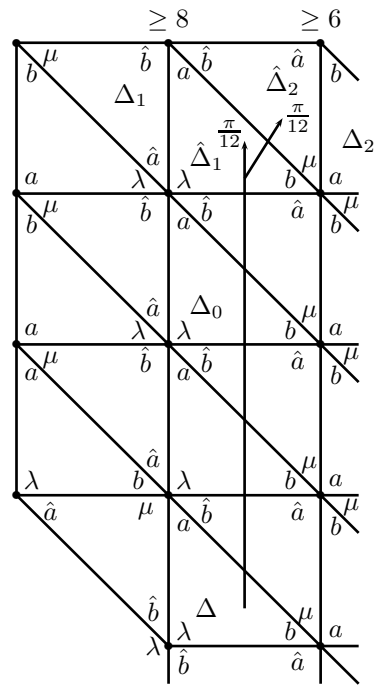
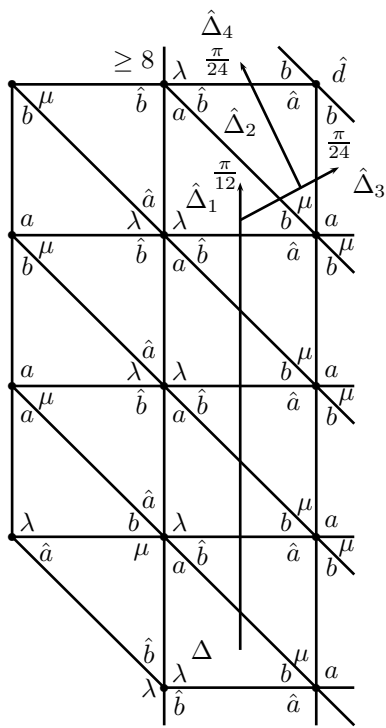
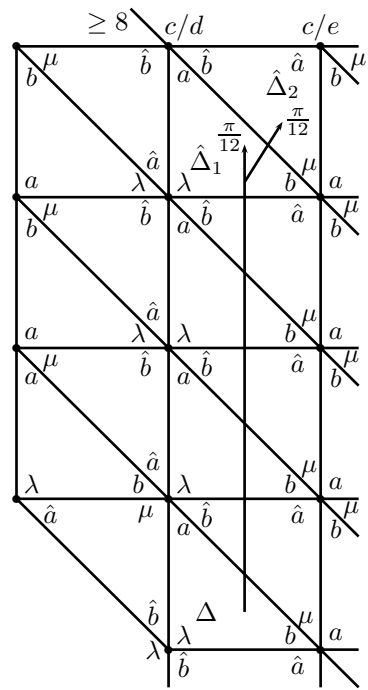
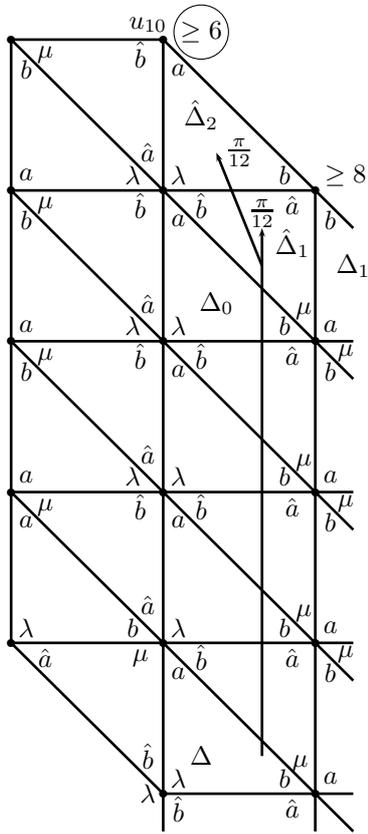
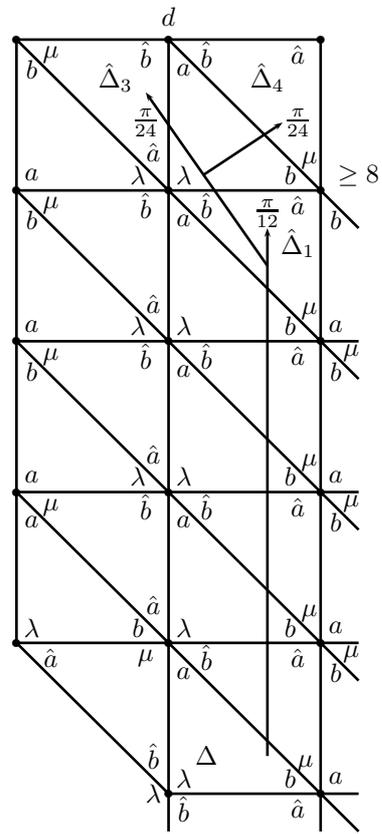
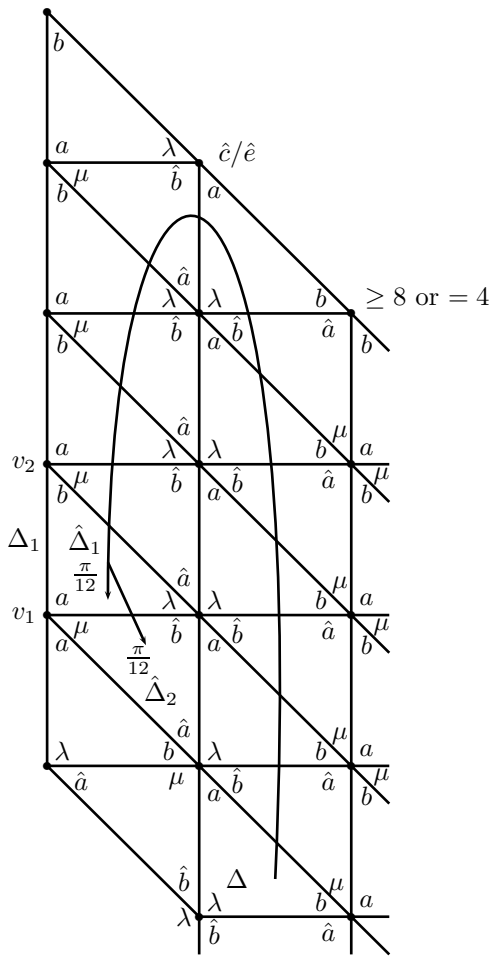
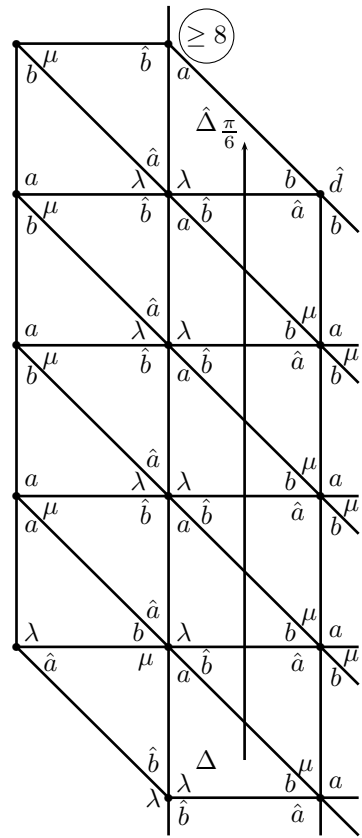
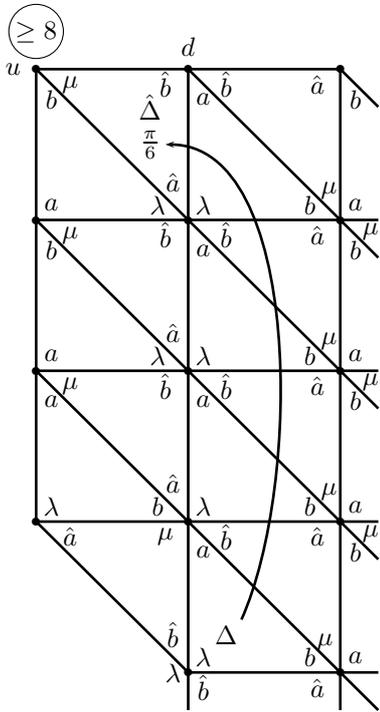
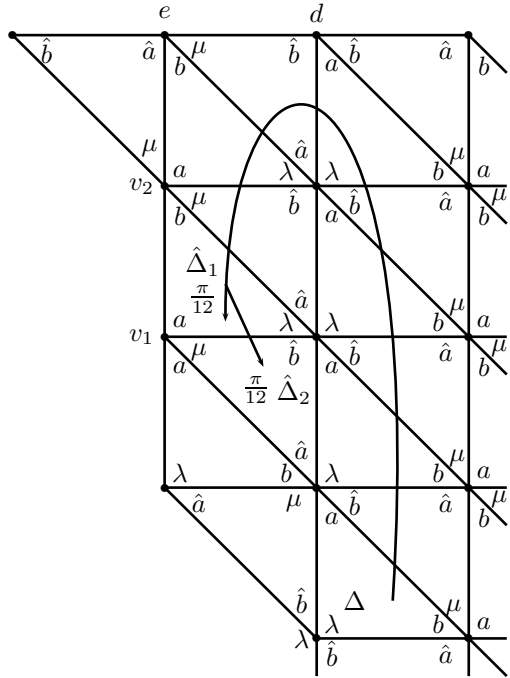
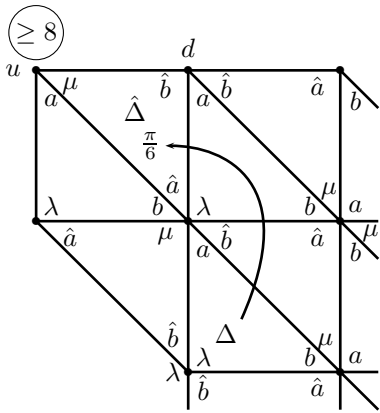
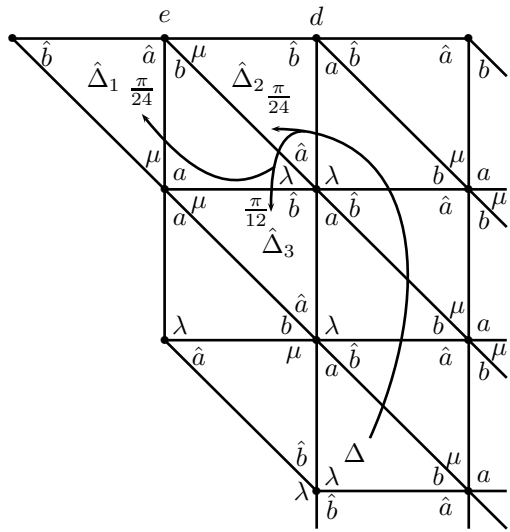
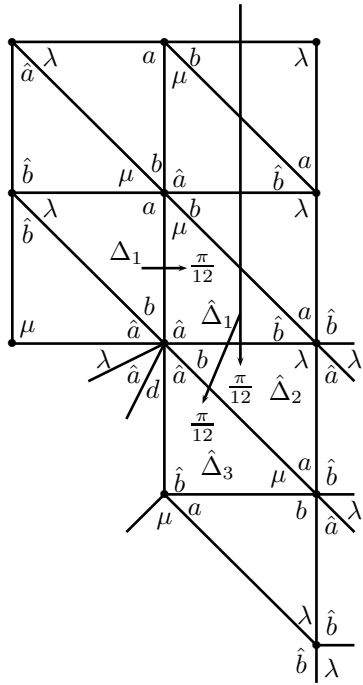
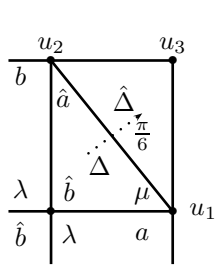
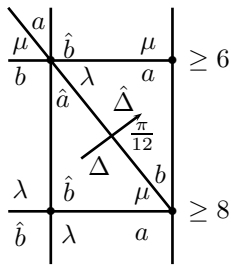
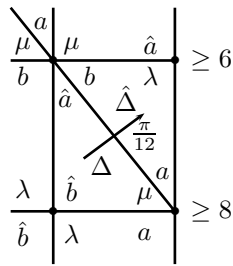
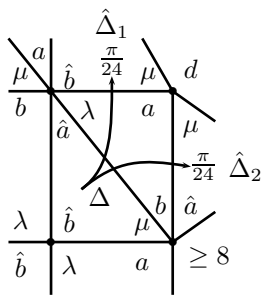
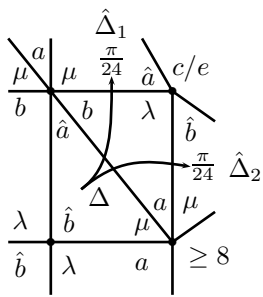
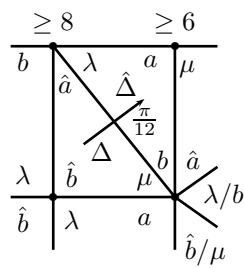
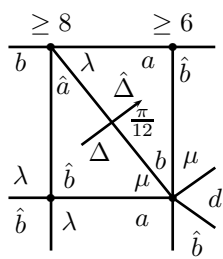
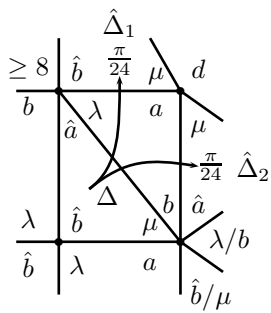
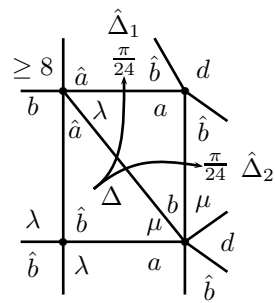
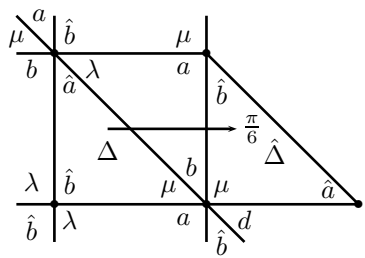
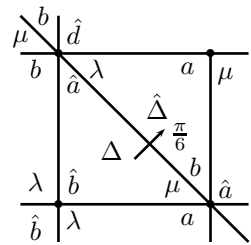
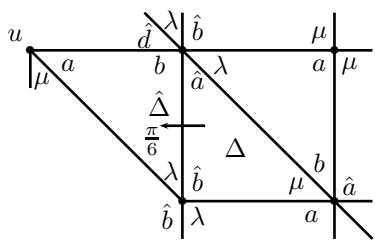
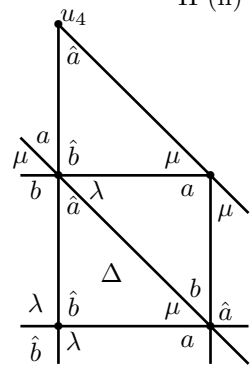
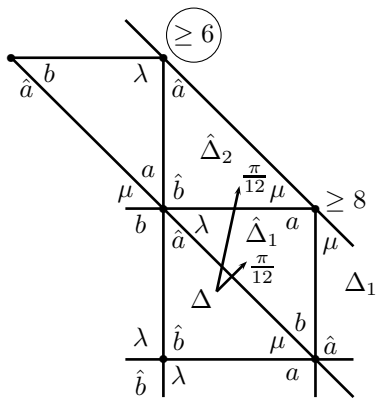
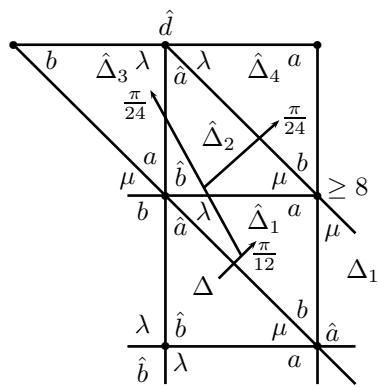
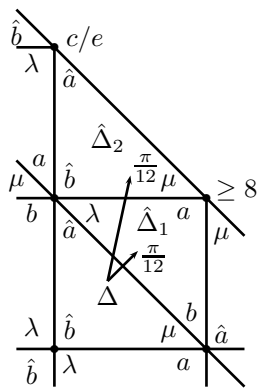
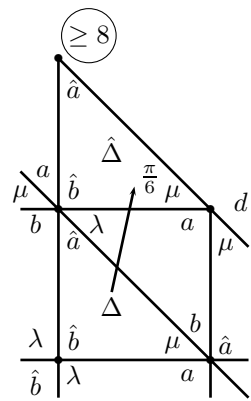
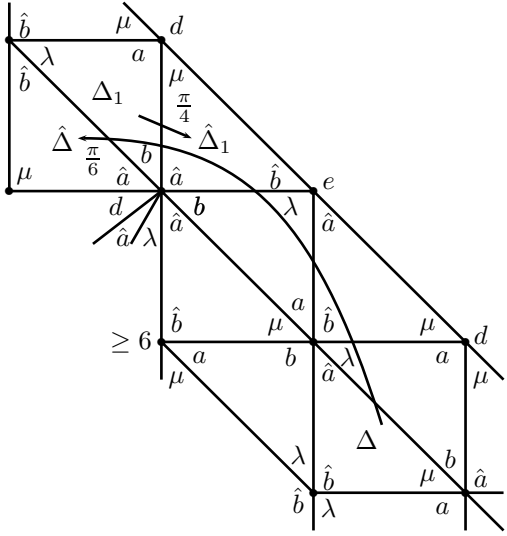
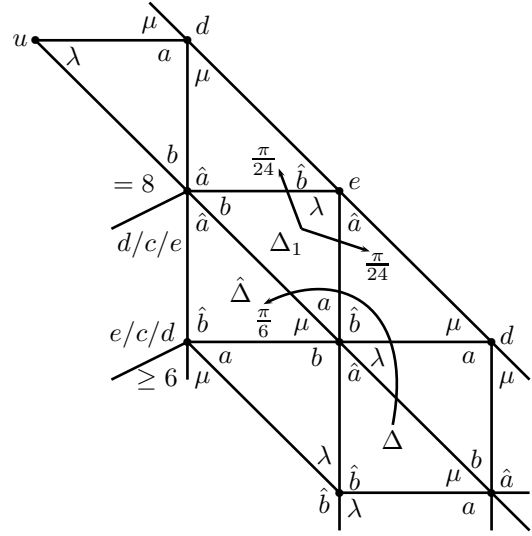
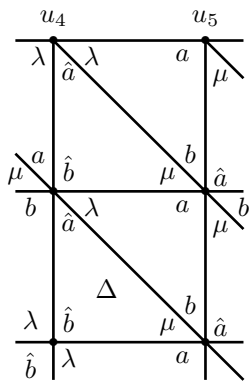
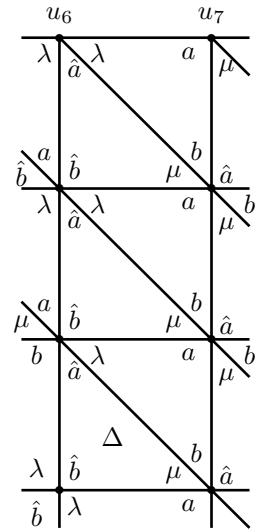
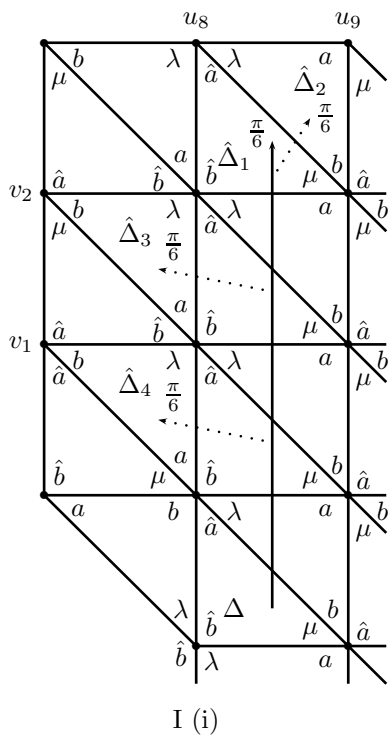
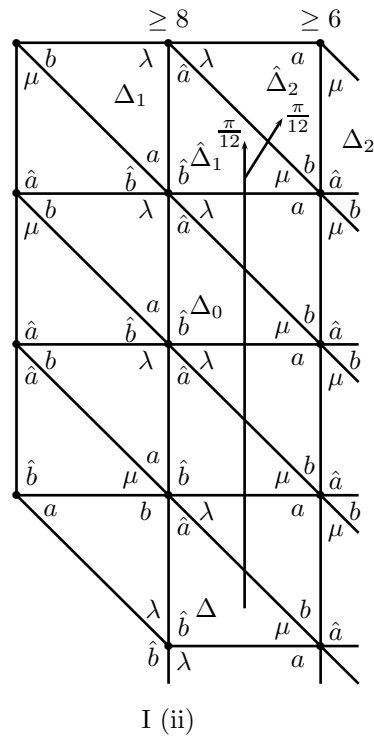
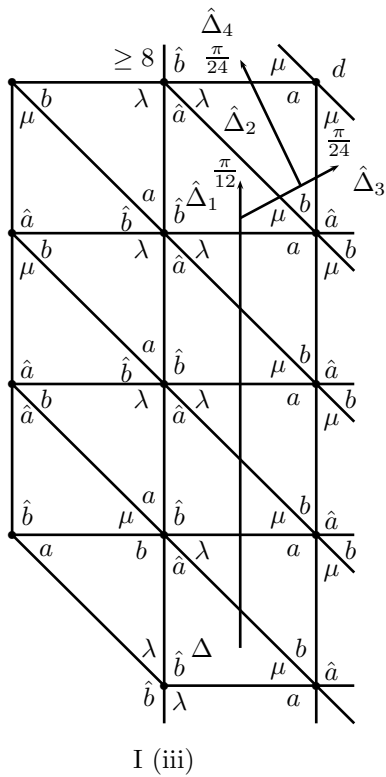
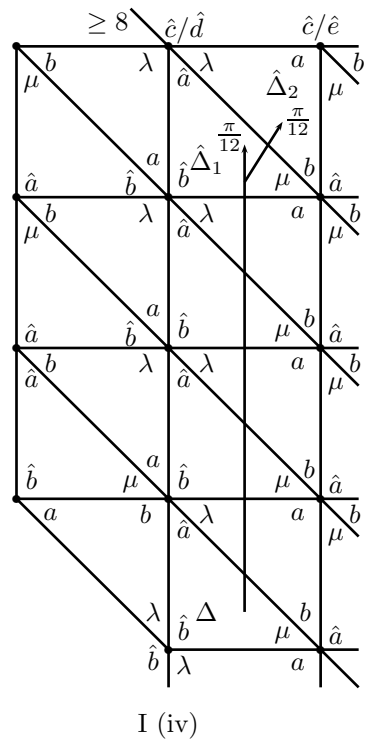
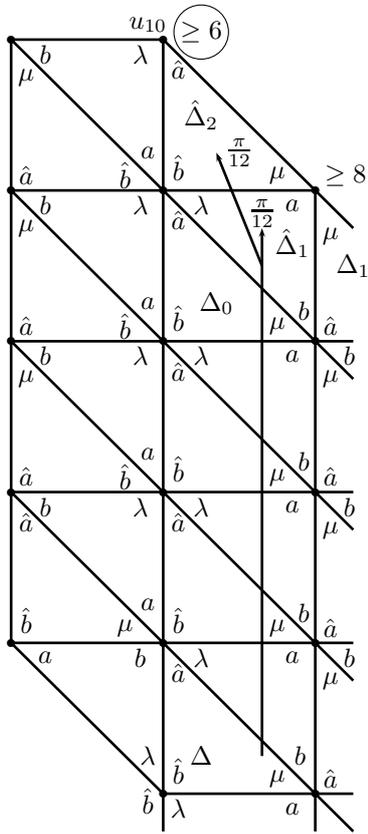
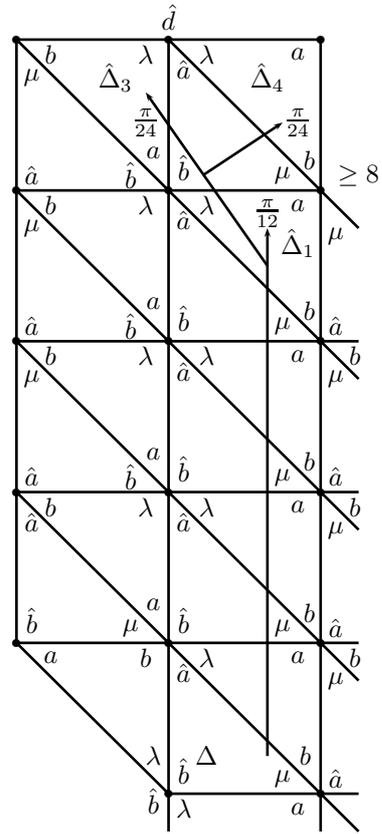
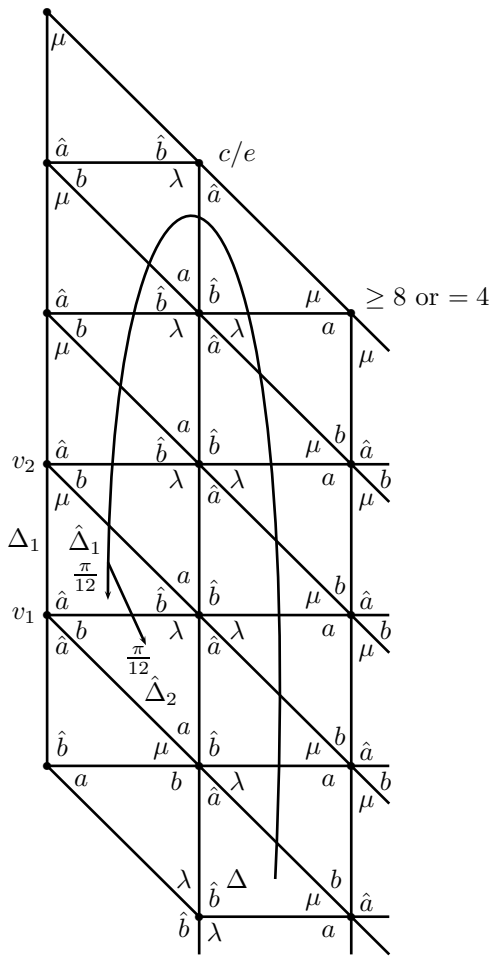
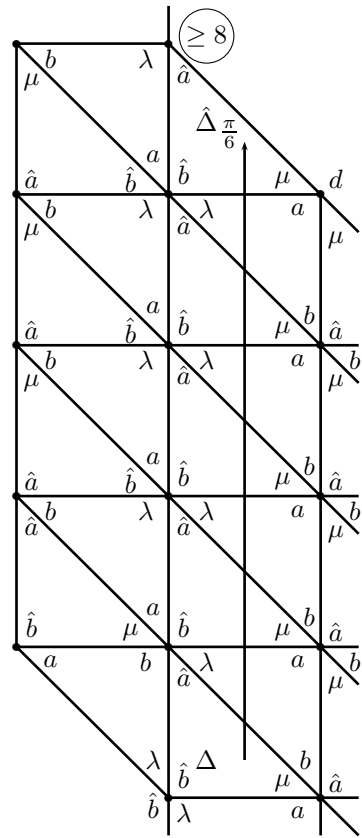
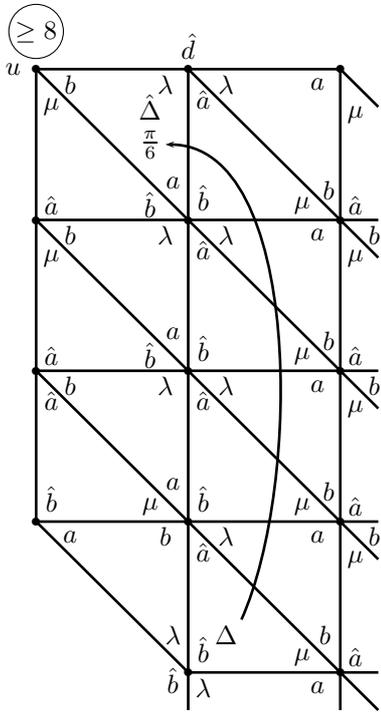
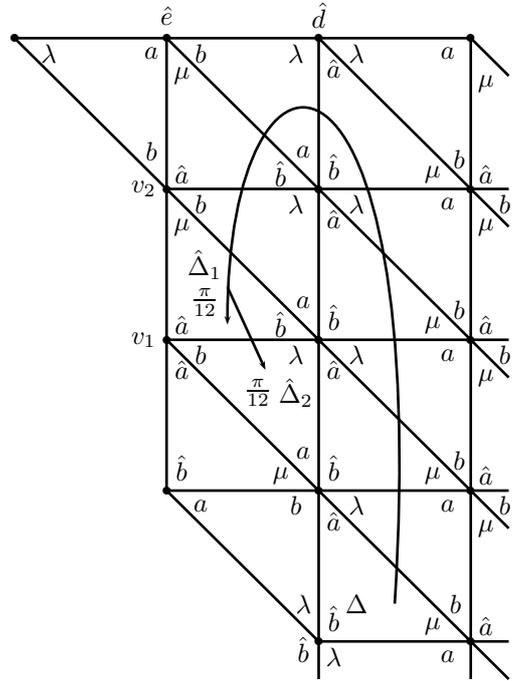
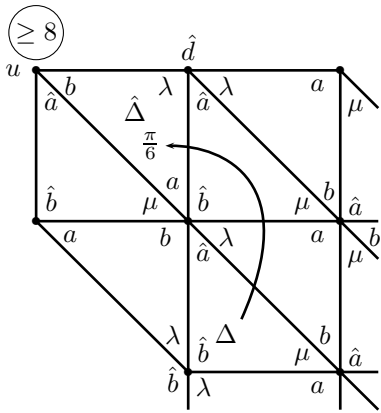
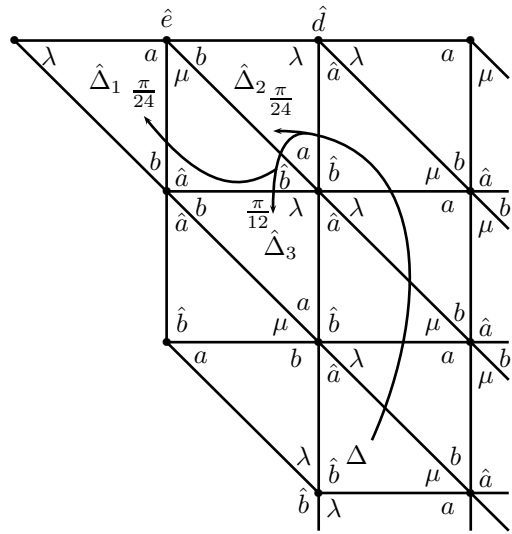
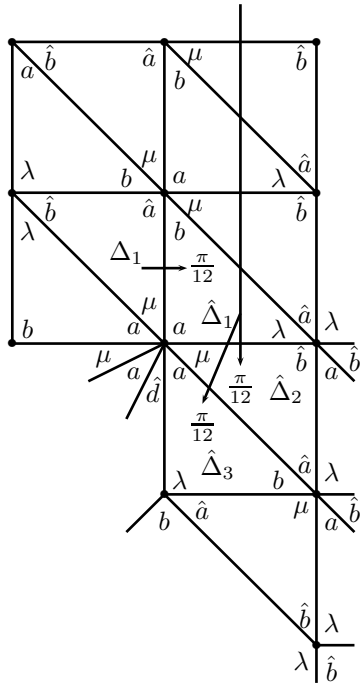
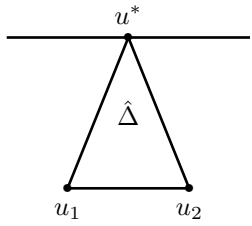
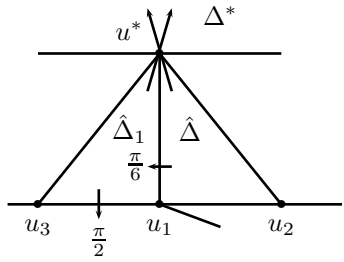
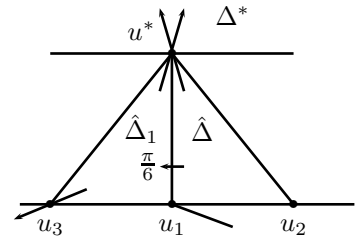
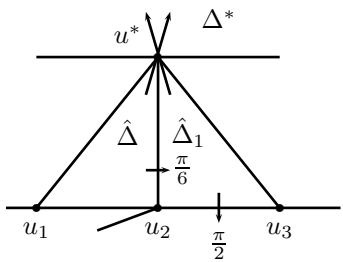
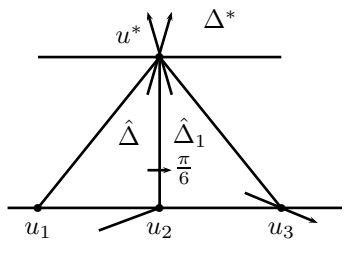
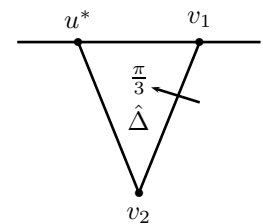
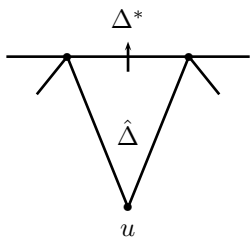
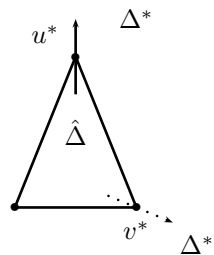
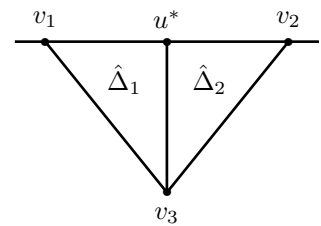
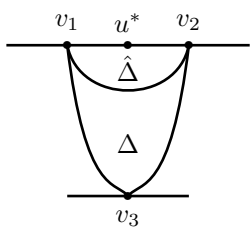
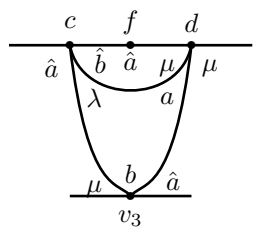
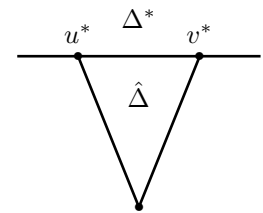
\end{document}